%% file: main.tex
\documentclass{article}

\usepackage{microtype}
\usepackage[nottoc]{tocbibind}
\usepackage[utf8]{inputenc}
\usepackage{graphicx}
\usepackage{float}
\usepackage{bm}
\usepackage{amsthm}
\usepackage{amsmath}
\usepackage{amssymb}
\usepackage{mathrsfs}
\usepackage{dsfont}
\usepackage[all]{xy}
\usepackage{tikz-cd}
\usepackage{stmaryrd}
\usepackage{enumitem}
\usepackage{tabularx}
\usepackage{placeins}
\usepackage{pdflscape}
\usepackage{adjustbox}
\usepackage{rotating}
\usepackage{hyperref}

\theoremstyle{definition}
\newtheorem{Definition}{Definition}[subsection]

\theoremstyle{plain}
\newtheorem{Theorem}[Definition]{Theorem}

\theoremstyle{plain}

\theoremstyle{plain}
\newtheorem{Proposition}[Definition]{Proposition}

\theoremstyle{plain}
\newtheorem{Lemma}[Definition]{Lemma}

\theoremstyle{plain}
\newtheorem{Corollary}[Definition]{Corollary}

\theoremstyle{definition}
\newtheorem{Example}[Definition]{Example}

\theoremstyle{definition}
\newtheorem{Notation}[Definition]{Notation}

\theoremstyle{remark}
\newtheorem{Remark}[Definition]{Remark}

\theoremstyle{plain}
\newcommand{\thistheoremname}{}
\newtheorem*{genericthm*}{\thistheoremname}
\newenvironment{namedthm*}[1]
  {\renewcommand{\thistheoremname}{#1}%
   \begin{genericthm*}}
  {\end{genericthm*}}

\title{The Morita Theory of Fusion 2-Categories}
\author{Thibault D. Décoppet}
\date{May 2023}

\begin{document}

\bibliographystyle{alpha}

     \maketitle
    \hspace{1cm}
    \begin{abstract}
         We develop the Morita theory of fusion 2-categories. In order to do so, we begin by proving that the relative tensor product of modules over a separable algebra in a fusion 2-category exists. We use this result to construct the Morita 3-category of separable algebras in a fusion 2-category. Then, we go on to explain how module 2-categories form a 3-category. After that, we define separable module 2-categories over a fusion 2-category, and prove that the Morita 3-category of separable algebras is equivalent to the 3-category of separable module 2-categories. As a consequence, we show that the dual tensor 2-category with respect to a separable module 2-category, that is the associated 2-category of module 2-endofunctors, is a multifusion 2-category. Finally, we give three equivalent characterizations of Morita equivalence between fusion 2-categories.
    \end{abstract}
    
\tableofcontents

\section*{Introduction}
\addcontentsline{toc}{section}{Introduction}

In the theory of fusion 1-categories, the notion of Morita equivalence plays an essential role. For instance, it is used in the study of the Drinfel'd center (see \cite{Mu2} and \cite{ENO1}), group-graded extensions (see \cite{ENO2}), and group-theoretical fusion 1-categories (see \cite{ENO3}). The first definition of Morita equivalence between fusion 1-categories, which uses Frobenius algebras, was introduced in \cite{FRS}, and \cite{Mu1}, where it was used to study conformal field theories, and subfactors, respectively. In \cite{EO}, the authors subsequently gave an equivalent characterization of Morita equivalence using the concept of the dual tensor 1-category with respect to a module 1-category from \cite{O2} (see also \cite{O1}). Categorifying the original definition of Morita equivalence between algebras introduced in \cite{Mor}, it also natural to say that two fusion 1-categories are Morita equivalent if the associated 2-categories of module 1-categories are equivalent. This last notion of Morita equivalence was proven to be equivalent to the previous ones in \cite{ENO2}.

In the present article, we study the Morita theory of fusion 2-categories. In particular, we categorify the concept of Morita equivalence between fusion 1-categories recalled above. Thus, we expect that this notion will play an analogous role in the theory of fusion 2-categories. In fact, so as to define Morita equivalence between fusion 2-categories, we need to thoroughly examine the properties separable algebras, separable module 2-categories, and dual tensor 2-categories. Such investigations have applications in other contexts. Firstly, it was conjectured in \cite{DR} that fusion 2-categories are the objects of a symmetric monoidal 4-category with duals (in the sense of \cite{L}). Proving this conjecture undoubtedly requires a detailed understanding of this 4-category. But, by analogy with the decategorified setting studied in \cite{DSPS13}, the 1-morphisms of the aforementioned 4-category are separable bimodule 2-categories. Secondly, at the moment, the only construction available to produce new fusion 2-categories out of the ones that are already known is the 2-Deligne tensor product introduced in \cite{D3}. Taking the dual tensor 2-category to a fusion 2-category with respect to a separable module 2-category provides a new method to build interesting fusion 2-categories. Thirdly, using separable module 2-categories, one can define Morita equivalence between separable algebras in an arbitrary fusion 2-category. This is an internalization of the concept of Morita equivalence between fusion 1-categories. Namely, separable algebras in the fusion 2-category $\mathbf{2Vect}$ of finite semisimple 1-categories are exactly multifusion 1-categories, and the associated notion of Morita equivalence coincides with the classical one. Further, it was shown in \cite{D7} that many familiar objects in the theory of fusion 1-categories such as $G$-graded fusion 1-categories over a finite group $G$ are separable algebras in certain fusion 2-categories. Thus, as a byproduct of our investigations, we recover the equivariant Morita theory of $G$-graded fusion 1-categories introduced in \cite{GJS}, and also obtain the correct version of Morita equivalence for fusion 1-categories with a $G$-action.

Let us now recall the equivalent characterizations of Morita equivalence between multifusion 1-categories in detail. Let $\mathcal{C}$ and $\mathcal{D}$ be two multifusion 1-categories over an algebraically closed field of characteristic zero. Further, let us write $\mathbf{Mod}(\mathcal{C})$ for the 2-category of finite semisimple right $\mathcal{C}$-module 1-categories, and likewise for $\mathbf{Mod}(\mathcal{D})$. Categorifying the classical notion of Morita equivalence for algebras, we say that $\mathcal{C}$ and $\mathcal{D}$ are Morita equivalent if the (linear) 2-categories $\mathbf{Mod}(\mathcal{C})$ and $\mathbf{Mod}(\mathcal{D})$ are equivalent. Alternatively, given $\mathcal{M}$ a finite semisimple left $\mathcal{C}$-module 1-category, we can consider $End_{\mathcal{C}}(\mathcal{M})$, the multifusion 1-category of left $\mathcal{C}$-module endofunctor of $\mathcal{M}$. Following \cite{EO}, we use $\mathcal{C}^*_{\mathcal{M}}$ to denote $ End_{\mathcal{C}}(\mathcal{M})$, and call it the dual tensor 1-category to $\mathcal{C}$ with respect to $\mathcal{M}$. Then, we say that $\mathcal{C}$ and $\mathcal{D}$ are Morita equivalent if there exists a faithful finite semisimple left $\mathcal{C}$-module 1-category $\mathcal{C}$ together with a monoidal equivalence between $\mathcal{C}^*_{\mathcal{M}}$ and $\mathcal{D}^{mop}$, that is $\mathcal{D}$ equipped with the opposite monoidal structure. It follows from \cite{ENO2} that this coincides with the notion of Morita equivalence recalled above. Moreover, it follows from \cite{O1} that there exists an algebra $A$ in $\mathcal{C}$ such that $\mathcal{M}$ is equivalent to $Mod_{\mathcal{C}}(A)$, the 1-category of right $A$-modules in $\mathcal{C}$. This implies that there is a monoidal equivalence between $End_{\mathcal{C}}(\mathcal{M})$ and $Bimod_{\mathcal{C}}(A)^{mop}$, the monoidal 1-category of $A$-$A$-bimodules in $\mathcal{C}$. Let us also note that, by \cite{ENO1}, the algebra $A$ is necessarily separable, i.e. $A$ is a special Frobenius algebra. It then follows that $\mathcal{C}$ and $\mathcal{D}$ are Morita equivalent if and only if there exists a faithful separable algebra $A$ in $\mathcal{C}$ together with an equivalence $\mathcal{D}\simeq Bimod_{\mathcal{C}}(A)$ of monoidal 1-categories. This recovers the notion of Morita equivalence introduced in \cite{FRS} and \cite{Mu1}. Let us also remark that, over an arbitrary field, the above discussion remains sensible provided that all the module 1-categories under consideration are assumed to be separable in the sense of \cite{DSPS13}, that is are equivalent to $Mod_{\mathcal{C}}(A)$ for some separable algebra $A$ in $\mathcal{C}$. In fact, they show that a finite semisimple left $\mathcal{C}$-module 1-category $\mathcal{M}$ is separable if and only if $\mathcal{C}^*_{\mathcal{M}}$ is a finite semisimple 1-category.

Our objective is to categorify the equivalent characterizations of Morita equivalence between fusion 1-categories given in the previous paragraph. More precisely, working momentarily over an algebraically closed field of characteristic zero, recall from \cite{DR} that a multifusion 2-category is a finite semisimple rigid monoidal 2-category, and that a fusion 2-category is a multifusion 2-category whose monoidal unit is a simple object. We fix a multifusion 2-category $\mathfrak{C}$ together with an algebra $A$ in $\mathfrak{C}$ that is separable, which implies that the 2-category $\mathbf{Bimod}_{\mathfrak{C}}(A)$ of $A$-$A$-bimodules in $\mathfrak{C}$ is finite semisimple (see \cite{D7}). We now wish to endow $\mathbf{Bimod}_{\mathfrak{C}}(A)$ with a monoidal structure. As expected, the desired monoidal structure is given by the relative tensor product of bimodules over the separable algebra $A$, which generalizes the relative tensor product of finite semisimple module 1-categories over a fusion 1-category introduced in \cite{ENO2}. We establish more generally the existence of the relative tensor product of modules over a separable algebra in any monoidal 2-category $\mathfrak{D}$ with monoidal product $\Box$ that is Karoubi complete in the sense of \cite{GJF}.

\begin{namedthm*}{Theorem \ref{thm:tensormodules}}
Let $B$ be a separable algebra in a Karoubi complete monoidal 2-category $\mathfrak{D}$. Then, the relative tensor product of any right $B$-module $M$, and left $B$-module $N$ in $\mathfrak{D}$ exists, and is given by the splitting of a 2-condensation monad on $M\Box N$.
\end{namedthm*}

\noindent In fact, elaborating on the above result, we construct the Morita 3-category $\mathbf{Mor}^{sep}(\mathfrak{D})$ of separable algebras, bimodules, and their morphisms in $\mathfrak{D}$. Related 3-categories have previously been considered in \cite{GJF} and \cite{H}. Further, we expect that the 3-category $\mathbf{Mor}^{sep}(\mathfrak{D})$ is equivalent to the 3-category $Kar(\mathrm{B}\mathfrak{D})$ of 3-condensation monads, condensation bimodules, and their morphisms in $\mathfrak{D}$ considered in \cite{GJF}. 

We then turn our attention towards the 2-category $\mathbf{End}_{\mathfrak{C}}(\mathfrak{M})$ of left $\mathfrak{C}$-module 2-endofunctors on the left $\mathfrak{C}$-module 2-category $\mathfrak{M}$. We show that this 2-category has a canonical monoidal structure given by composition. More generally, for any fixed monoidal 2-category $\mathfrak{D}$, we will construct a 3-category $\mathbf{LMod}(\mathfrak{D})$ of left $\mathfrak{D}$-module 2-categories, left $\mathfrak{D}$-module 2-functors, left $\mathfrak{D}$-module 2-natural transformations, and left $\mathfrak{D}$-module modification, by promoting the 3-category of 2-categories considered in \cite{Gur}. Further, if $\mathfrak{D}$ is rigid, we will show that if a left $\mathfrak{D}$-module 2-functor has a 2-adjoint as a plain 2-functor, it has a 2-adjoint as a $\mathfrak{D}$-module 2-functor. In particular, for any left $\mathfrak{D}$-module 2-category $\mathfrak{N}$, the monoidal 2-category $\mathbf{End}_{\mathfrak{C}}(\mathfrak{N})$ is rigid if every (plain) 2-endofunctor on $\mathfrak{N}$ has a 2-adjoint.

Now, it was shown in \cite{D4} that the 2-category $\mathbf{Mod}_{\mathfrak{C}}(A)$ of right $A$-module in $\mathfrak{C}$ admits a canonical left $\mathfrak{C}$-module structure. By analogy with the decategorified setting, we wish to compare the monoidal 2-categories $\mathbf{Bimod}_{\mathfrak{C}}(A)$ and $\mathbf{End}_{\mathfrak{C}}(\mathbf{Mod}_{\mathfrak{C}}(A))$. We will do so in more generality by working over an arbitrary field, and letting $\mathfrak{C}$ be a compact semisimple tensor 2-category in the sense of \cite{D5}. Over an algebraically closed field of characteristic zero, this recovers precisely the notion of a multifusion 2-category recalled above. Under these hypotheses, we say that a left $\mathfrak{C}$-module 2-category $\mathfrak{M}$ is separable if it is equivalent as a left $\mathfrak{C}$-module 2-category to $\mathbf{Mod}_{\mathfrak{C}}(A)$ for some separable algebra $A$ in $\mathfrak{C}$. In addition, we write $\mathbf{LMod}^{sep}(\mathfrak{C})$ for the full sub-3-category of $\mathbf{LMod}(\mathfrak{C})$ on the separable module 2-categories. We then prove the following twice categorified version of the classical Eilenberg-Watts theorem.

\begin{namedthm*}{Theorem \ref{thm:equivalencealgebrasmodules2categories}}
Let $\mathfrak{C}$ be a compact semisimple tensor 2-category. There is a linear 3-functor, contravariant on 1-morphisms, $$\mathbf{Mod}_{\mathfrak{C}}:\mathbf{Mor}^{sep}(\mathfrak{C})\rightarrow\mathbf{LMod}^{sep}(\mathfrak{C})$$ that sends a separable algebra in $\mathfrak{C}$ to the associated separable left $\mathfrak{C}$-module 2-category of right modules. Moreover, this 3-functor is an equivalence.
\end{namedthm*}

\noindent The above theorem is an internalization of corollary 3.1.5 of \cite{D5} stating that the Morita 3-category of separable multifusion 1-categories is equivalent to the 3-category of locally separable compact semisimple 2-categories (see also theorem 3.2.2 of \cite{D1}). Namely, this corollary is recovered by taking $\mathfrak{C}=\mathbf{Vect}$ over a perfect field, in which case $\mathbf{Mor}^{sep}(\mathbf{2Vect})$ is the underlying 3-category of the symmetric monoidal 3-category $\mathbf{TC}^{sep}$ of separable tensor 1-categories considered in \cite{DSPS13}. In addition, let us mention that the finite semisimple case of theorem 4.16 of \cite{GJS} is recovered as a consequence of the above theorem for $\mathfrak{C}=\mathbf{2Vect}_G$, the fusion 2-category of 2-vector spaces graded by the finite group $G$, over an algebraically closed field of characteristic zero. 

Let us now assume that the compact semisimple tensor 2-category $\mathfrak{C}$ is locally separable, a mild technical condition, which is always satisfied over an algebraically closed field of characteristic zero. Further, for any monoidal 2-category $\mathfrak{D}$, let us use $\mathfrak{D}^{mop}$ to denote $\mathfrak{D}$ equipped with the opposite monoidal structure. Then, bringing together the various results of this article, we obtain the following theorem.

\begin{namedthm*}{Theorem \ref{thm:bimodrigidmonoidal}}
Let $\mathds{k}$ be a perfect field, and $A$ a separable algebra in a locally separable compact semisimple tensor 2-category $\mathfrak{C}$. Then, $$\mathbf{End}_{\mathfrak{C}}(\mathbf{Mod}_{\mathfrak{C}}(A))\simeq \mathbf{Bimod}_{\mathfrak{C}}(A)^{mop}$$ is a compact semisimple tensor 2-category.
\end{namedthm*}

\noindent In particular, given a separable module 2-category $\mathfrak{M}$, we call $\mathbf{End}_{\mathfrak{C}}(\mathfrak{M})$ the dual tensor 2-category to $\mathfrak{C}$ with respect to $\mathfrak{M}$, which we denote by $\mathfrak{C}^*_{\mathfrak{M}}$. Finally, using the main result of \cite{D4}, we obtain three equivalent characterizations of Morita equivalence between locally separable compact semisimple tensor 2-categories.

\begin{namedthm*}{Theorem \ref{thm:Moritaequivalence}}
For any two locally separable compact semisimple tensor 2-categories $\mathfrak{C}$ and $\mathfrak{D}$ over a perfect field $\mathds{k}$, the following are equivalent:
\begin{enumerate}
    \item The 3-categories $\mathbf{LMod}^{sep}(\mathfrak{C})$ and $\mathbf{LMod}^{sep}(\mathfrak{D})$ are equivalent.
    \item There exists a faithful separable left $\mathfrak{C}$-module 2-category $\mathfrak{M}$, and an equivalence of monoidal 2-categories $\mathfrak{D}^{mp}\simeq \mathfrak{C}_{\mathfrak{M}}^*$.
    \item There exists a faithful separable algebra $A$ in $\mathfrak{C}$, and an equivalence of monoidal 2-categories $\mathfrak{D}\simeq \mathbf{Bimod}_{\mathfrak{C}}(A)$.
\end{enumerate}
If either of the above conditions is satisfied, we say that $\mathfrak{C}$ and $\mathfrak{D}$ are Morita equivalent.
\end{namedthm*}

\noindent We end by examining some examples. Over an algebraically closed field of characteristic zero, we show that, for any finite group $G$, the fusion 2-category $\mathbf{2Vect}_G$ is Morita equivalent to $\mathbf{2Rep}(G)$, the fusion 2-category of 2-representations of $G$. Additionally, we explain how the concept of Morita equivalence between fusion 2-categories recovers the notion of Witt equivalence between non-degenerate braided fusion 1-categories considered in \cite{DMNO}.

\subsection*{Acknowledgments}

I am indebted to Christopher Douglas for his valuable feedback on the contents of the present article. In addition, I would like to thank the referee for their judicious suggestions. 

\input{Preliminaries}

\input{AlgebrasModules}

\input{TensorProduct}

\input{Module2Categories}

\input{SeparableM2C}

\appendix

\input{Appendix}

\bibliography{bibliography.bib}

\end{document}

%% file: Preliminaries.tex
\section{Preliminaries}

\subsection{Graphical Conventions}\label{sub:calculus}

The main objects of study of the present article are (weak) 2-categories with additional structures. In this context, it is convenient to use the graphical calculus originally developed in \cite{GS}, and subsequently modified in \cite{D4}. More precisely, we use string diagrams, in which regions correspond to objects, strings to 1-morphisms, and coupons to 2-morphisms. Our diagrams are to be read from top to bottom, which gives the composition of 1-morphisms, and from left to right, which gives the composition of 2-morphisms. We use the symbol $1$ to denote the identity 1-morphism on an object, but will omit it from the notations if it is not necessary. To illustrate our conventions, let $\mathfrak{C}$ be a 2-category, and let $f:A\rightarrow B$, and $g,h:B\rightarrow C$ be 1-morphisms. Given a 2-morphism $\xi:g\Rightarrow h$, the composite 2-morphism $\xi \circ f$ is represented in our graphical calculus by the following diagram:

$$\includegraphics[width=30mm]{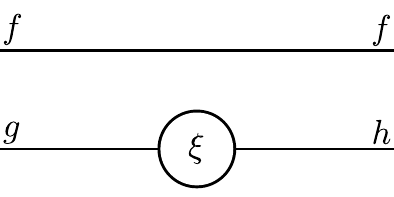}.$$

Throughout, we will work with a monoidal 2-category $\mathfrak{C}$ in the sense of \cite{SP}. In particular, we write $\Box:\mathfrak{C}\times\mathfrak{C}\rightarrow \mathfrak{C}$ for the monoidal product of $\mathfrak{C}$ and $I$ for its monoidal unit. In fact, we will most often assume that $\mathfrak{C}$ is strict cubical, i.e. satisfies definition 2.26 of \cite{SP}, which is not a loss of generality thanks to \cite{Gur}. More precisely, a strict cubical monoidal 2-category is a strict 2-category $\mathfrak{C}$, such that the monoidal product $\Box$ is strictly associative and the unit $I$ is strict. We will therefore systematically omit $I$ from the notations in this case. In addition, the 2-functor $\Box$ is strict in either variable separately. In general, the 2-functor $\Box$ is not strict though. In detail, given pairs of composable 1-morphisms $f_1, f_2$ and $g_1, g_2$ in $\mathfrak{C}$, the 2-isomorphism $$\phi^{\Box}_{(f_2,g_2), (f_1,g_1)}:(f_2\Box g_2)\circ (f_1\Box g_1)\cong (f_2\circ f_1)\Box (g_2\circ g_1)$$ witnessing that $\Box$ preserves the composition of 1-morphisms, called the interchanger, is not trivial. Nevertheless, the strict cubical hypothesis guarantees that $\phi^{\Box}_{(f_2,g_2), (f_1,g_1)}$ is trivial when either $f_2=1$ or $g_1=1$. Given $f$ and $g$ two 1-morphisms in $\mathfrak{C}$, the 2-isomorphism $$\phi^{\Box}_{(f,1),(1,g)}:(f\Box 1)\circ (1\Box g)\cong (1\Box g)\circ (f\Box 1)$$ will be depicted using the diagram below on the left, and its inverse using the diagram on the right: $$\begin{tabular}{c c c c}
\includegraphics[width=20mm]{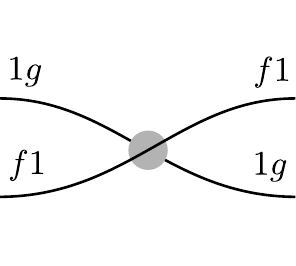},\ \ \ \   & \ \ \ \  \includegraphics[width=20mm]{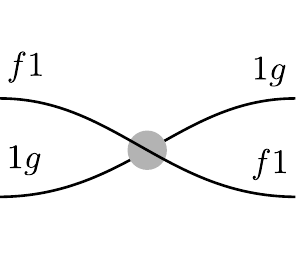}.
\end{tabular}$$

\newlength{\prelim}
\settoheight{\prelim}{\includegraphics[width=30mm]{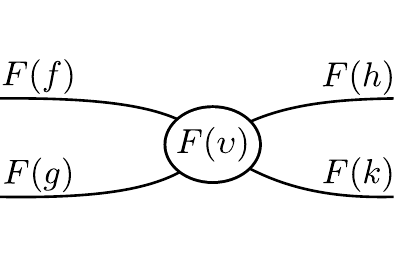}}

\noindent In particular, note that we have omitted the symbol $\Box$. We will systemically do so in order to improve the readability of our diagrams.

In section \ref{sec:M2C}, we will also consider 2-functors and 2-natural transformations, and we now recall from \cite{D4} how to extend the above graphical calculus to these objects. In detail, given $F:\mathfrak{A}\rightarrow \mathfrak{B}$ a (weak) 2-functor, we write $\phi^F_A:Id_{F(A)}\cong F(Id_A)$ for the 2-isomorphism witnessing that $F$ preserves the identity 1-morphism on the object $A$ in $\mathfrak{A}$, and $\phi^F_{g,f}:F(g)\circ F(f)\cong F(g\circ f)$ for the 2-isomorphism witnessing that $F$ preserves the composition of the two composable 1-morphisms $g$ and $f$ in $\mathfrak{A}$. These 2-isomorphisms satisfy well-known compatibility conditions. Now, given any 2-morphism $\upsilon: g\circ f \Rightarrow k\circ h$ in $\mathfrak{A}$, we set: $$\includegraphics[width=30mm]{Pictures/prelim/calculus/Fupsilon.pdf} \raisebox{0.43\prelim}{$\ := (\phi^F_{k,h})^{-1} \cdot F(\upsilon)\cdot \phi^F_{g,f}.$}$$ We extend this convention in the obvious way to the image under a 2-functor of a general 2-morphism, and note that it is well-defined thanks to the coherence axioms for a 2-functor.

Now, let $F,G:\mathfrak{A}\rightarrow \mathfrak{B}$ be two 2-functors, and let $\tau:F\Rightarrow G$ be 2-natural transformation. This means that, for every object $A$ in $\mathfrak{A}$, we have a 1-morphism $\tau_A:F(A)\rightarrow G(A)$, and for every 1-morphism $f:A\rightarrow B$ in $\mathfrak{A}$, we have a 2-isomorphism $$\begin{tikzcd}[sep=tiny]
F(A) \arrow[ddd, "F(f)"']\arrow[rrr, "\tau_A"]  &                                        &    & F(B) \arrow[ddd, "G(f)"]  \\
 &  &    & \\
  &  &  &  \\
F(B)\arrow[rrr, "\tau_B"']\arrow[Rightarrow, rrruuu, "\tau_f", shorten > = 2ex, shorten < = 2ex]                                            &                                        &    &  G(B), 
\end{tikzcd}$$ The collection of these 2-isomorphisms has to satisfy the obvious coherence relations. In our graphical language, we will depict the 2-isomorphism $\tau_f$ using the following diagram on the left, and its inverse using the diagram on the right: $$\begin{tabular}{c c c c}
\includegraphics[width=20mm]{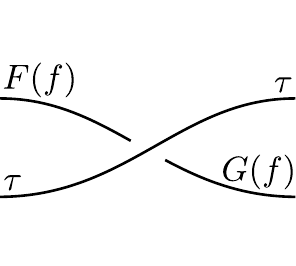},\ \ \ \   & \ \ \ \  \includegraphics[width=20mm]{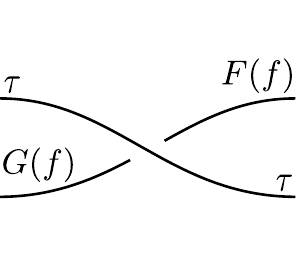}.
\end{tabular}$$

\subsection{2-Condensation Monads}\label{sub:2condensations}

We review the definition of a 2-condensation monad introduced in \cite{GJF} as a categorification of the notions of an idempotent (see also \cite{DR}). More precisely, we recall the unpacked version of this definition given in section 1.1 of \cite{D1}.

\begin{Definition}
A 2-condensation monad in a 2-category $\mathfrak{C}$ is an object $A$ of $\mathfrak{C}$ equipped with a 1-morphism $e:A\rightarrow A$ and two 2-morphisms $\mu:e\circ e\Rightarrow e$ and $\delta:e\Rightarrow e\circ e$ such that $\mu$ is associative, $\delta$ is coassociative, the Frobenius relations hold (i.e. $\delta$ is a 2-morphism of $e$-$e$-bimodules) and $\mu\cdot \delta=Id_e$.
\end{Definition}

\noindent Categorifying the notion of split surjection, $\cite{GJF}$ gave the definition of a 2-condensation, which we recall below. Further, we also review the definition of the splitting of a 2-condensation monad by a 2-condensation, which is spelled out in \cite{D1}.

\begin{Definition}
A 2-condensation in a 2-category $\mathfrak{C}$ is a pair of objects $A,B$ in $\mathfrak{C}$ together with two 1-morphisms $f:A\rightarrow B$ and $g:B\rightarrow A$ and two 2-morphisms $\phi:f\circ g\Rightarrow Id_B$ and $\gamma:Id_B\Rightarrow f\circ g$ such that $\phi\cdot\gamma=Id_{Id_B}$.
\end{Definition}

\begin{Definition}
Let $\mathfrak{C}$ be a 2-category, and $(A,e,\mu^e,\delta^e)$ a 2-condensation monad in $\mathfrak{C}$. A splitting of $(A,e,\mu^e,\delta^e)$ is a 2-condensation $(A,B,f,g,\phi,\gamma)$ together with a 2-isomorphism $\theta:g\circ f\cong e$ such that $$\mu^e= \theta\cdot (g\circ \phi\circ f)\cdot (\theta^{-1}\circ \theta^{-1})\textrm{ and }\delta^e= (\theta\circ\theta)\cdot (g\circ \gamma\circ f)\cdot \theta^{-1}.$$
\end{Definition}

\begin{Remark}
Let $\mathfrak{C}$ be a 2-category whose $Hom$-categories are idempotent complete. It was shown in theorem 2.3.2 of \cite{GJF} that the 2-category of splittings of a fixed 2-condensation monads in $\mathfrak{C}$ is either empty or a contractible 2-groupoid.
\end{Remark}

Following \cite{GJF}, we will call a 2-category locally idempotent complete if its $Hom$-categories are idempotent complete, that is idempotents splits. Further, when working over a fixed field $\mathds{k}$, we will call a $\mathds{k}$-linear 2-category locally Cauchy complete if its $Hom$-categories are Cauchy complete, that is they have direct sums and idempotents splits.

\begin{Definition}
A locally idempotent complete 2-category is Karoubi complete if every 2-condensation monad splits. A locally Cauchy complete $\mathds{k}$-linear 2-category is Cauchy complete if it is Karoubi complete and has direct sums for objects.
\end{Definition}

\begin{Remark}
It is always possible to Karoubi complete an arbitrary locally idempotent complete 2-category (see \cite{DR} and \cite{GJF}). Further, this process satisfies a precise 3-universal property as explained in \cite{D1}.
\end{Remark}

\subsection{Compact Semisimple 2-Categories}\label{sub:SS2CF2C}

Let $\mathds{k}$ be a field. We now review the definition of a semisimple 2-category, given in \cite{DR} over algebraically closed field of characteristic zero. We will then recall the notion of a compact semisimple 2-category introduced in \cite{D5}.

\begin{Definition}
A $\mathds{k}$-linear 2-category is semisimple if it is locally semisimple, has right and left adjoints for 1-morphisms, and is Cauchy complete.
\end{Definition}

An object $C$ of a semisimple 2-category $\mathfrak{C}$ is called simple if the identity 1-morphism $Id_C$ is a simple object of the 1-category $End_{\mathfrak{C}}(C)$. We say that two simple object $C$, $D$ of $\mathfrak{C}$ are in the same connected component if there exists a non-zero 1-morphism between them. As explained in section 1 of \cite{D5}, this defines an equivalence relation on the set of simple object, whose equivalence classes are called the connected components of $\mathfrak{C}$.

\begin{Definition}
A semisimple $\mathds{k}$-linear 2-category is compact if it is locally finite semisimple and has finitely many connected component.
\end{Definition}

As was shown in \cite{D5}, the notion of compact semisimple 2-category is the appropriate categorification of the definition of a finite semisimple 1-category. Namely, following \cite{DR}, a finite semisimple 2-category is a semisimple 2-categories which is locally finite semisimple and has finitely many equivalence classes of simple objects. However, it was proven in \cite{D5} that, over a general field, there does not exist any finite semisimple 2-category, but there always exists compact semisimple 2-categories. Let us note that, over algebraically closed fields or real closed fields, they do show that every compact semisimple 2-category is in fact finite.

Finally, we recall the definitions of a tensor 2-category and of a fusion 2-category, as introduced in \cite{DR} over algebraically closed fields of characteristic zero. We proceed to give some examples.

\begin{Definition}
A tensor 2-category is a rigid monoidal $\mathds{k}$-linear 2-category. A fusion 2-category is finite semisimple tensor 2-category, whose monoidal unit is simple.
\end{Definition}

\begin{Example}
A perfect ($\mathds{k}$-linear) 1-category is a finite semisimple ($\mathds{k}$-linear) 1-category, for which the algebra of endomorphisms of any object is separable. Note that if $\mathds{k}$ is algebraically closed or has characteristic zero, then every finite semisimple 1-category is perfect. We write $\mathbf{2Vect}$ for the 2-category of perfect finite semisimple 1-categories, also called perfect 2-vector spaces. The Deligne tensor product endows $\mathbf{2Vect}$ with the structure of a fusion 2-category.
\end{Example}

\begin{Example}\label{ex:twistedgroupgraded2vectorspaces}
Let $G$ be a finite group. We use $\mathbf{2Vect}_G$ to denote the compact semisimple 2-category of $G$-graded perfect 2-vector spaces. The convolution product turns $\mathbf{2Vect}_G$ into a compact semisimple tensor 2-category. Furthermore, given a 4-cocycle $\pi$ for $G$ with coefficients in $\mathds{k}^{\times}$, we can form the fusion 2-category $\mathbf{2Vect}^{\pi}_G$ by twisting the structure 2-isomorphisms of $\mathbf{2Vect}_G$ using $\pi$ (see construction 2.1.16 of \cite{DR} or \cite{Delc}).
\end{Example}

\begin{Example}\label{ex:connectedfusion2categories}
Let us fix $\mathcal{C}$ a finite semisimple tensor 1-category (over $\mathds{k}$). Following \cite{DSPS13}, we say that a finite semisimple right $\mathcal{C}$-module 1-category is separable if it is equivalent to the 1-category of left modules over a separable algebra in $\mathcal{C}$. If $\mathds{k}$ has characteristic zero, every finite semisimple $\mathcal{C}$-module 1-category is separable. We write $\mathbf{Mod}(\mathcal{C})$ for the compact semisimple 2-category of separable right $\mathcal{C}$-module 1-categories. If $\mathcal{B}$ is a braided finite semisimple tensor 1-category with braiding $\beta$, then the relative Deligne tensor product over $\mathcal{B}$ endows the 2-category $\mathbf{Mod}(\mathcal{B})$ with a rigid monoidal structure, so that $\mathbf{Mod}(\mathcal{B})$ is a compact semisimple tensor 2-category (see \cite{D5}).
\end{Example}

\begin{Example}\label{ex:2representations}
Let $G$ be a finite group whose order is coprime to $char(\mathds{k})$. We write $\mathrm{B}G$ for the 2-category with one object $*$, and $End_{\mathrm{B}G}(*)=G$. We may consider the compact semisimple 2-category $\mathbf{Fun}(\mathrm{B}G,\mathbf{2Vect})$ of (finite perfect) 2-representations of $G$, denoted by $\mathbf{2Rep}(G)$. Said differently, the objects of $\mathbf{2Rep}(G)$ are perfect 2-vector spaces equipped with a $G$-action. The symmetric monoidal structure of $\mathbf{2Vect}$ endows $\mathbf{2Rep}(G)$ with the structure of a symmetric compact semisimple 2-category. More precisely, given $V$ and $W$ two 2-vector spaces with a $G$-action, their monoidal product is given by the Deligne tensor product $V\boxtimes W$ endowed with the diagonal $G$-action. The compact semisimple 2-category $\mathbf{2Rep}(G)$ is fact rigid as can be seen either directly or from lemma \ref{lem:2repG} below.
\end{Example}

The next lemma gives an alternative description of the symmetric monoidal 2-category $\mathbf{2Rep}(G)$ of perfect 2-representations of a finite group $G$. To this end, let us write $\mathbf{Rep}(G)$ for the symmetric fusion 1-category of finite dimensional representations of $G$.

\begin{Lemma}\label{lem:2repG}
Let $G$ be a finite group whose order is coprime to $char(\mathds{k})$. The symmetric monoidal compact semisimple 2-categories $\mathbf{Mod}(\mathbf{Rep}(G))$ and $\mathbf{2Rep}(G)$ are equivalent. In particular, $\mathbf{2Rep}(G)$ is rigid.
\end{Lemma}
\begin{proof}
This follows from a slight elaboration on theorem 8.5 of \cite{Gre}. For completeness, we give a proof using the theory of compact semisimple tensor 2-categories. By definition, the monoidal unit $I$ of $\mathbf{2Rep}(G)$ is $\mathbf{Vect}$, the 1-category of finite $\mathds{k}$-vector spaces, equipped with the trivial $G$-action, and inspection shows that $End_{\mathbf{2Rep}(G)}(I)\cong \mathbf{Rep}(G)$ as symmetric finite semisimple tensor 1-categories. Finally, note that the compact semisimple 2-category $\mathbf{2Rep}(G)$ is a connected, so that the desired equivalence of symmetric monoidal compact semisimple 2-categories follows from proposition 3.3.4 of \cite{D5}.
\end{proof}

%% file: AlgebrasModules.tex
\section{Algebras and Modules}

We review some key definitions using our graphical calculus. More precisely, we begin recalling the definition of an algebra in a (strict cubical) monoidal 2-category. We go on to review the definitions of right and left modules as well as that of bimodules. We end this section by recollecting the definitions of rigid and separable algebras, and giving plenty of examples in fusion 2-categories.

\subsection{Algebras}

Throughout, we work with a fixed strict cubical monoidal 2-category $\mathfrak{C}$. We begin by recalling the definition of an algebra (also called pseudo-monoid in \cite{DS}) in $\mathfrak{C}$ in the form of definition 1.2.1 of \cite{D7}. For the definition of an algebra in an arbitrary monoidal 2-category expressed using our graphical language, we refer the reader to definition 3.1.1 of \cite{D4}.

\begin{Definition}\label{def:algebra}
An algebra in $\mathfrak{C}$ consists of:
\begin{enumerate}
    \item An object $A$ of $\mathfrak{C}$;
    \item Two 1-morphisms $m:A\Box A\rightarrow A$ and $i:I\rightarrow A$;
    \item Three 2-isomorphisms
\end{enumerate}
\begin{center}
\begin{tabular}{@{}c c c@{}}
$\begin{tikzcd}[sep=small]
A \arrow[rrrr, equal] \arrow[rrdd, "i1"'] &  & {} \arrow[dd, Rightarrow, "\lambda"', near start, shorten > = 1ex] &  & A \\
                                   &  &                           &  &   \\
                                   &  & AA, \arrow[rruu, "m"']     &  &  
\end{tikzcd}$

&

$\begin{tikzcd}[sep=small]
AAA \arrow[dd, "1m"'] \arrow[rr, "m1"]    &  & AA \arrow[dd, "m"] \\
                                            &  &                      \\
AA \arrow[rr, "m"'] \arrow[rruu, Rightarrow, "\mu", shorten > = 2.5ex, shorten < = 2.5ex] &  & A,                   
\end{tikzcd}$

&

$\begin{tikzcd}[sep=small]
                                  &  & AA \arrow[rrdd, "m"] \arrow[dd, Rightarrow, "\rho", shorten > = 1ex, shorten < = 2ex] &  &   \\
                                  &  &                                             &  &   \\
A \arrow[rruu, "1i"] \arrow[rrrr,equal] &  & {}                                          &  & M,
\end{tikzcd}$

\end{tabular}
\end{center}

satisfying:

\begin{enumerate}
\item [a.] We have
\end{enumerate}

\settoheight{\prelim}{\includegraphics[width=52.5mm]{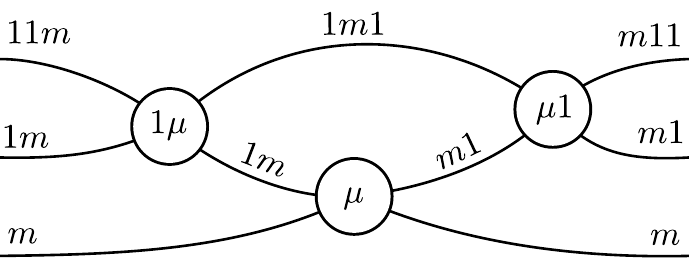}}

\begin{equation}\label{eqn:algebraassociativity}
\begin{tabular}{@{}ccc@{}}

\includegraphics[width=52.5mm]{Pictures/prelim/algebra/associativity1.pdf} & \raisebox{0.45\prelim}{$=$} &
\includegraphics[width=40mm]{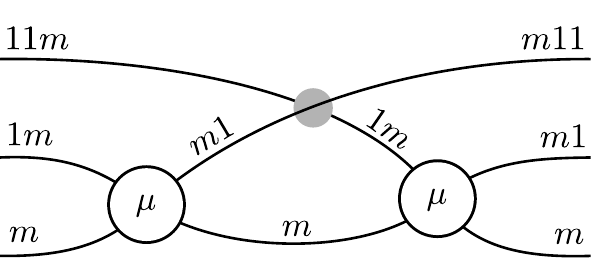},

\end{tabular}
\end{equation}

\begin{enumerate}
\item [b.] We have:
\end{enumerate}

\settoheight{\prelim}{\includegraphics[width=22.5mm]{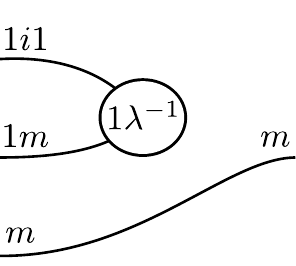}}

\begin{equation}\label{eqn:algebraunitality}
\begin{tabular}{@{}ccc@{}}

\includegraphics[width=22.5mm]{Pictures/prelim/algebra/unitality1.pdf} & \raisebox{0.45\prelim}{$=$} &

\includegraphics[width=37.5mm]{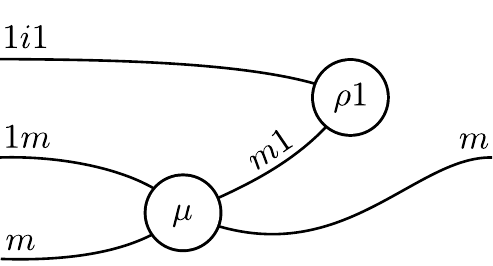}.

\end{tabular}
\end{equation}
\end{Definition}

We will make use of the following coherence results for algebras derived in section 6.3 of \cite{Hou}. We will also use the analogue of equation (\ref{eqn:coherenceleft}) for $\rho$, which follows from lemma \ref{lem:coherenceright} below.

\begin{Lemma}
Given any algebra $A$, the following two equalities hold:

\settoheight{\prelim}{\includegraphics[width=37.5mm]{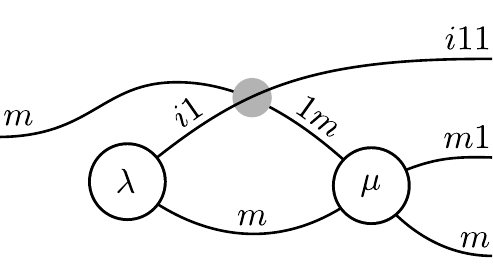}}

\begin{equation}\label{eqn:coherenceleft}
\begin{tabular}{@{}ccc@{}}
\includegraphics[width=37.5mm]{Pictures/prelim/coherenceleft1.pdf}&
\raisebox{0.45\prelim}{$=$} &
\includegraphics[width=22.5mm]{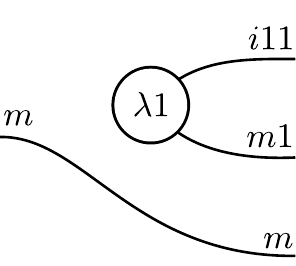},
\end{tabular}
\end{equation}

\settoheight{\prelim}{\includegraphics[width=30mm]{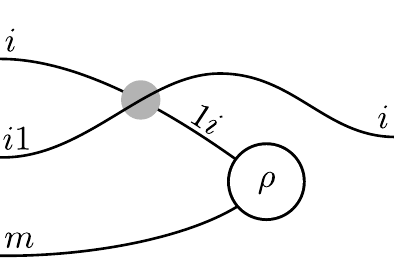}}

\begin{equation}\label{eqn:coherencemiddle}
\begin{tabular}{@{}ccc@{}}
\includegraphics[width=30mm]{Pictures/prelim/coherencemixte1.pdf}&
\raisebox{0.45\prelim}{$=$} &
\includegraphics[width=22.5mm]{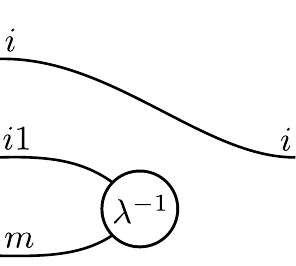}.
\end{tabular}
\end{equation}

\end{Lemma}

\subsection{Modules}

Let us fix an algebra $A$ in the strict cubical monoidal 2-category $\mathfrak{C}$. We now recall the notion of a right $A$-module in $\mathfrak{C}$ given in definition 1.2.3 of \cite{D7}. We invite the reader to consult definition 3.2.1 of \cite{D4} for a version of this definition in a general monoidal 2-category.

\begin{Definition}\label{def:module}
A right $A$-module in $\mathfrak{C}$ consists of:
\begin{enumerate}
    \item An object $M$ of $\mathfrak{C}$;
    \item A 1-morphism $n^M:M\Box A\rightarrow M$;
    \item Two 2-isomorphisms
\end{enumerate}
\begin{center}
\begin{tabular}{@{}c c@{}}
$\begin{tikzcd}[sep=small]
MAA \arrow[dd, "1m"'] \arrow[rr, "n^M1"]    &  & MA \arrow[dd, "n^M"] \\
                                            &  &                      \\
MA \arrow[rr, "n^M"'] \arrow[rruu, Rightarrow, "\nu^M", shorten > = 2.5ex, shorten < = 2.5ex] &  & M,                   
\end{tikzcd}$

&

$\begin{tikzcd}[sep=small]
                                  &  & MA \arrow[rrdd, "n^M"] \arrow[dd, Rightarrow, "\rho^M", shorten > = 1ex, shorten < = 2ex] &  &   \\
                                  &  &                                             &  &   \\
M \arrow[rruu, "1i"] \arrow[rrrr,equal] &  & {}                                          &  & M,
\end{tikzcd}$
\end{tabular}
\end{center}

satisfying:

\begin{enumerate}
\item [a.] We have
\end{enumerate}

\settoheight{\prelim}{\includegraphics[width=52.5mm]{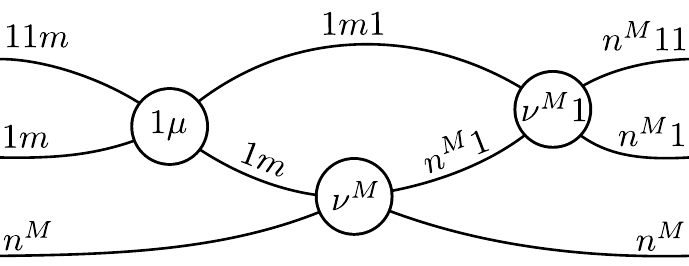}}

\begin{equation}\label{eqn:moduleassociativity}
\begin{tabular}{@{}ccc@{}}

\includegraphics[width=52.5mm]{Pictures/prelim/module/associativity1.pdf} & \raisebox{0.45\prelim}{$=$} &
\includegraphics[width=45mm]{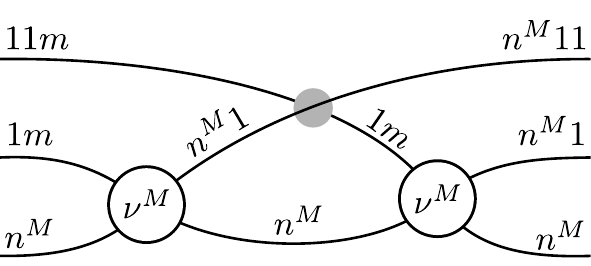},

\end{tabular}
\end{equation}

\begin{enumerate}
\item [b.] We have:
\end{enumerate}

\settoheight{\prelim}{\includegraphics[width=22.5mm]{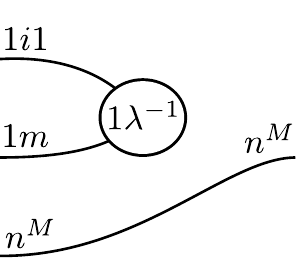}}

\begin{equation}\label{eqn:moduleunitality}
\begin{tabular}{@{}ccc@{}}

\includegraphics[width=22.5mm]{Pictures/prelim/module/unitality1.pdf} & \raisebox{0.45\prelim}{$=$} &

\includegraphics[width=37.5mm]{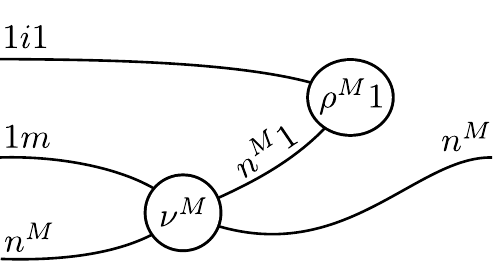}.

\end{tabular}
\end{equation}
\end{Definition}

For later use, let us recall the following coherence result established in lemma 1.2.8 of \cite{D7}.

\begin{Lemma}\label{lem:coherenceright}
Given any right $A$-module $M$, we have the following equality:

\settoheight{\prelim}{\includegraphics[width=37.5mm]{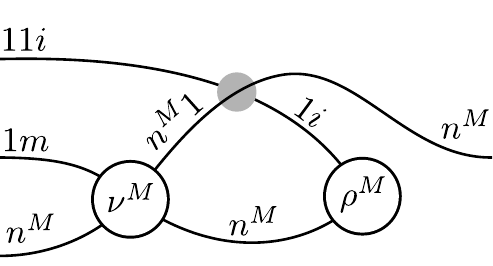}}

\begin{equation}\label{eqn:coherenceright}
\begin{tabular}{@{}ccc@{}}
\includegraphics[width=37.5mm]{Pictures/prelim/module/coherenceright1.pdf}&
\raisebox{0.45\prelim}{$=$} &
\includegraphics[width=22.5mm]{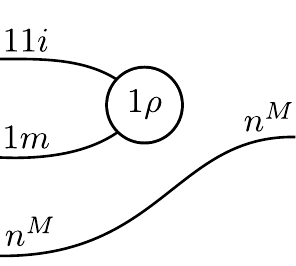}.
\end{tabular}
\end{equation}
\end{Lemma}

Finally, let us recall definitions 3.2.6 and 3.2.7 of \cite{D4}.

\begin{Definition}\label{def:modulemap}
Let $M$ and $N$ be two right $A$-modules. A right $A$-module 1-morphism consists of a 1-morphism $f:M\rightarrow N$ in $\mathfrak{C}$ together with an invertible 2-morphism

$$\begin{tikzcd}[sep=small]
MA \arrow[dd, "f1"'] \arrow[rr, "n^M"]    &  & M \arrow[dd, "f"] \\
                                            &  &                      \\
NA \arrow[rr, "n^N"'] \arrow[rruu, Rightarrow, "\psi^f", shorten > = 2.5ex, shorten < = 2.5ex] &  & N,                   
\end{tikzcd}$$

subject to the coherence relations:

\begin{enumerate}
\item [a.] We have:
\end{enumerate}

\settoheight{\prelim}{\includegraphics[width=52.5mm]{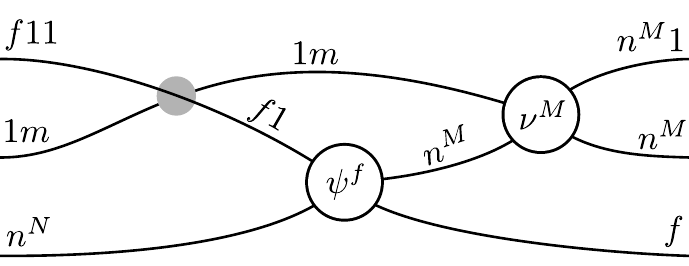}}

\begin{equation}\label{eqn:modulemapassociativity}
\begin{tabular}{@{}ccc@{}}

\includegraphics[width=52.5mm]{Pictures/prelim/module/map1.pdf} & \raisebox{0.45\prelim}{$=$} &

\includegraphics[width=52.5mm]{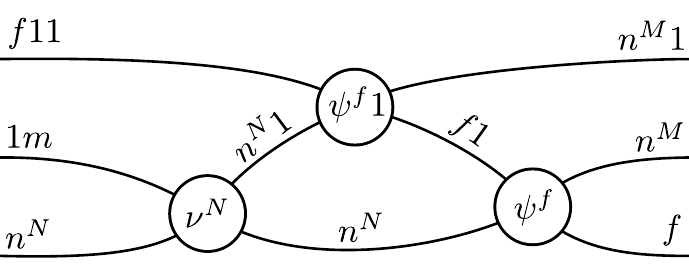},

\end{tabular}
\end{equation}

\begin{enumerate}
\item [b.] We have:
\end{enumerate}

\settoheight{\prelim}{\includegraphics[width=30mm]{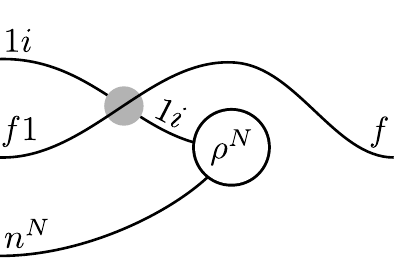}}

\begin{equation}\label{eqn:modulemapunitality}
\begin{tabular}{@{}ccc@{}}

\includegraphics[width=30mm]{Pictures/prelim/module/map3.pdf} & \raisebox{0.45\prelim}{$=$} &

\includegraphics[width=30mm]{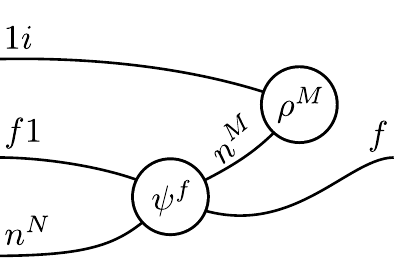}.

\end{tabular}
\end{equation}
\end{Definition}

\begin{Definition}\label{def:moduleintertwiner}
Let $M$ and $N$ be two right $A$-modules, and $f,g:M\rightarrow M$ two right $A$-module 1-morphisms. A right $A$-module 2-morphism $f\Rightarrow g$ is a 2-morphism $\gamma:f\Rightarrow g$ in $\mathfrak{C}$ that satisfies the following equality:

\settoheight{\prelim}{\includegraphics[width=30mm]{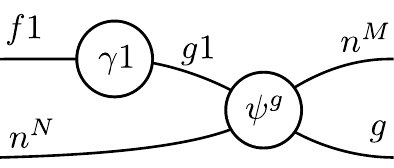}}

\begin{center}
\begin{tabular}{@{}ccc@{}}

\includegraphics[width=30mm]{Pictures/prelim/module/2morphism1.pdf} & \raisebox{0.45\prelim}{$=$} &

\includegraphics[width=30mm]{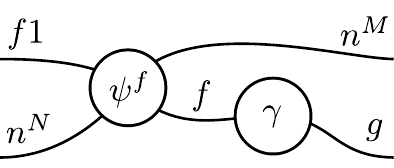}.

\end{tabular}
\end{center}
\end{Definition}

The above structures can be assembled into a 2-category as was proven in lemma 3.2.10 of \cite{D4}. In fact, as we have assumed that $\mathfrak{C}$ is strict cubical, this 2-category is strict.

\begin{Lemma}
Right $A$-modules, right $A$-module 1-morphisms, and right $A$-module 2-morphisms in $\mathfrak{C}$ form a strict 2-category, which we denote by $\mathbf{Mod}_{\mathfrak{C}}(A)$.
\end{Lemma}

Let us now recall the definition of left $A$-module in $\mathfrak{C}$ from definition A.1.1 of \cite{D7}.

\begin{Definition}\label{def:leftmodule}
A left $A$-module in $\mathfrak{C}$ consists of:
\begin{enumerate}
    \item An object $M$ of $\mathfrak{C}$;
    \item A 1-morphism $l^M:A\Box M\rightarrow M$;
    \item Two 2-isomorphisms
\end{enumerate}
\begin{center}
\begin{tabular}{@{}c c@{}}
$\begin{tikzcd}[sep=small]
M \arrow[rrrr, equal] \arrow[rrdd, "i1"'] &  & {} \arrow[dd, Rightarrow, "\lambda^M"', near start, shorten > = 1ex] &  & M \\
                                   &  &                           &  &   \\
                                   &  & AM, \arrow[rruu, "l^M"']     &  &  
\end{tikzcd}$&

$\begin{tikzcd}[sep=small]
AAM \arrow[dd, "1l^M"'] \arrow[rr, "m1"]    &  & AM \arrow[dd, "l^M"] \\
                                            &  &                      \\
AM \arrow[rr, "l^M"'] \arrow[rruu, Rightarrow, "\kappa^M", shorten > = 2.5ex, shorten < = 2.5ex] &  & M,                   
\end{tikzcd}$
\end{tabular}
\end{center}

satisfying:

\begin{enumerate}
\item [a.] We have:
\end{enumerate}

\settoheight{\prelim}{\includegraphics[width=52.5mm]{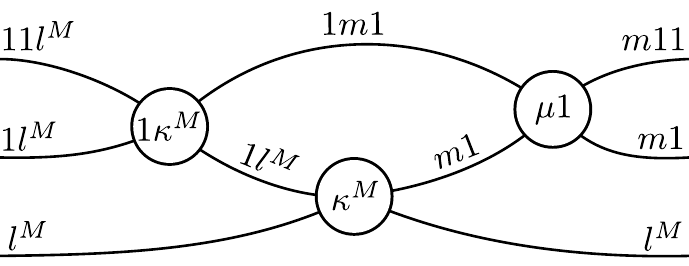}}

\begin{equation}\label{eqn:leftmoduleassociativity}
\begin{tabular}{@{}ccc@{}}

\includegraphics[width=52.5mm]{Pictures/prelim/leftmodule/associativity1.pdf} & \raisebox{0.45\prelim}{$=$} &
\includegraphics[width=45mm]{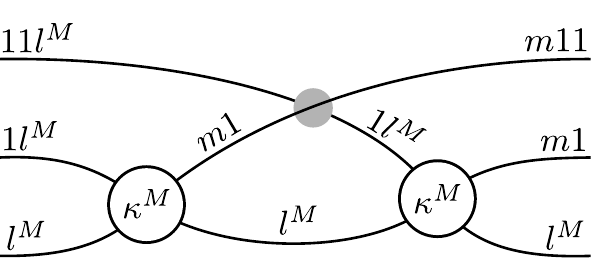},

\end{tabular}
\end{equation}

\begin{enumerate}
\item [b.] We have:
\end{enumerate}

\settoheight{\prelim}{\includegraphics[width=22.5mm]{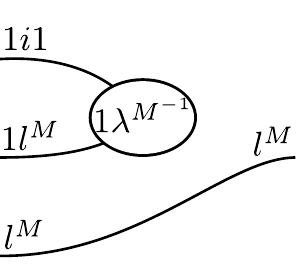}}

\begin{equation}\label{eqn:leftmoduleunitality}
\begin{tabular}{@{}ccc@{}}

\includegraphics[width=22.5mm]{Pictures/prelim/leftmodule/unitality1.pdf} & \raisebox{0.45\prelim}{$=$} &

\includegraphics[width=37.5mm]{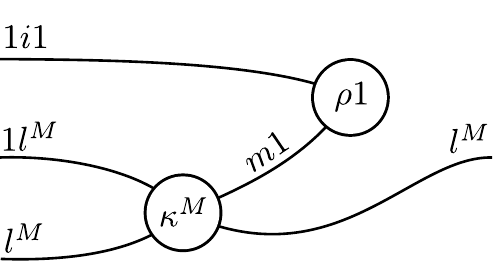}.

\end{tabular}
\end{equation}
\end{Definition}

\begin{Definition}\label{def:leftmodulemap}
Let $M$ and $N$ be two left $A$-modules. A left $A$-module 1-morphism consists of a 1-morphism $f:M\rightarrow N$ in $\mathfrak{C}$ together with an invertible 2-morphism

$$\begin{tikzcd}[sep=small]
AM \arrow[dd, "1f"'] \arrow[rr, "l^M"]    &  & M \arrow[dd, "f"] \\
                                            &  &                      \\
AN \arrow[rr, "l^N"'] \arrow[rruu, Rightarrow, "\chi^f", shorten > = 2.5ex, shorten < = 2.5ex] &  & N,                   
\end{tikzcd}$$

subject to the coherence relations:

\begin{enumerate}
\item [a.] We have:
\end{enumerate}

\settoheight{\prelim}{\includegraphics[width=45mm]{Pictures/prelim/module/map1.pdf}}

\begin{equation}\label{eqn:leftmodulemapassociativity}
\begin{tabular}{@{}ccc@{}}

\includegraphics[width=45mm]{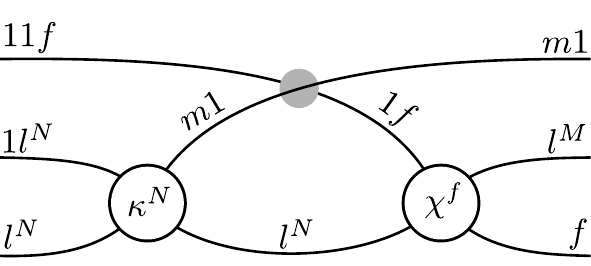} & \raisebox{0.45\prelim}{$=$} &

\includegraphics[width=45mm]{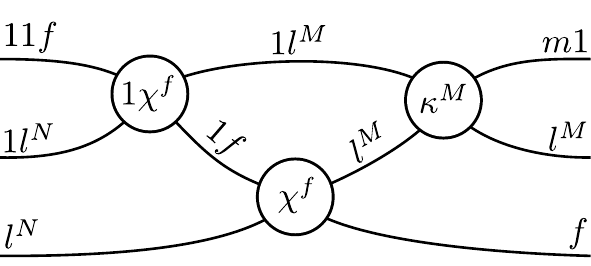},

\end{tabular}
\end{equation}

\begin{enumerate}
\item [b.] We have:
\end{enumerate}

\settoheight{\prelim}{\includegraphics[width=30mm]{Pictures/prelim/module/map3.pdf}}

\begin{equation}\label{eqn:leftmodulemapunitality}
\begin{tabular}{@{}ccc@{}}

\includegraphics[width=30mm]{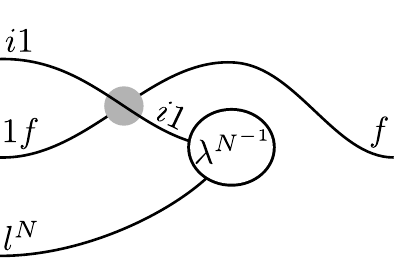} & \raisebox{0.45\prelim}{$=$} &

\includegraphics[width=30mm]{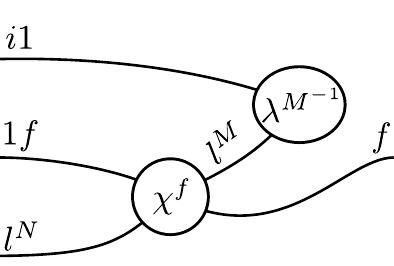}.

\end{tabular}
\end{equation}
\end{Definition}

\begin{Definition}\label{def:leftmoduleintertwiner}
Let $M$ and $N$ be two left $A$-modules, and $f,g:M\rightarrow M$ two left $A$-module 1-morphisms. A left $A$-module 2-morphism $f\Rightarrow g$ is a 2-morphism $\gamma:f\Rightarrow g$ in $\mathfrak{C}$ that satisfies the following equality:

\settoheight{\prelim}{\includegraphics[width=30mm]{Pictures/prelim/module/2morphism1.pdf}}

\begin{center}
\begin{tabular}{@{}ccc@{}}

\includegraphics[width=30mm]{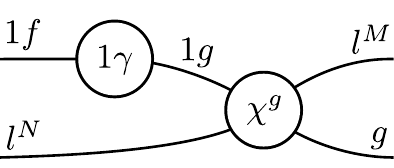} & \raisebox{0.45\prelim}{$=$} &

\includegraphics[width=30mm]{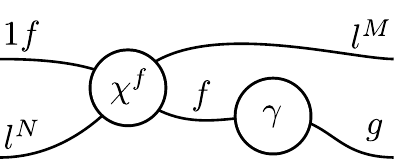}.

\end{tabular}
\end{center}
\end{Definition}

A slight variant of the proof of lemma 3.2.10 of \cite{D4} shows that left $A$-modules and their morphisms can be assembled into a 2-category. As $\mathfrak{C}$ is strict cubical, this 2-category is strict.

\begin{Lemma}
Left $A$-modules, left $A$-module 1-morphisms, and left $A$-module 2-morphisms in $\mathfrak{C}$ form a strict 2-category, which we denote by $\mathbf{LMod}_{\mathfrak{C}}(A)$.
\end{Lemma}

\subsection{Bimodules}

Let $(A,m^A,i^A,\lambda^A,\mu^A,\rho^A)$ and $(B,m^B,i^B,\lambda^B,\mu^B,\rho^B)$ be algebras in the strict cubical monoidal 2-category $\mathfrak{C}$. We now review the notion of an $A$-$B$-bimodule in $\mathfrak{C}$.

\begin{Definition}\label{def:bimodule}
An $A$-$B$-bimodule in $\mathfrak{C}$ consists of:
\begin{enumerate}
    \item An object $P$ of $\mathfrak{C}$;
    \item The data $(P,l^P,\lambda^P,\kappa^P)$ of a left $A$-module structure on $P$;
    \item The data $(P,n^P,\nu^P,\rho^P)$ of a right $B$-module structure on $P$;
    \item A 2-isomorphism
\end{enumerate}

$$\begin{tikzcd}[sep=small]
APB \arrow[dd, "1n^P"'] \arrow[rr, "l^P1"]    &  & PB \arrow[dd, "n^P"] \\
                                            &  &                      \\
AP \arrow[rr, "l^P"'] \arrow[rruu, Rightarrow, "\beta^P", shorten > = 2.5ex, shorten < = 2.5ex] &  & M,                   
\end{tikzcd}$$

satisfying:

\begin{enumerate}
\item [a.] We have:
\end{enumerate}

\settoheight{\prelim}{\includegraphics[width=52.5mm]{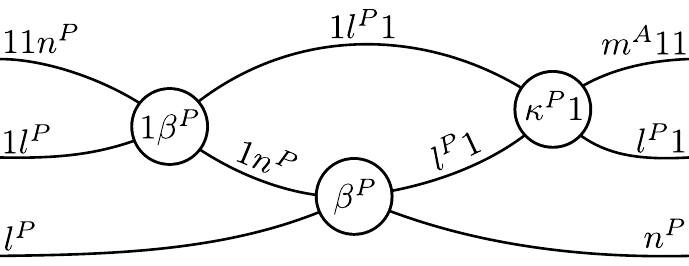}}

\begin{equation}\label{eqn:leftbimoduleassociativity}
\begin{tabular}{@{}ccc@{}}

\includegraphics[width=52.5mm]{Pictures/prelim/leftmodule/bimodule1.pdf} & \raisebox{0.45\prelim}{$=$} &
\includegraphics[width=45mm]{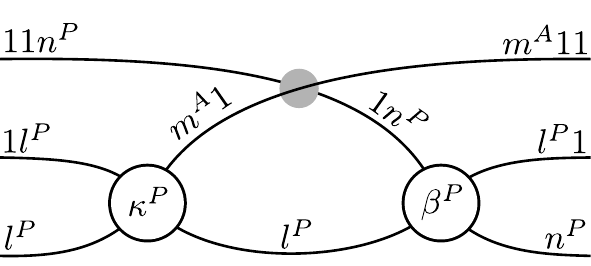},

\end{tabular}
\end{equation}

\begin{enumerate}
\item [b.] We have:
\end{enumerate}

\settoheight{\prelim}{\includegraphics[width=52.5mm]{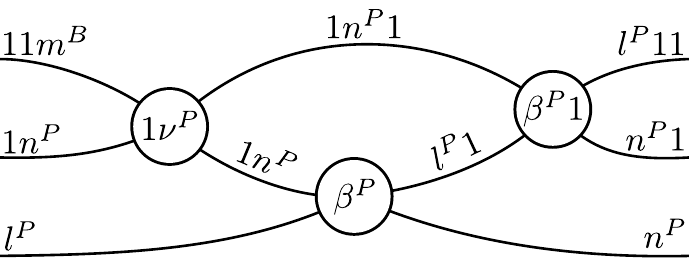}}

\begin{equation}\label{eqn:rightbimoduleassociativity}
\begin{tabular}{@{}ccc@{}}

\includegraphics[width=52.5mm]{Pictures/prelim/leftmodule/bimodule3.pdf} & \raisebox{0.45\prelim}{$=$} &
\includegraphics[width=45mm]{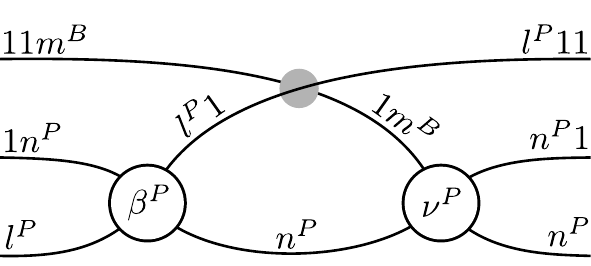}.

\end{tabular}
\end{equation}
\end{Definition}

\begin{Definition}\label{def:bimodulemap}
Let $P$ and $Q$ be two $A$-$B$-bimodules in $\mathfrak{C}$. An $A$-$B$-bimodule 1-morphism consists of a 1-morphism $f:P\rightarrow Q$ in $\mathfrak{C}$ together with the data $(f,\chi^f)$ of a left $A$-module structure and $(f, \psi^f)$ of a right $B$-module structure satisfying:

\settoheight{\prelim}{\includegraphics[width=45mm]{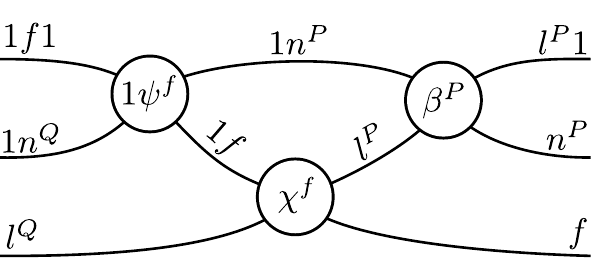}}

\begin{equation}\label{eqn:bimodulemapassociativity}
\begin{tabular}{@{}ccc@{}}

\includegraphics[width=45mm]{Pictures/prelim/leftmodule/bimodulemap1.pdf} & \raisebox{0.45\prelim}{$=$} &

\includegraphics[width=52.5mm]{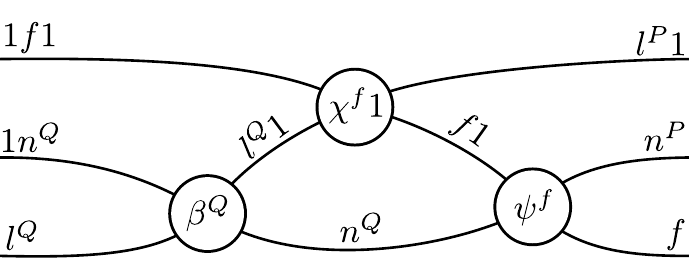}.

\end{tabular}
\end{equation}
\end{Definition}

\begin{Definition}\label{def:bimoduleintertwiner}
Let $P$ and $Q$ be two $A$-$B$-bimodules, and $f,g:P\rightarrow Q$ two $A$-$B$-bimodules 1-morphisms in $\mathfrak{C}$. An $A$-$B$-bimodule 2-morphism $f\Rightarrow g$ is a 2-morphism $\gamma:f\Rightarrow g$ in $\mathfrak{C}$, which is both a left $A$-module 2-morphism and a right $B$-module 2-morphism.
\end{Definition}

A slight elaboration on the proof of lemma 3.2.10 of \cite{D4} proves that $A$-$B$-bimodules in $\mathfrak{C}$ and their morphisms can be assembled into a 2-category. Further, as $\mathfrak{C}$ is strict cubical, this 2-category is in fact strict.

\begin{Lemma}
Given two algebras $A$ and $B$ in $\mathfrak{C}$, $A$-$B$-bimodules, $A$-$B$-bimodule 1-morphisms, and $A$-$B$-bimodule 2-morphisms form a strict 2-category, which we denote by $\mathbf{Bimod}_{\mathfrak{C}}(A,B)$.
\end{Lemma}

\subsection{Rigid and Separable Algebras}

Let $\mathfrak{C}$ be a strict cubical monoidal 2-category. A rigid algebra in $\mathfrak{C}$ is an algebra $A$ whose multiplication 1-morphism $m:A\Box A\rightarrow A$ has a right adjoint as an $A$-$A$-bimodule 1-morphism. In particular, we wish to emphasize that this is a property of an algebra, and not additional structure. Let us also remark that this notion was first introduced in \cite{G}, and was first considered in the study of fusion 2-categories in \cite{JFR}. Before giving examples of this notion in the next section, we review the unpacked version of this definition given in section 2.1 of \cite{D7}.

\begin{Definition}
A rigid algebra in $\mathfrak{C}$ consists of:
\begin{enumerate}
    \item An algebra $A$ in $\mathfrak{C}$ as in definition \ref{def:algebra};
    \item A right adjoint $m^*:A\rightarrow A\Box A$  in $\mathfrak{C}$ to the multiplication map $m$ with unit $\eta^m$ and counit $\epsilon^m$ (depicted below as a cup and a cap);
    \item Two 2-isomorphisms
\end{enumerate}
    
\begin{center}
\begin{tabular}{cc}
$\begin{tikzcd}[sep=tiny]
AA \arrow[rrr, "m"] \arrow[ddd, "m^*1"'] &                                     &    & A \arrow[ddd, "m^*"] \\
                                                     &  &  {}    &                              \\
                                                   &       {} \arrow[ur, "\psi^r", Rightarrow]                              &  &                              \\
AAA \arrow[rrr, "1m"']                                &                                     &    & AA,                         
\end{tikzcd}$ & 
$\begin{tikzcd}[sep=tiny]
AA \arrow[rrr, "m"] \arrow[ddd, "1m^*"'] &                                     &    & A \arrow[ddd, "m^*"] \\
                                                     &  &  {}    &                              \\
                                                   &       {} \arrow[ur, "\psi^l", Rightarrow]                              &  &                              \\
AAA \arrow[rrr, "m1"']                                &                                     &    & AA;                       
\end{tikzcd}$
\end{tabular}
\end{center}
    
\begin{enumerate}
\item[] satisfying:
\item [a.] The 2-morphism $\psi^l$ endow $m^*$ with the structure of a left $A$-module 1-morphism:
\end{enumerate}

\settoheight{\prelim}{\includegraphics[width=45mm]{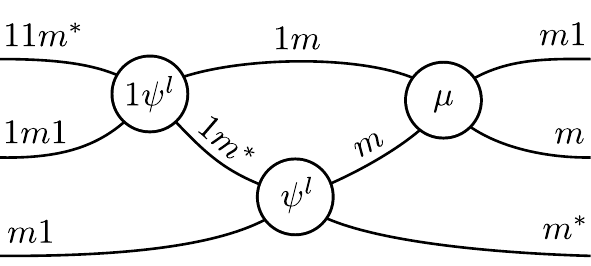}}

\begin{equation}\label{eqn:rigidleftassociativity}
\begin{tabular}{@{}ccc@{}}
\includegraphics[width=45mm]{Pictures/prelim/separablealgebra/leftmodulemap1.pdf}&
\raisebox{0.45\prelim}{$=$} &
\includegraphics[width=45mm]{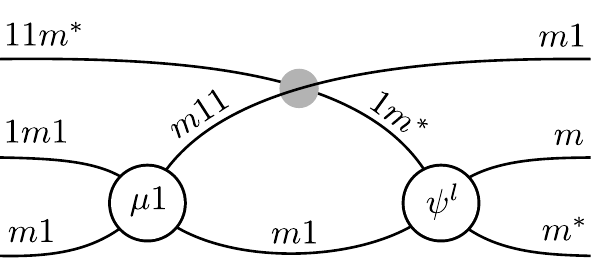},
\end{tabular}
\end{equation}

\settoheight{\prelim}{\includegraphics[width=30mm]{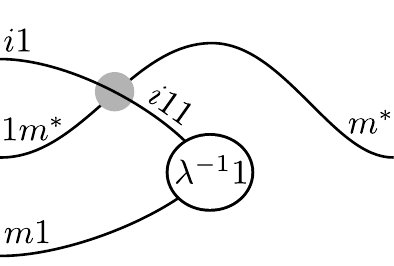}}

\begin{equation}\label{eqn:rigidleftunit}
\begin{tabular}{@{}ccc@{}}
\includegraphics[width=30mm]{Pictures/prelim/separablealgebra/leftmodulemap3.pdf}&
\raisebox{0.45\prelim}{$=$} &
\includegraphics[width=30mm]{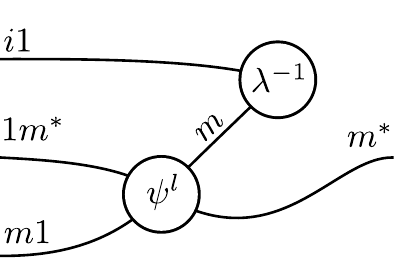},
\end{tabular}
\end{equation}

\begin{enumerate}
\item [b.] The 2-morphism $\psi^r$ endow $m^*$ with the structure of a right $A$-module 1-morphism:
\end{enumerate}

\settoheight{\prelim}{\includegraphics[width=45mm]{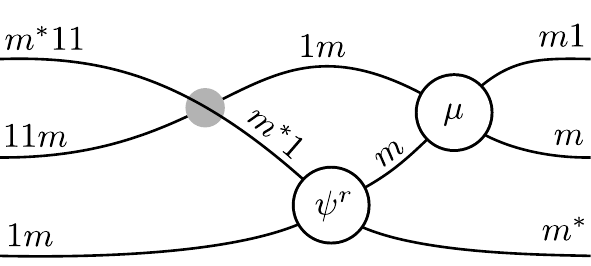}}

\begin{equation}\label{eqn:rigidrightassociativity}
\begin{tabular}{@{}ccc@{}}
\includegraphics[width=45mm]{Pictures/prelim/separablealgebra/rightmodulemap1.pdf}&
\raisebox{0.45\prelim}{$=$} &
\includegraphics[width=45mm]{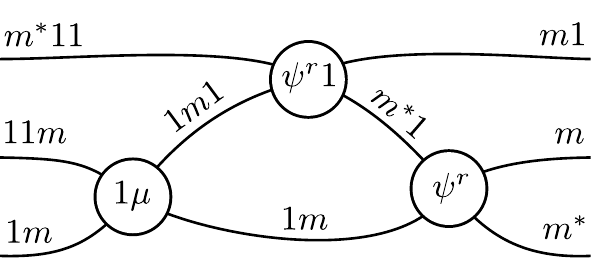},
\end{tabular}
\end{equation}

\settoheight{\prelim}{\includegraphics[width=30mm]{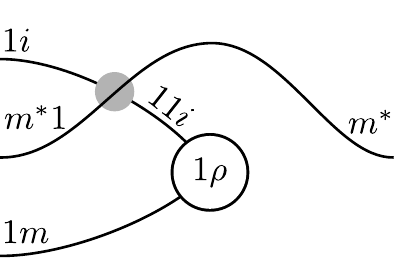}}

\begin{equation}\label{eqn:rigidrightunit}
\begin{tabular}{@{}ccc@{}}
\includegraphics[width=30mm]{Pictures/prelim/separablealgebra/rightmodulemap3.pdf}&
\raisebox{0.45\prelim}{$=$} &
\includegraphics[width=30mm]{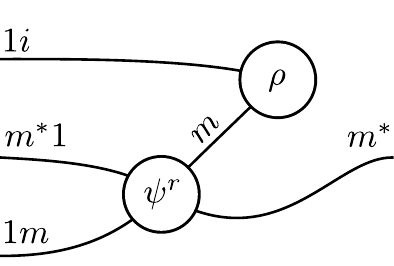},
\end{tabular}
\end{equation}

\begin{enumerate}
\item [c.] The structures of left and right $A$-module 1-morphisms on $m^*$ constructed above are compatible, i.e. they turn $m^*$ into an $A$-$A$-bimodule 1-morphism:
\end{enumerate}

\settoheight{\prelim}{\includegraphics[width=45mm]{Pictures/prelim/separablealgebra/leftmodulemap1.pdf}}

\begin{equation}\label{eqn:rigidbimodule}
\begin{tabular}{@{}ccc@{}}
\includegraphics[width=45mm]{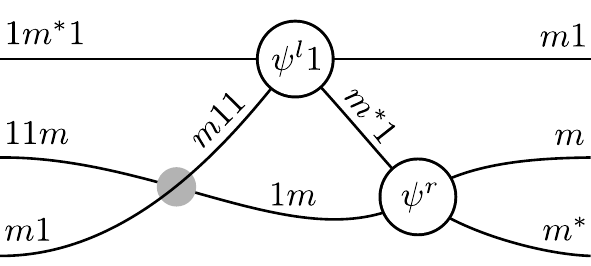}&
\raisebox{0.45\prelim}{$=$} &
\includegraphics[width=45mm]{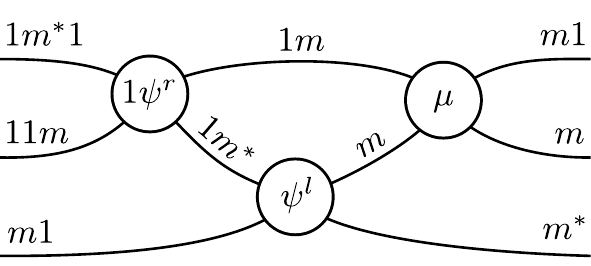},
\end{tabular}
\end{equation}

\begin{enumerate}
\item [d.] The 2-morphism $\epsilon^m$, depicted below as a cap, is an $A$-$A$-bimodule 2-morphism:
\end{enumerate}

\settoheight{\prelim}{\includegraphics[width=22.5mm]{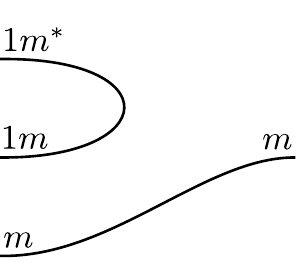}}

\begin{equation}\label{eqn:epsilonleft}
\begin{tabular}{@{}ccc@{}}
\includegraphics[width=22.5mm]{Pictures/prelim/separablealgebra/epsilon1.pdf}&
\raisebox{0.45\prelim}{$=$} &
\includegraphics[width=45mm]{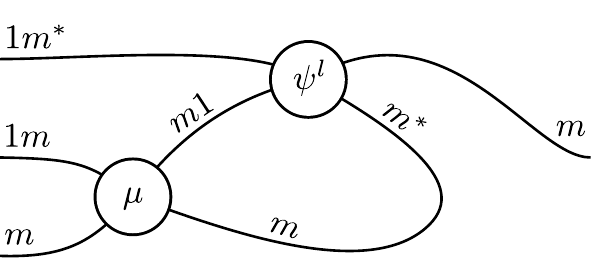},
\end{tabular}
\end{equation}

\settoheight{\prelim}{\includegraphics[width=22.5mm]{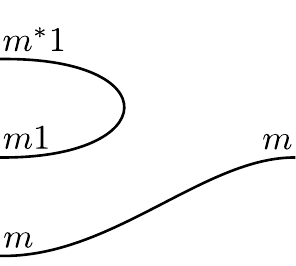}}

\begin{equation}\label{eqn:epsilonright}
\begin{tabular}{@{}ccc@{}}
\includegraphics[width=22.5mm]{Pictures/prelim/separablealgebra/epsilon3.pdf}&
\raisebox{0.45\prelim}{$=$} &
\includegraphics[width=45mm]{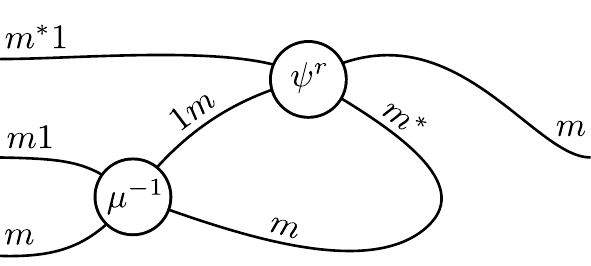},
\end{tabular}
\end{equation}

\begin{enumerate}
\item [e.] The 2-morphism $\eta^m$, depicted below as a cup, is an $A$-$A$-bimodule 2-morphism:
\end{enumerate}

\settoheight{\prelim}{\includegraphics[width=45mm]{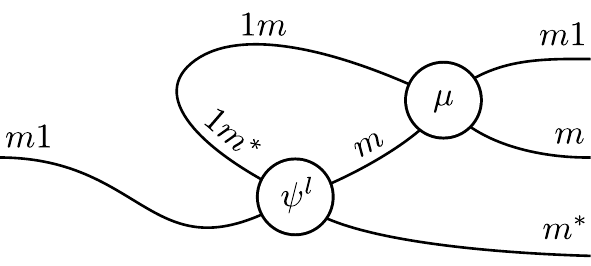}}

\begin{equation}\label{eqn:etaleft}
\begin{tabular}{@{}ccc@{}}
\includegraphics[width=45mm]{Pictures/prelim/separablealgebra/eta1.pdf}&
\raisebox{0.45\prelim}{$=$} &
\includegraphics[width=22.5mm]{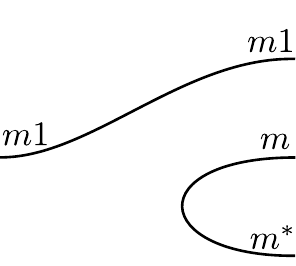},
\end{tabular}
\end{equation}

\settoheight{\prelim}{\includegraphics[width=45mm]{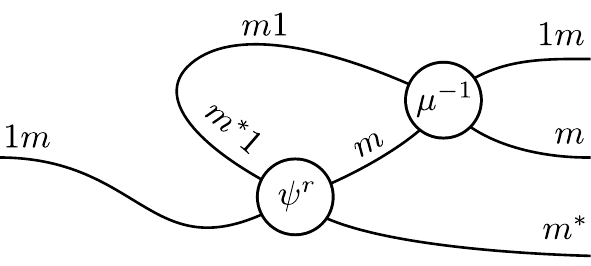}}

\begin{equation}\label{eqn:etaright}
\begin{tabular}{@{}ccc@{}}
\includegraphics[width=45mm]{Pictures/prelim/separablealgebra/eta3.pdf}&
\raisebox{0.45\prelim}{$=$} &
\includegraphics[width=22.5mm]{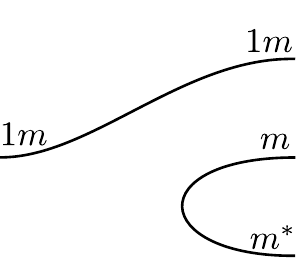}.
\end{tabular}
\end{equation}

\end{Definition}

Following \cite{JFR}, a rigid algebra $A$ in $\mathfrak{C}$ is called separable if the $A$-$A$-bimodule 2-morphism $\epsilon^m:m\circ m^*\Rightarrow Id_A$ as in the above definition has a section as an $A$-$A$-bimodule 2-morphism. Let us again highlight that being separable is a property of an algebra. We now recall the detailed definition of a separable algebra given in definition 2.1.2 of \cite{D7}.

\begin{Definition}
A separable algebra in $\mathfrak{C}$ is a rigid algebra $A$ in $\mathfrak{C}$ equipped with a 2-morphism $\gamma^m:Id_A\Rightarrow m\circ m^*$ such that:

\begin{enumerate}
\item [a.] The 2-morphism $\gamma^m$ is a section of $\epsilon^m$, i.e. $\epsilon^m\cdot\gamma^m = Id_{Id_A}$,
\item [b.] The 2-morphism $\gamma^m$ is an $A$-$A$-bimodule 2-morphism:
\end{enumerate}

\settoheight{\prelim}{\includegraphics[width=45mm]{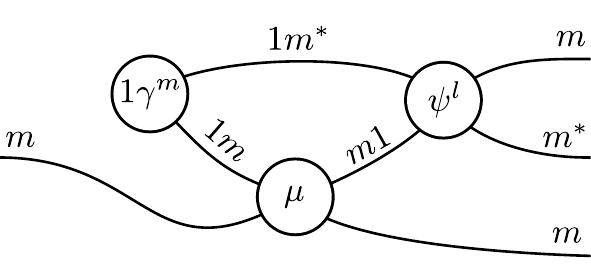}}

\begin{equation}\label{eqn:gammaleft}
\begin{tabular}{@{}ccc@{}}
\includegraphics[width=45mm]{Pictures/prelim/separablealgebra/gamma1.pdf}&
\raisebox{0.45\prelim}{$=$} &
\includegraphics[width=22.5mm]{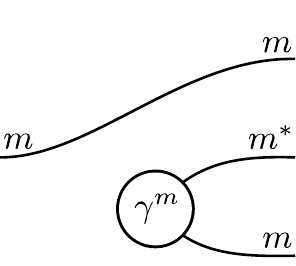},
\end{tabular}
\end{equation}

\settoheight{\prelim}{\includegraphics[width=45mm]{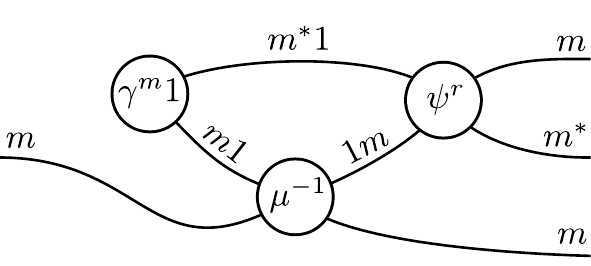}}

\begin{equation}\label{eqn:gammaright}
\begin{tabular}{@{}ccc@{}}
\includegraphics[width=45mm]{Pictures/prelim/separablealgebra/gamma3.pdf}&
\raisebox{0.45\prelim}{$=$} &
\includegraphics[width=22.5mm]{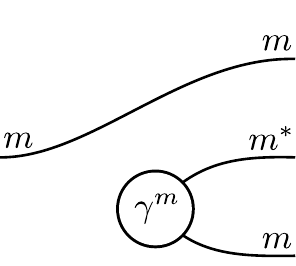}.
\end{tabular}
\end{equation}
\end{Definition}

Let $\mathds{k}$ be a field, and let assume that $\mathfrak{C}$ is a monoidal compact semisimple $\mathds{k}$-linear 2-category. The properties of rigid and separable algebras in $\mathfrak{C}$ have been investigated in details in \cite{D7}. In particular, theorem 3.1.6 of \cite{D7} shows that if $A$ is a rigid algebra in $\mathfrak{C}$, then $A$ is separable if and only if $\mathbf{Bimod}_{\mathfrak{C}}(A)$ is compact semisimple. Further, if either of these conditions is satisfied, both $\mathbf{Mod}_{\mathfrak{C}}(A)$, and $\mathbf{LMod}_{\mathfrak{C}}(A)$ are compact semisimple 2-categories.

\subsection{Examples}

Let $\mathds{k}$ be a field. Following \cite{D7}, we examine rigid and separable algebras in some of the examples of compact semisimple tensor 2-categories given in section \ref{sub:SS2CF2C}. We emphasize that these compact semisimple tensor 2-categories are not strict cubical monoidal 2-category, so that we really have to use the fully weak definition of an algebra in a monoidal 2-category.

\begin{Example}
Algebras in $\mathbf{2Vect}$ are precisely perfect monoidal ($\mathds{k}$-linear) 1-categories, and rigid algebras are precisely perfect tensor 1-categories, i.e. perfect monoidal 1-categories whose objects have right and left duals. Corollary 3.1.7 of \cite{D7} shows that a perfect tensor 1-category $\mathcal{C}$ yields a separable algebra in $\mathbf{2Vect}$ if and only if its Drinfel'd center $\mathcal{Z}(\mathcal{C})$ is a finite semisimple 1-category. If $\mathds{k}$ has characteristic zero, it follows from corollary 2.6.8 of \cite{DSPS13} that finite semisimple tensor 1-categories give all separable algebras in $\mathbf{2Vect}$.
\end{Example}

\begin{Example}\label{ex:algebrastwistedgroupgraded}
Let $G$ be a finite group. Algebras in $\mathbf{2Vect}_G$ are precisely perfect $G$-graded monoidal 1-categories, and rigid algebras are exactly perfect $G$-graded tensor 1-categories. If $\mathds{k}$ has characteristic zero, it is straightforward to check that finite semisimple $G$-graded tensor categories yield all separable algebras in $\mathbf{2Vect}_G$. More generally, given a 4-cocycle for $G$ with coefficients in $\mathds{k}^{\times}$, algebras in $\mathbf{2Vect}_G^{\pi}$ should be thought of as perfect $\pi$-twisted $G$-graded monoidal 1-categories. If $H\subseteq G$ is a subgroup, and $\gamma$ is a 3-cochain for $H$ such that $d\gamma = \pi|_H$, we can consider the algebra $\mathbf{Vect}_H^{\gamma}$ in $\mathbf{2Vect}_G^{\pi}$. It follows from corollary 3.3.7 of \cite{D7} that $\mathbf{Vect}_H^{\gamma}$ yields a rigid algebra in $\mathbf{2Vect}_G^{\pi}$, which is separable if and only if the characteristic of $\mathds{k}$ does not divide the order of $H$.
\end{Example}

\begin{Example}\label{ex:algebrasModB}
Let $\mathcal{B}$ be a braided finite semisimple tensor 1-category. In the terminology of \cite{BJS}, a $\mathcal{B}$-central monoidal 1-category is a monoidal 1-category $\mathcal{C}$ equipped with a braided monoidal functor $F_{\mathcal{C}}:\mathcal{B}\rightarrow \mathcal{Z}(\mathcal{C})$ to the Drinfel'd center of $\mathcal{C}$. Note that this induces in particular a right $\mathcal{B}$-module structure on $\mathcal{C}$. This notion has also appeared under different names in \cite{DGNO}, \cite{HPT} and \cite{MPP}. It follows from proposition 3.2 of \cite{BJS} that algebras in $\mathbf{Mod}(\mathcal{B})$ correspond exactly to finite semisimple $\mathcal{B}$-central monoidal 1-categories, which are separable as right $\mathcal{B}$-module 1-categories. Moreover, by lemma 2.1.4 of \cite{D7}, every $\mathcal{B}$-central finite semisimple tensor 1-category, which is separable as right $\mathcal{B}$-module 1-category, is a rigid algebra in $\mathbf{Mod}(\mathcal{B})$. If $\mathds{k}$ has characteristic zero, if follows from proposition 3.3.3 of \cite{D7} that every $\mathcal{B}$-central finite semisimple tensor 1-categories with simple monoidal unit yields a separable algebra in $\mathbf{Mod}(\mathcal{B})$. For completeness, let us also mention that algebra 1-homomorphisms in $\mathbf{Mod}(\mathcal{B})$ are exactly the monoidal functors over $\mathcal{B}$ described in definition 2.7 of \cite{DNO}.
\end{Example}

\begin{Example}
Let $G$ be a finite group of order coprime to $char(\mathds{k})$. Algebras in $\mathbf{2Rep}(G)$ are given exactly by perfect monoidal 1-categories with a $G$-action. Further, algebra 1-homomorphisms are monoidal functors preserving the $G$-actions, and monoidal natural transformations preserving the $G$-action. Moreover, rigid algebras $\mathbf{2Rep}(G)$ are precisely perfect tensor 1-categories with a $G$-action, and it follows from lemma 3.3.5 of \cite{D7} that such a rigid algebra is separable if and only if the underlying perfect tensor 1-category is separable.
\end{Example}

\begin{Remark}
Lemma \ref{lem:2repG} has one particularly noteworthy consequence, which we now explain. As $\mathbf{Mod}(\mathbf{Rep}(G))$ and $\mathbf{2Rep}(G)$ are equivalent as symmetric monoidal 2-categories, the associated (symmetric monoidal) 2-categories of algebra, algebra 1-homomorphisms and algebra 2-homomorphisms are equivalent. In particular, this induces an equivalence between the full sub-2-categories on the rigid algebras. If we assume that $\mathds{k}$ is an algebraically closed field of characteristic zero, we therefore get an equivalence between the 2-category of multifusion 1-categories with a $G$-action, and the 2-category of $\mathbf{Rep}(G)$-central multifusion 1-categories. In the theory of fusion 1-categories, this is a well-known result (see theorem 4.18 of \cite{DGNO}). In addition, we also get an equivalence between the (symmetric monoidal) 2-categories of braided rigid algebras. That is there is an equivalence between the 2-categories of braided multifusion 1-categories with a braided $G$-action and braided multifusion 1-categories equipped with a braided functor from $\mathbf{Rep}(G)$ into its M\"uger center. This is also a classical result (see proposition 4.22 of \cite{DGNO}).
\end{Remark}

%% file: TensorProduct.tex
\section{The Relative Tensor Product over Separable Algebras}\label{sec:tensorproduct}

Throughout this section, we work with a fixed monoidal 2-category $\mathfrak{C}$, which we assume to be strict cubical without loss of generality. Our first goal is to explain the 2-universal property of the relative tensor product of a right and a left module over an arbitrary algebra $A$. We then prove that if $\mathfrak{C}$ is Karoubi complete and $A$ is separable, then the relative tensor product over $A$ always exists. Using this fact, we construct the Morita 3-category of separable algebras, bimodules, and their morphisms in $\mathfrak{C}$.

\subsection{Definition and Existence}

Let $A$ be an algebra in $\mathfrak{C}$. We fix $M$ a right $A$-module in $\mathfrak{C}$, $N$ a left $A$-module in $\mathfrak{C}$. We begin by defining $A$-balanced morphisms out of the pair $(M,N)$.

\begin{Definition}
Let $C$ be an object of $\mathfrak{C}$. An $A$-balanced 1-morphism $(M,N)\rightarrow C$ consists of:
\begin{enumerate}
    \item A 1-morphism $f:M\Box N\rightarrow C$ in $\mathfrak{C}$;
    \item A 2-isomorphism
\end{enumerate}

$$\begin{tikzcd}[sep=small]
MAN \arrow[dd, "1l^N"'] \arrow[rr, "n^M1"]    &  & MN \arrow[dd, "f"] \\
                                            &  &                      \\
MN \arrow[rr, "f"'] \arrow[rruu, Rightarrow, "\beta^f", shorten > = 2.5ex, shorten < = 2.5ex] &  & C;                 
\end{tikzcd}$$

satisfying:

\begin{enumerate}
\item [a.] We have:
\end{enumerate}

\newlength{\tensor}
\settoheight{\tensor}{\includegraphics[width=52.5mm]{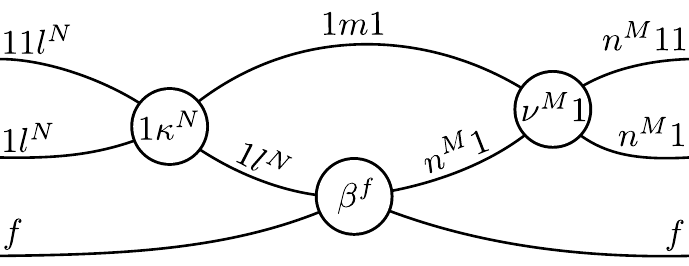}}

\begin{equation}\label{eqn:balancedassociativity}
\begin{tabular}{@{}ccc@{}}

\includegraphics[width=52.5mm]{Pictures/tensor/balanced/associativity1.pdf} & \raisebox{0.45\tensor}{$=$} &
\includegraphics[width=45mm]{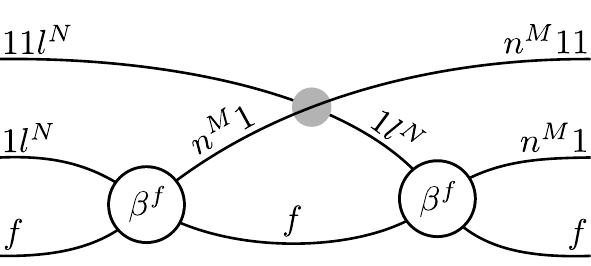},

\end{tabular}
\end{equation}

\begin{enumerate}
\item [b.] We have:
\end{enumerate}

\settoheight{\tensor}{\includegraphics[width=22.5mm]{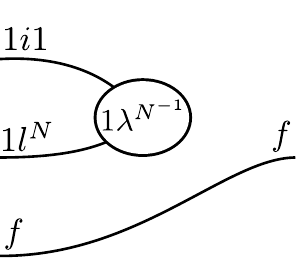}}

\begin{equation}\label{eqn:balancedunitality}
\begin{tabular}{@{}ccc@{}}

\includegraphics[width=22.5mm]{Pictures/tensor/balanced/unitality1.pdf} & \raisebox{0.45\tensor}{$=$} &

\includegraphics[width=37.5mm]{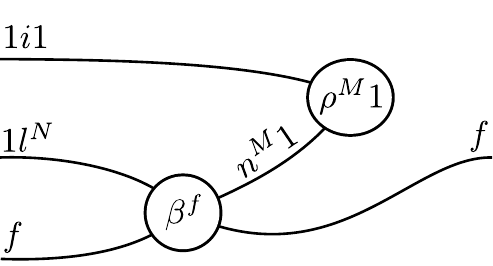}.

\end{tabular}
\end{equation}
\end{Definition}

\begin{Definition}
Let $C$ be an object of $\mathfrak{C}$, and $f,g:(M,N)\rightarrow C$ be two $A$-balanced 1-morphisms. An $A$-balanced 2-morphism $f\Rightarrow g$ is a 2-morphism $\gamma:f\Rightarrow g$ in $\mathfrak{C}$ such that
\settoheight{\tensor}{\includegraphics[width=30mm]{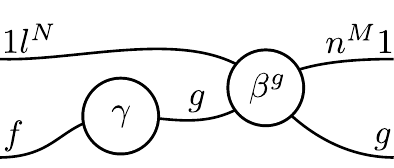}}

\begin{center}
\begin{tabular}{@{}ccc@{}}

\includegraphics[width=30mm]{Pictures/tensor/balanced/2morphism1.pdf} & \raisebox{0.45\tensor}{$=$} &

\includegraphics[width=30mm]{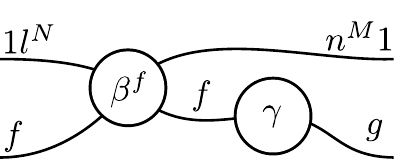}.

\end{tabular}
\end{center}
\end{Definition}

Using the above definitions of $A$-balanced morphisms, we can give the definition of the relative tensor product over $A$.

\begin{Definition}\label{def:relativetensor}
The relative tensor product of $M$ and $N$ over $A$, if it exists, is an object $M\Box_A N$ of $\mathfrak{C}$ together with an $A$-balanced 1-morphism $t_A:(M,N)\rightarrow M\Box_A N$ satisfying the following 2-universal property:
\begin{enumerate}
    \item For every $A$-balanced 1-morphism $f:(M,N)\rightarrow C$, there exists a 1-morphism $\widetilde{f}:M\Box_A N\rightarrow C$ in $\mathfrak{C}$ and an $A$-balanced 2-isomorphism $\xi:\widetilde{f}\circ t_A\cong f$.
    \item For any 1-morphisms $g,h:M\Box_A N\rightarrow C$ in $\mathfrak{C}$, and any $A$-balanced 2-morphism $\gamma:g\circ t_A\Rightarrow h\circ t_A$, there exists a unique 2-morphism $\zeta:g\Rightarrow h$ such that $\zeta\circ t_A = \gamma$.
\end{enumerate}
\end{Definition}

\begin{Remark}\label{rem:relativetensorreformulation}
Observe that, for any object $C$ in $\mathfrak{C}$, $A$-balanced 1-morphisms and 2-morphisms out of $(M,N)$ form a 1-category, which we denote by $Bal_A(M,N;C)$. Furthermore, this assignment is functorial in $M$, $N$, and $C$. Definition \ref{def:relativetensor} may be rephrased as asserting that precomposition with $t_A$ induces an equivalence of 1-categories $$Hom_{\mathfrak{C}}(M\Box_A N,C)\simeq Bal_A(M,N;C),$$ which is natural in the object $C$ in $\mathfrak{C}$. Let us also note that it follows readily from the definition that the 2-category of relative tensor products $M\Box_A N$ is either empty or a contractible 2-groupoid.
\end{Remark}

\begin{Remark}
Over an algebraically closed field, with $\mathfrak{C}=\mathbf{2Vect}$, and $\mathcal{C}$ a multifusion 1-category, then definition \ref{def:relativetensor} recovers the relative tensor product over $\mathcal{C}$ given in definition 3.3 of \cite{ENO2}. As $\mathcal{C}$ is automatically separable in this case, theorem \ref{thm:tensormodules} below recovers the well-known statement that the relative tensor product of two finite semisimple $\mathcal{C}$-module 1-categories exists and is a finite semisimple 1-category. Other particular cases of definition \ref{def:relativetensor} have already appeared in definition 3.2 \cite{DSPS14} and definition 3.3 of \cite{BZBJ}.
\end{Remark}

\begin{Theorem}\label{thm:tensormodules}
Let $A$ be a separable algebra in a Karoubi complete monoidal 2-category $\mathfrak{C}$. Then, the relative tensor product of any right $A$-module $M$ and any left $A$-module $N$ in $\mathfrak{C}$ exists, and is given by the splitting of a 2-condensation monad on $M\Box N$.
\end{Theorem}
\begin{proof}
Let us consider the 2-condensation monad $(M\Box N, e, \mu, \delta)$ in $\mathfrak{C}$ given by $$e:=\big(M\Box l^N\big)\circ \big(n^M\Box A\Box N\big)\circ \big(M\Box (m^*\circ i)\Box N\big),$$ and

\settoheight{\tensor}{\includegraphics[width=105mm]{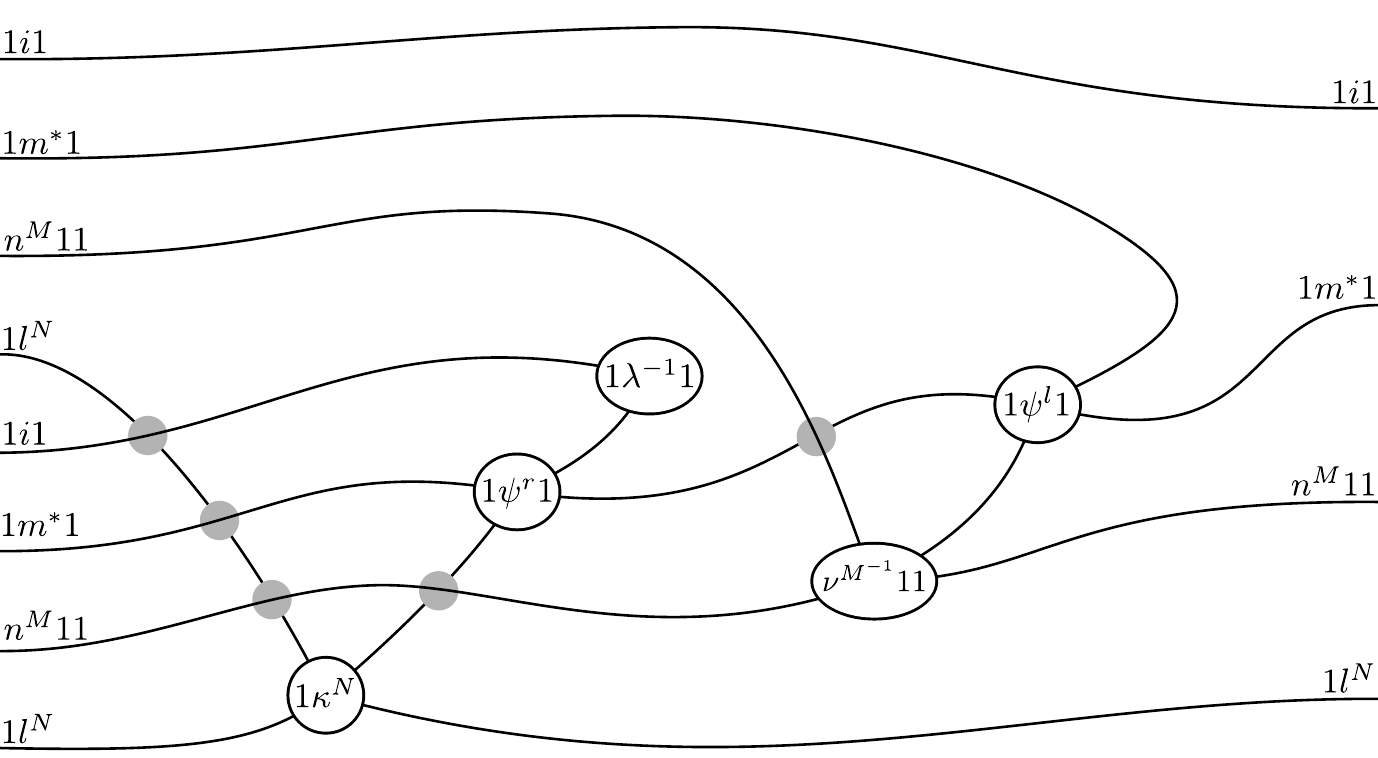}}

\begin{center}
\begin{tabular}{@{}cc@{}}
\raisebox{0.45\tensor}{$\mu:=$} &
\includegraphics[width=105mm]{Pictures/tensor/condensation/mu.pdf},
\end{tabular}
\end{center}

\begin{center}
\begin{tabular}{@{}cc@{}}
\raisebox{0.45\tensor}{$\delta:=$} &
\includegraphics[width=105mm]{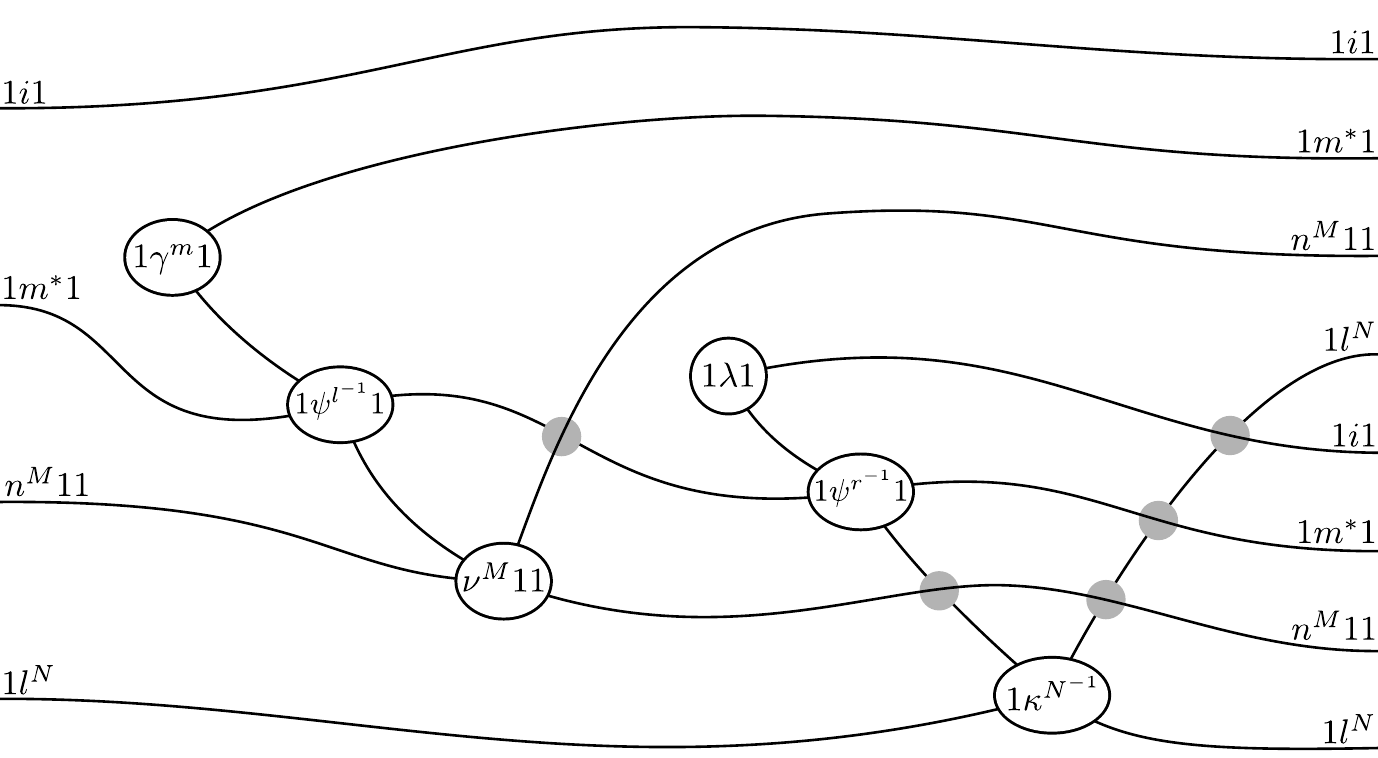}.
\end{tabular}
\end{center}

Clearly, $\mu\cdot\delta = Id_e$. We now prove that $\mu$ is associative using the diagrams depicted in section \ref{sub:condensation}. Figure \ref{fig:muassociativity1} depicts the composite $\mu\cdot (e \circ \mu)$. We begin by moving the two indicated coupons labeled $1\kappa^N$ and $\nu^{M^{-1}}11$ to the left along the corresponding arrows, which brings us to figure \ref{fig:muassociativity2}. We then using equation (\ref{eqn:leftmoduleassociativity}) on the blue coupons and equation (\ref{eqn:moduleassociativity}) on the green coupons to arrive at figure \ref{fig:muassociativity3}. We go on by moving the coupon labeled $11\kappa^N$ up, as well as the coupons labeled $\nu^{M^{-1}}1111$, $\nu^{M^{-1}}111$, and $1\mu^{-1}111$ to the left along the green arrow. Having arrived at figure \ref{fig:muassociativity4}, we move the coupon labeled $\nu^{M^{-1}}11111$ up, and that labelled $\nu^{M^{-1}}1111$ down. Further, we also move the left most cap along the red arrow, and in doing so, use equations (\ref{eqn:epsilonright}) and (\ref{eqn:epsilonleft}), which brings us to figure \ref{fig:muassociativity5}. Now, we use equations (\ref{eqn:rigidbimodule}) on the blue coupons to get to figure \ref{fig:muassociativity6}. We then move the coupon labeled $1\mu^{-1}1111$ to the right, as well as the coupon labeled $11\mu1$ up in order to apply equation (\ref{eqn:rigidrightassociativity}) to the green coupons, and use equation (\ref{eqn:rigidleftassociativity}) on the red coupons, bringing us to figure \ref{fig:muassociativity7}. We can then make use of equation (\ref{eqn:coherenceleft}) on the blue coupons, and cancel the green coupons to arrive at figure \ref{fig:muassociativity8}. Finally, reorganising the diagram along the depicted arrows leads us to figure \ref{fig:muassociativity9}, which represents $\mu\cdot (\mu \circ e)$. Thence, we have established the associativity of $\mu$ as desired. The coassociativity of $\delta$ can be proven similarly.

Let us now move on to proving that $(\mu\circ e)\cdot (e\circ \delta)=\delta \cdot \mu$ using diagrams depicted in section \ref{sub:condensation}. Figure \ref{fig:frobenius1} depicts the left hand-side of this equality. By moving the coupons labeled $1\kappa^{N^{-1}}$ and $\nu^{M}11$ to the right, we arrive at figure \ref{fig:frobenius2}. Then, applying equation (\ref{eqn:leftmoduleassociativity}) to the blue coupons, and equation (\ref{eqn:moduleassociativity}) to the green ones, we get to contemplate figure \ref{fig:frobenius3}. We proceed to move some coupons along the depicted arrows, and use equation (\ref{eqn:rigidrightassociativity}) on the blue coupons, and equation (\ref{eqn:rigidleftassociativity}) on the green coupons, which brings us to figure \ref{fig:frobenius4}. Using equation (\ref{eqn:coherenceleft}) on the blue coupons, and moving the coupons labeled $1\psi^{r^{-1}}1$ and $1\kappa^{N^{-1}}$ to the right yields the diagram given in figure \ref{fig:frobenius5}. Then, we first apply equation (\ref{eqn:rigidbimodule}) to the blue coupons, and then equation (\ref{eqn:algebraassociativity}) on the green coupon together with the coupon labeled $1\mu1$, which was just created. This brings us to figure \ref{fig:frobenius6}. Finally, using in succession equation (\ref{eqn:epsilonright}) on the blue coupons, equation (\ref{eqn:epsilonleft}) on the green coupons, and equation (\ref{eqn:coherenceleft}) on the red coupons, leads us to figure \ref{fig:frobenius7}, which depicts $\delta \cdot \mu$. This proves the desired equality. The equality $(e\circ \mu)\cdot (\delta\circ e)=\delta \cdot \mu$ can be proven using a similar argument.

In order to prove that the relative tensor product of $M$ and $N$ over $A$ exists, we will use the reformulation given in remark \ref{rem:relativetensorreformulation}. To this end, recall that 2-condensation monads are preserved by all 2-functors, so that applying $Hom_{\mathfrak{C}}(-,C)$ to $(M\Box N, e, \mu, \delta)$ yields a 2-condensation monad on the 1-category $Hom_{\mathfrak{C}}(M\Box N,C)$. In fact, this yields a 2-condensation monad on the 2-functor $Hom_{\mathfrak{C}}(M\Box N,-)$. We claim that $Bal_A(M,N;C)$ is a splitting this 2-condensation monad. Namely, let $U:Bal_A(M,N;C)\rightarrow Hom_{\mathfrak{C}}(M\Box N,C)$ be the forgetful functor, and $E:Hom_{\mathfrak{C}}(M\Box N,C)\rightarrow Bal_A(M,N;C)$ be the functor given by $g\mapsto g\circ e$, with $A$-balanced structure on the composite $g\circ e$ supplied by the 2-isomorphism $\beta^{g\circ e}$ given by

\settoheight{\tensor}{\includegraphics[width=105mm]{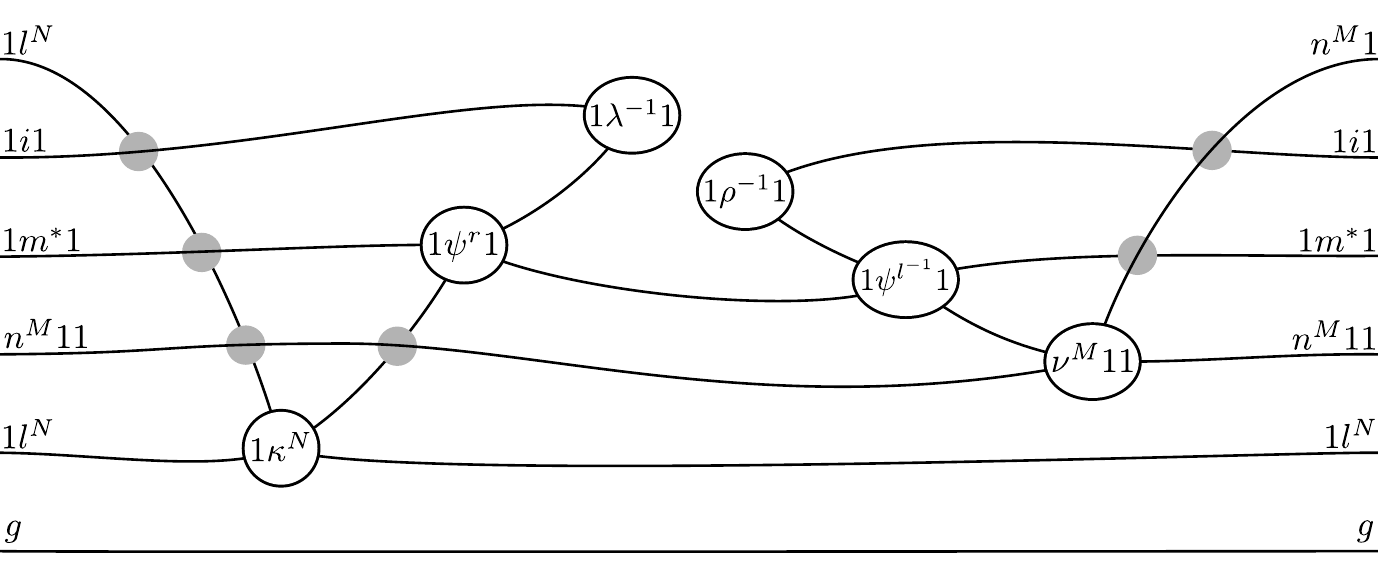}}

\begin{center}
\begin{tabular}{@{}cc@{}}
\raisebox{0.45\tensor}{$\beta^{g\circ e}:=$} &
\includegraphics[width=105mm]{Pictures/tensor/condensation/Abalancedge.pdf}.
\end{tabular}
\end{center}

\noindent The fact that this defines an $A$-balanced structure can be seen as follows. Let us start with the right hand-side of equation (\ref{eqn:balancedassociativity}) for $\beta^{g\circ e}$. We begin by applying equation (\ref{eqn:rigidbimodule}) after having moved some coupons, then we use equation (\ref{eqn:algebraunitality}). We continue by appealing to equations (\ref{eqn:moduleassociativity}) and (\ref{eqn:leftmoduleassociativity}), followed by (\ref{eqn:rigidrightassociativity}) and (\ref{eqn:rigidleftassociativity}). At last, we can use equations (\ref{eqn:coherenceleft}) and (\ref{eqn:coherenceright}) for $A$ as well as reorganise the string diagram to get to the left hand-side of (\ref{eqn:balancedassociativity}). Equation (\ref{eqn:balancedunitality}) for $\beta^{g\circ e}$ follows similarly. Now, observe that both $U$ and $A$ are 2-natural in $C$. Further, let us define natural transformations $p:E\circ U\Rightarrow Id$ and $s:Id\Rightarrow E\circ U$ by

\settoheight{\tensor}{\includegraphics[width=67.5mm]{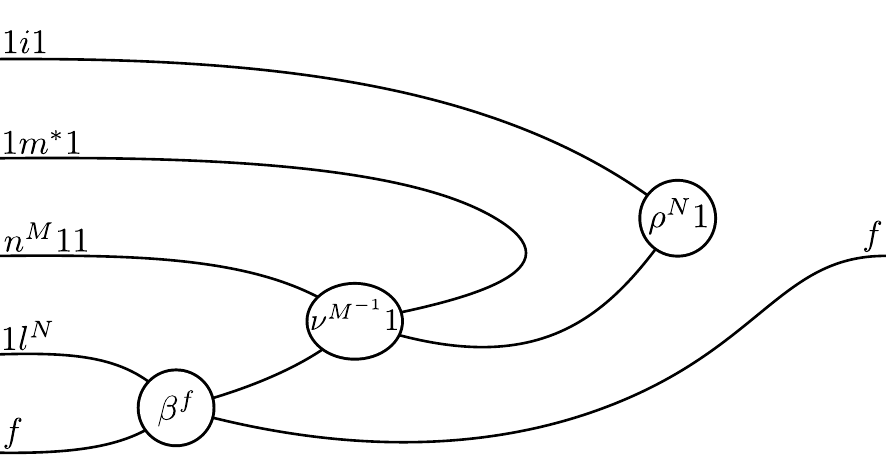}}

\begin{center}
\begin{tabular}{@{}cc@{}}
\raisebox{0.45\tensor}{$p_f:=$} &
\includegraphics[width=67.5mm]{Pictures/tensor/condensation/pf.pdf},
\end{tabular}
\end{center}

\begin{center}
\begin{tabular}{@{}cc@{}}
\raisebox{0.45\tensor}{$s_f:=$} &
\includegraphics[width=67.5mm]{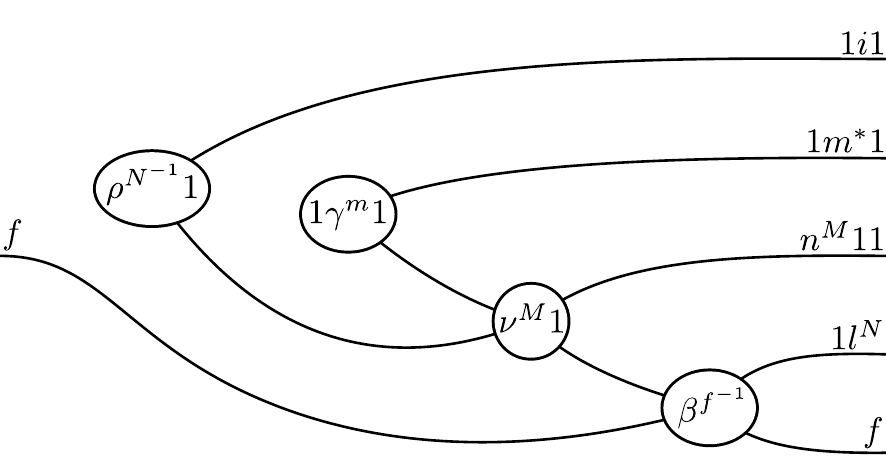},
\end{tabular}
\end{center}

\noindent for every $A$-balanced 1-morphism $f:(M,N)\rightarrow C$. Again, note that $s$ and $p$ are 2-natural in $C$. Further, we have $p\cdot s = Id$, so that $$(Hom_{\mathfrak{C}}(M\Box N,-), Bal_A(M,N;-),E,U,p,s)$$ is a 2-condensation. It remains to check that this splits the 2-condensation monad on $Hom_{\mathfrak{C}}(M\Box N,-)$ induced by  $(M\Box N, e, \mu, \delta)$. To see this, it is enough to prove that for every 1-morphism $g:M\Box N\rightarrow C$ in $\mathfrak{C}$, we have $p_{g\circ e} = g \cdot \mu$ and $s_{g\circ e} = g \cdot \delta$. The first equality follows by applying equations (\ref{eqn:moduleassociativity}), (\ref{eqn:rigidleftassociativity}), followed by equation (\ref{eqn:coherenceright}) for $A$, and then by using successively equations (\ref{eqn:moduleunitality}), (\ref{eqn:rigidleftunit}), and (\ref{eqn:coherencemiddle}). The second equality is obtained in a similar fashion.

Finally, as $\mathfrak{C}$ is Karoubi complete, the 2-condensation monad $(M\Box N, e, \mu, \delta)$ admits a splitting in $\mathfrak{C}$, which we denote by $M\Box_A N$. Now, the splitting of a 2-condensation monad is preserved by any 2-functor, so that $Hom_{\mathfrak{C}}(M\Box_A N,-)$ is also a splitting of the 2-condensation monad on $Hom_{\mathfrak{C}}(M\Box N,-)$ induced by $(M\Box N, e, \mu, \delta)$. But, the 2-category of splittings of a 2-condensation monad is a contractible 2-groupoid, so that we get the desired equivalence.
\end{proof}

\begin{Remark}\label{rem:relativetensorcodescent}
In the language of \cite{CMV}, the 1-category $Bal_A(M,N;C)$ is the pseudo-coequalizer for the descent object $$\begin{tikzcd}
{Hom_{\mathfrak{C}}(MN,C)} \arrow[rr, shift right=2] \arrow[rr, shift left=2] &  & {Hom_{\mathfrak{C}}(MAN,C)} \arrow[ll] \arrow[rr] \arrow[rr, shift left=2] \arrow[rr, shift right=2] &  & {Hom_{\mathfrak{C}}(MAAN,C)}
\end{tikzcd}$$
obtained by applying $Hom_{\mathfrak{C}}(-,C)$ to the canonical codescent object
$$\begin{tikzcd}
M\Box A\Box A\Box N \arrow[rr, shift right=2] \arrow[rr, shift left=2] \arrow[rr] &  & M\Box A\Box N \arrow[rr, shift left=2] \arrow[rr, shift right=2] &  & M\Box N. \arrow[ll]
\end{tikzcd}$$
Theorem \ref{thm:tensormodules} shows that the 2-functor $Bal_A(M,N;-)$ is corepresented by $M\Box_A N$, so that $M\Box_A N$ is the pseudo-coequalizer of the above codescent object.
\end{Remark}

Thanks to the definition of the relative tensor product using a 2-universal property, the following result is an immediate consequence of the above theorem.

\begin{Corollary}\label{cor:tensorfunctoriality}
If $\mathfrak{C}$ is a Karoubi complete 2-category, and $A$ is a separable algebra, the relative tensor product over $A$ defines a 2-functor $$\Box_A:\mathbf{Mod}_{\mathfrak{C}}(A)\times \mathbf{LMod}_{\mathfrak{C}}(A)\rightarrow \mathfrak{C}.$$
\end{Corollary}

\begin{Remark}
For completeness, let us note that if $\mathfrak{C}$ is a linear monoidal 2-category, then it follows from the 2-universal property of $\Box_A$ and the fact that $\Box$ is a bilinear 2-functor that $\Box_A$ is a bilinear 2-functor.
\end{Remark}

\subsection{The Morita 3-Category}

Our goal is now to explain how to construct the Morita 3-category of separable algebras in a Karoubi complete monoidal 2-category. In order to do so, we need to generalize the setup of the previous section to bimodules.

\begin{Definition}
Let $A$, $B$, $C$ be algebras in $\mathfrak{C}$, and let $M$ be an $A$-$B$-bimodule, $N$ be a $B$-$C$-bimodule, and $P$ be an $A$-$C$-bimodule. A $B$-balanced $A$-$C$-bimodule 1-morphism $(M,N)\rightarrow P$ is an $A$-$C$-bimodule 1-morphism $f:M\Box N\rightarrow P$ together with an $A$-$C$-bimodule 2-isomorphism $\beta^f:f\circ (M\Box l^N)\cong f\circ (n^M\Box N)$ providing $f$ with an $A$-balanced structure. A $B$-balanced $A$-$C$-bimodule 2-morphism is an $A$-$C$-bimodule 2-morphism that is also $B$-balanced.
\end{Definition}

\begin{Proposition}\label{prop:relativetensorbimodules}
Let $A$, $B$, $C$ be algebras in $\mathfrak{C}$, with $B$ separable. Let $M$ be an $A$-$B$-bimodule, and $N$ be a $B$-$C$-bimodule, the relative tensor product $t_B:M\Box N\rightarrow M\Box_B N$ can be endowed with an $A$-$C$-bimodule structure such that it is 2-universal with respect to $B$-balanced $A$-$C$-bimodule morphisms.
\end{Proposition}
\begin{proof}
Note that if $M$ and $N$ are bimodules in the proof of theorem \ref{thm:tensormodules}, then the 2-condensation monad $(M\Box N, e, \mu, \delta)$ in $\mathfrak{C}$ can be upgraded to a 2-condensation monad in $\mathbf{Bimod}_{\mathfrak{C}}(A,C)$. The remainder of the proof can be straightforwardly adapted to accommodate for the bimodule case. The only noteworthy change is that one needs to use the fact that $\mathbf{Bimod}_{\mathfrak{C}}(A,C)$ is Karoubi complete, which follows from the proof of proposition 3.3.5 of \cite{D4} as $\mathfrak{C}$ is Karoubi complete. In particular, this constructs a 2-universal $B$-balanced $A$-$C$-bimodule 1-morphism $\widetilde{t}_B:M\Box N\rightarrow M\Box_B N$. But, as splittings of 2-condensation monads are preserved by all 2-functors, the underlying $B$-balanced 1-morphism $\widetilde{t}_B:M\Box N\rightarrow M\Box_B N$ in $\mathfrak{C}$ satisfies the 2-universal property of definition \ref{def:relativetensor}. This finishes the proof of the proposition.
\end{proof}

\begin{Remark}
Let us sketch an alternative proof of proposition \ref{prop:relativetensorbimodules}. It follows from the construction of theorem \ref{thm:tensormodules} and the fact that 2-condensation are preserved by all 2-functors that $A\Box t_B:M\Box N\rightarrow A\Box (M\Box_B N)$ is 2-universal with respect to $B$-balanced 1-morphisms. The 2-universal property of the relative tensor product over $B$ can then be used repeatedly to endow $t_B:M\Box N\rightarrow M\Box_B N$ with a left $A$-module structure. Similarly, we can construct a right $C$-module structure on $t_B$, which is compatible with the left $A$-module structure. Finally, one can directly check that the $B$-balanced $A$-$C$-bimodules 1-morphism $t_B$ is 2-universal with respect to $B$-balanced $A$-$C$-bimodule morphisms.
\end{Remark}

\begin{Corollary}\label{cor:relativetensorbimodules}
Let $A$, $B$, $C$ be arbitrary algebras in $\mathfrak{C}$ with $B$ separable. The relative tensor product over $B$ induces a 2-functor $$\Box_B:\mathbf{Bimod}_{\mathfrak{C}}(A,B)\times \mathbf{Bimod}_{\mathfrak{C}}(B,C)\rightarrow \mathbf{Bimod}_{\mathfrak{C}}(A,C).$$
\end{Corollary}

We now prove a unitality property of the relative tensor product that will play a crucial role later on.

\begin{Lemma}\label{lem:relativetensorA}
Let $A$ and $B$ be arbitrary algebras in $\mathfrak{C}$. There is a 2-natural adjoint equivalence $$l^{\mathbf{M}}_P:A\Box_A P\simeq P$$ for any $A$-$B$-bimodule $P$ in $\mathfrak{C}$.
\end{Lemma}
\begin{proof}
Let $P$, $Q$ be two $A$-$B$-bimodule in $\mathfrak{C}$. Observe that $l^P:A\Box P\rightarrow P$ is an $A$-balanced $A$-$B$-bimodule 1-morphism via $\beta^{l^P}:= \kappa^P$. We claim that this 1-morphism satisfies the 2-universal property defining the relative tensor product. Namely, given $f:A\Box P\rightarrow Q$ an $A$-balanced $A$-$B$-bimodule 1-morphism, we define $g$ as the composite right $B$-module 1-morphism $$g:P\xrightarrow{i\Box P}A\Box P\xrightarrow{f} Q.$$ In addition, the 2-isomorphism

\settoheight{\tensor}{\includegraphics[width=52.5mm]{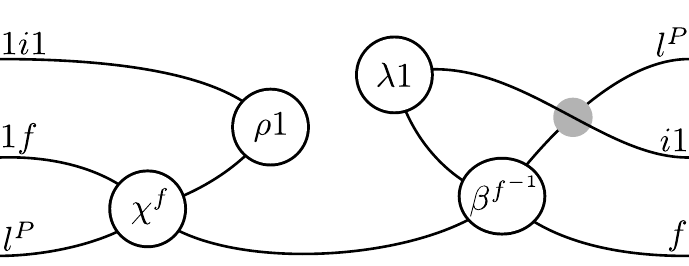}}

\begin{center}
\begin{tabular}{@{}cc@{}}

 \raisebox{0.45\tensor}{$\chi^g:=$} &

\includegraphics[width=52.5mm]{Pictures/tensor/unit/chig.pdf}

\end{tabular}
\end{center}

\noindent endows $g$ with a compatible left $A$-module structure. Further, it follows from the definitions that the 2-isomorphism $\xi:g\circ l^P\cong f$ given by

\settoheight{\tensor}{\includegraphics[width=30mm]{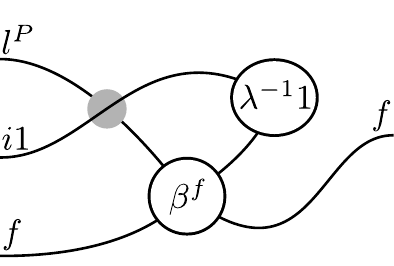}}

\begin{center}
\begin{tabular}{@{}cc@{}}

 \raisebox{0.45\tensor}{$\xi:=$} &

\includegraphics[width=30mm]{Pictures/tensor/unit/xi.pdf}

\end{tabular}
\end{center}

\noindent is an $A$-balanced $A$-$B$-bimodule 2-morphism as desired. Now, let $g,h:P\rightarrow Q$ be two $A$-$B$-bimodule 1-morphisms, and $\gamma: g\circ l^P\Rightarrow h\circ l^P$ be an $A$-balanced $A$-$B$-bimodule 2-morphisms, then it is not hard to check that $\zeta:= \gamma \circ (i\Box P)$ is an $A$-$B$-bimodule 2-morphism satisfying $\zeta \circ l^p =\gamma$. This finishes the proof of the claim. Finally, using the 2-universal property of the relative tensor product, one can readily construct the desired adjoint 2-natural equivalence $l^{\mathbf{M}}$.
\end{proof}

For our purposes, it is also necessary to examine the relative tensor product of multiple bimodules.

\begin{Definition}\label{def:BCbalanced}
Let $A$, $B$, $C$, $D$ be algebras in $\mathfrak{C}$, and let $M$ be an $A$-$B$-bimodule, $N$ be a $B$-$C$-bimodule, $P$ be a $C$-$D$-bimodule, and $Q$ an $A$-$D$-bimodule in $\mathfrak{C}$. A $(B,C)$-balanced $A$-$D$-bimodule 1-morphism $(M,N,P)\rightarrow Q$ consists of:
\begin{enumerate}
    \item An $A$-$D$-bimodule 1-morphism $f:M\Box N\Box P\rightarrow Q$,
    \item Two $A$-$D$-bimodule 2-isomorphisms $\beta^f_B$ and $\beta^f_C$ given by
\end{enumerate}
$$\begin{tabular}{@{} c c @{}}
$\begin{tikzcd}[sep=small]
MBNP \arrow[dd, "1l^N1"'] \arrow[rr, "n^M11"]    &  & MNP \arrow[dd, "f"] \\
                                            &  &                      \\
MNP \arrow[rr, "f"'] \arrow[rruu, Rightarrow, "\beta^f_B", shorten > = 2.5ex, shorten < = 2.5ex] &  & Q,                   
\end{tikzcd}$  &  $\begin{tikzcd}[sep=small]
MBNP \arrow[dd, "11l^P"'] \arrow[rr, "1n^N1"]    &  & MNP \arrow[dd, "f"] \\
                                            &  &                      \\
MNP \arrow[rr, "f"'] \arrow[rruu, Rightarrow, "\beta^f_C", shorten > = 2.5ex, shorten < = 2.5ex] &  & Q,                   
\end{tikzcd}$
\end{tabular}$$

satisfying:
\begin{enumerate}
\item [a.] The 2-isomorphism $\beta^f_B$ endows $f:(M,N\Box P)\rightarrow Q$ with a $B$-balanced structure,
\item [b.] The 2-isomorphism $\beta^f_C$ endows $f:(M\Box N, P)\rightarrow Q$ with a $C$-balanced structure,
\item [c.] The 2-isomorphisms $\beta^f_B$ and $\beta^f_C$ commute in the sense that
\end{enumerate}

\settoheight{\tensor}{\includegraphics[width=45mm]{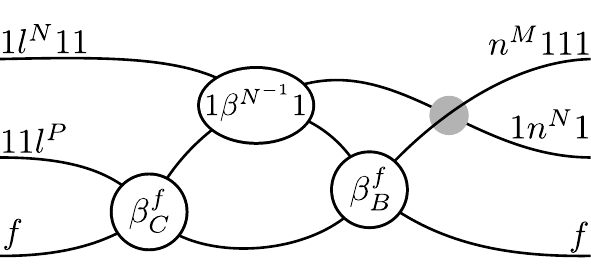}}

\begin{equation}\label{eqn:BCbalancedcompatibility}
\begin{tabular}{@{}ccc@{}}

\includegraphics[width=45mm]{Pictures/tensor/unit/multibalanced1.pdf} & \raisebox{0.45\tensor}{$=$} &

\includegraphics[width=45mm]{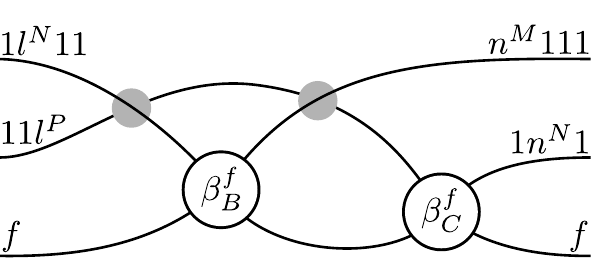}.

\end{tabular}
\end{equation}

A $(B,C)$-balanced $A$-$D$-bimodule 2-morphism is an $A$-$D$-bimodule 2-morphism that is both $B$-balanced and $C$-balanced.
\end{Definition}

\begin{Lemma}\label{lem:relativetensorassociativity}
Let $A$, $B$, $C$, $D$ be algebras in $\mathfrak{C}$, and let $M$ be an $A$-$B$-bimodule, $N$ be a $B$-$C$-bimodule, and $P$ be a $C$-$D$-bimodule. If $\mathfrak{C}$ is Karoubi complete, and $B$, $C$ are separable algebras, then both $M\Box_B(N\Box_C P)$ and $(M\Box_B N)_CP$ are 2-universal with respect to $(B,C)$-balanced $A$-$D$-bimodule morphisms. In particular, there exists an adjoint 2-natural equivalence $$\alpha^{\mathbf{M}}_{M,N,P}:(M\Box_B N)\Box_CP\simeq M\Box_B(N\Box_C P).$$
\end{Lemma}
\begin{proof}
Let us show that the $(B,C)$-balanced $A$-$D$-bimodule 1-morphism $$M\Box N\Box P\xrightarrow{M\Box t_C} M\Box (N\Box_C P)\xrightarrow{t_B} M\Box_B (N\Box_C P)$$ is 2-universal with respect to $(B,C)$-balanced $A$-$D$-bimodule morphisms. In order to prove this, note that the $C$-balanced $A$-$D$-bimodule 1-morphism $M\Box t_C:M\Box N\Box P\rightarrow M\Box (N\Box_C P)$ is 2-universal with respect to $C$-balanced $A$-$D$-bimodule morphism.

Now, let $Q$ be an $A$-$D$-bimodule, and let $f:(M,N,P)\rightarrow Q$ be a $(B,C)$-balanced $A$-$D$-bimodule 1-morphism. This gives us the solid arrow part of the diagram below
$$\begin{tikzcd}[sep=small]
                                                            &  & M\Box (N\Box_C P) \arrow[rrd, "t_B"] \arrow[dd, "f'", dotted] &  &                                                          \\
M\Box N\Box P \arrow[rrd, "f"'] \arrow[rru, "M\Box t_C"] &  &                                                                  &  & M\Box_B (N\Box_C P) \arrow[lld, "\widetilde{f}", dotted] \\
                                                            &  & Q.                                                               &  &                                                         
\end{tikzcd}$$
\noindent As $f$ is in particular a $C$-balanced $A$-$D$-bimodule 1-morphisms, there exists an $A$-$D$-bimodule 1-morphism $f':M\Box (N\Box_C P)\rightarrow Q$, and a $C$-balanced $A$-$D$-bimodule 2-isomorphism $\xi_C:f'\circ (M\Box t_C)\cong f$. But, thanks to equation (\ref{eqn:BCbalancedcompatibility}) and the 2-universal property of $M\Box t_C$, the $B$-balanced structure of $f$ induces a $B$-balanced structure on $f'$. Thus, there exists an $A$-$D$-bimodule 1-morphisms $\widetilde{f}:M\Box_B (N\Box_C P)\rightarrow Q$, and a $B$-balanced $A$-$D$-bimodule 2-isomorphism $\xi_B:\widetilde{f}\circ t_B\cong f'$. It follows from the definitions that the composite $A$-$D$-bimodule 2-isomorphism $\xi_C\cdot(\xi_B\circ (M\Box t_C))$ is $(B,C)$-balanced, so that $\widetilde{f}$ is the sought after factorization of $f$.

Finally, let $g,h:M\Box_B(N\Box_CP)\rightarrow Q$ be two $A$-$D$-bimodule 1-morphisms, and $\gamma:g\circ t_B\circ (M\Box t_C)\Rightarrow h\circ t_B\circ (M\Box t_C)$ a $(B,C)$-balanced $A$-$D$-bimodule 2-morphism. It follows immediately from the 2-universal property of $M\Box t_C$ that there exists an $A$-$D$-bimodule 2-morphism $\zeta':g\circ t_B\Rightarrow h\circ t_B$ such that $\zeta\circ (M\Box t_C)=\gamma$. But, using the 2-universal property of $M\Box t_C$ again together with the fact that $\gamma$ is $B$-balanced, we find that $\zeta'$ is necessarily $B$-balanced. Thence, by the 2-universal property of $t_B$, there exists an $A$-$D$-bimodule 2-morphism $\zeta:g\Rightarrow h$ such that $\zeta\circ t_B = \zeta'$. Putting everything together, we find that $\zeta \circ t_B \circ (M\Box t_C) = \gamma$ as desired. This proves that $M\Box_B (N\Box_C P)$ is 2-universal with respect to $(B,C)$-balanced $A$-$D$-bimodule morphisms.

One proceeds analogously to show that the $(B,C)$-balanced $A$-$D$-bimodule 1-morphism $$M\Box N\Box P\xrightarrow{t_B \Box P} (M\Box_B N)\Box P\xrightarrow{t_C} (M\Box_B N)\Box_C P$$ is 2-universal with respect to $(B,C)$-balanced $A$-$D$-bimodule morphisms. The second part of the statement then follows readily by appealing to the 2-universal property.
\end{proof}

We are now ready to explain the main construction of this section.

\begin{Theorem}\label{thm:separablealgebra3category}
Let $\mathfrak{C}$ be a Karoubi complete monoidal 2-category. Separable algebras in $\mathfrak{C}$, bimodules, bimodule 1-morphisms, and bimodule 2-morphisms form a 3-category, which we denote by $\mathbf{Mor}^{sep}(\mathfrak{C})$.
\end{Theorem}
\begin{proof}
Let $A$, $B$, $C$, be separable algebras in $\mathfrak{C}$. We set $$Hom_{\mathbf{Mor}^{sep}(\mathfrak{C})}(B,A):=\mathbf{Bimod}(A,B).$$ Then, the bilinear 2-functor $$\Box_B:\mathbf{Bimod}_{\mathfrak{C}}(A,B)\times \mathbf{Bimod}_{\mathfrak{C}}(B,C)\rightarrow \mathbf{Bimod}_{\mathfrak{C}}(A,C)$$ of corollary \ref{cor:relativetensorbimodules} provides us with the necessary composition 2-functor. Further, the identity 1-morphism on the algebra $A$ is given by the canonical $A$-$A$-bimodule $A$. It remains to prove that these operations can be made suitably coherent in the sense of definition 4.1 of \cite{Gur}. Firstly, note that lemma \ref{lem:relativetensorA} provides us with an adjoint 2-natural equivalence $l^{\mathbf{M}}$. Using a similar argument, one can construct a 2-natural equivalence $r^{\mathbf{M}}$ given on the $A$-$B$-bimodule $P$ by $r^{\mathbf{M}}_P:P\Box_BB\simeq P$. Moreover, lemma \ref{lem:relativetensorassociativity} provides us with an adjoints 2-natural equivalence $\alpha^{\mathbf{M}}$ witnessing associativity of the composition of 1-morphisms.

Secondly, we have to supply invertible modifications $\lambda^{\mathbf{M}}$, $\mu^{\mathbf{M}}$, $\rho^{\mathbf{M}}$, and $\pi^{\mathbf{M}}$ between specific composites of $l^{\mathbf{M}}$, $r^{\mathbf{M}}$, and $\alpha^{\mathbf{M}}$. Let us explain how to construct $\lambda^{\mathbf{M}}$. Let $M$ be an $A$-$B$-bimodule and $N$ a $B$-$C$-bimodule in $\mathfrak{C}$, and consider the diagram
$$\begin{tikzcd}
                                                                                                          &  & A\Box M\Box N \arrow[lldd] \arrow[rrdd] \arrow[ddd]           &  &          \\
                                                                                                          &  &                                                               &  &          \\
(A\Box_A M)\Box_B N \arrow[rrd, "{\alpha^{\mathbf{M}}_{A,M,N}}"'] \arrow[rrrr, "l^{\mathbf{M}}_M\Box_BN", near end] &  &                                                               &  & M\Box_BN \\
                                                                                                          &  & A\Box_A (M\Box_B N), \arrow[rru, "l^{\mathbf{M}}_{M\Box_BN}"']\arrow[from= uuu,crossing over] &  &         
\end{tikzcd}$$
\noindent where the three unlabeled arrows are the canonical $(A,B)$-balanced $A$-$B$-bimodule 1-morphisms, and the three top triangles are filled by canonical $(A,B)$-balanced $A$-$B$-bimodule 2-isomorphisms. Thanks to the 2-universal property of $A\Box M\Box N\rightarrow (A\Box_A M)\Box_B N$, there exists an $A$-$B$-bimodule 2-isomorphism $$\lambda^{\mathbf{M}}_{M,N}:l^{\mathbf{M}}_M\Box_B N\cong l^{\mathbf{M}}_{M\Box_B N}\circ \alpha^{\mathbf{M}}_{A,M,N}.$$ Using the 2-universal property again, it is easy to check that these 2-isomorphisms define an invertible modification. The invertible modifications $\mu^{\mathbf{M}}$ and $\rho^{\mathbf{M}}$ are constructed similarly.

It remains to construct the invertible modification $\pi^{\mathbf{M}}$. Given separable algebras $A$, $B$, $C$, $D$, $E$, one defines $(B,C,D)$-balanced $A$-$E$-bimodule morphisms by adapting definition \ref{def:BCbalanced} in the obvious way. Following the proof of lemma \ref{lem:relativetensorassociativity}, one then shows that for any $A$-$B$-bimodule $M$, $B$-$C$-bimodule $N$, $C$-$D$-bimodule $P$, and $D$-$E$-bimodule $Q$, the canonical $(B,C,D)$-balanced $A$-$E$-bimodule 1-morphisms to the different ways of parenthesising $M\Box_B N\Box_C P\Box_D Q$ are all 2-universal with respect $(B,C,D)$-balanced $A$-$E$-bimodule morphisms. Analogously to the above arguments, $\pi^{\mathbf{M}}$ is constructed using this 2-universal property.

Finally, one has to check that the equation between these invertible modifications given in definition 4.1 of \cite{Gur} are satisfied. All of them follow readily from the 2-universal property of the relative tensor product over either three or four algebras.
\end{proof}

\begin{Remark}
Over a perfect field, the 3-category $\mathbf{Mor}^{sep}(\mathbf{2Vect})$ constructed above is the underlying 3-category of the symmetric monoidal 3-category $\mathbf{TC}^{sep}$ of separable multifusion 1-categories considered in \cite{DSPS13}. Over an algebraically closed field of characteristic zero, and given $\mathcal{B}$ a braided fusion 1-category, the 3-category $\mathbf{Mor}^{sep}(\mathbf{Mod}(\mathcal{B}))$ corresponds to the $Hom$-3-category from $\mathcal{B}$ to $\mathbf{Vect}$ in the symmetric monoidal 4-category $\mathbf{BrFus}$ of braided fusion 1-categories considered in \cite{BJS}.
\end{Remark}

\begin{Remark}
Let $\mathfrak{C}$ be a Karoubi complete monoidal 2-category. In \cite{GJF}, the authors outlined the construction of a 3-category $Kar(\mathrm{B}\mathfrak{C})$ of 3-condensation monads, condensation bimodules, and their morphisms. Using variants of the results proven in section 3 of \cite{GJF}, we expect that one can prove that the 3-category $\mathbf{Mor}^{sep}(\mathfrak{C})$ considered above is equivalent to $Kar(\mathrm{B}\mathfrak{C})$. In particular, this would show that $\mathbf{Mor}^{sep}(\mathfrak{C})$ satisfies a 4-universal property. (We refer the reader to \cite{D1} for a precise discussion of the 2-categorical case.)
\end{Remark}

\begin{Remark}
Our proof of theorem \ref{thm:separablealgebra3category} also applies to other setups. Namely, given any monoidal 2-category $\mathfrak{C}$ and any set $\mathscr{A}$ of algebras in $\mathfrak{C}$ such that for any algebras $A$, $B$, and $C$ in $\mathscr{A}$ the relative tensor product over $B$ of any $A$-$B$-bimodule and $B$-$C$-bimodule exists. The above proof constructs a 3-category $\mathbf{Mor}^{\mathscr{A}}(\mathfrak{C})$ of algebras in $\mathscr{A}$, bimodules between them, and their bimodule morphisms. In particular, if every codescent diagram admits a pseudo-coequalizer in $\mathfrak{C}$, and that $\Box$ commutes with them, then it follows from remark \ref{rem:relativetensorcodescent} that the relative tensor product over any algebra in $\mathfrak{C}$ exists. In this case, we can therefore consider the 3-category $\mathbf{Mor}(\mathfrak{C})$ of all algebras in $\mathfrak{C}$, bimodules and their bimodule morphims. We note that this last example has already been thoroughly examined in \cite{H} in an $\infty$-categorical context.
\end{Remark}

For later use, let us also record the following corollary.

\begin{Corollary}\label{cor:algebrabimodule}
Let $A$ be an algebra in $\mathfrak{C}$. Giving an algebra $B$ in the monoidal 2-category $\mathbf{Bimod}_{\mathfrak{C}}(A)$ is equivalent to giving an algebra $B$ in $\mathfrak{C}$ together with an algebra 1-homomorphism $A\rightarrow B$. Furthermore, if $\mathfrak{C}$ is compact semisimple, then the algebra $B$ in $\mathbf{Bimod}_{\mathfrak{C}}(A)$ is rigid, respectively separable, if and only if the underlying algebra $B$ in $\mathfrak{C}$ is rigid, respectively separable.
\end{Corollary}
\begin{proof}
Inspection of the proof of theorem \ref{thm:separablealgebra3category} shows that the forgetful 2-functor $\mathbf{Bimod}_{\mathfrak{C}}(A)\rightarrow \mathfrak{C}$ is lax monoidal. This yields the forward direction of the first part. The backward direction follows by the 2-universal property of the balanced tensor product over $A$. For the second, note that, if we write $\mathfrak{D}:=\mathbf{Bimod}_{\mathfrak{C}}(A)$, then it follows 2-universal property of the relative tensor product over $A$ that the forgetful 2-functor $\mathbf{Bimod}_{\mathfrak{D}}(B)\rightarrow \mathbf{Bimod}_{\mathfrak{C}}(B)$ is an equivalence. The result then follows from theorems 2.2.8 and 3.1.6 of \cite{D7}.
\end{proof}

%% file: Module2Categories.tex
\section{Module 2-Categories}\label{sec:M2C}

We recall the definitions of a module 2-category, module 2-functor, module 2-natural transformation, and module modification and show that, over a fixed monoidal 2-category, these objects assemble into a 3-category. We then review the definition of a 2-adjunction between two 2-functors, and explain how this concepts interacts with that of a module 2-functor over a rigid monoidal 2-category. Theses results are quite technical in nature, but will play a determining role in the last part of the present article.

\subsection{The 3-Category of Module 2-Categories}

Let $\mathfrak{C}$ be a cubical monoidal 2-category.  Our goal is to construct a 3-category whose objects are left $\mathfrak{C}$-module 2-categories in the sense of definition 2.1.3 of \cite{D4}. Now, it follows from proposition 2.2.8 of \cite{D4} that every pair $(\mathfrak{C},\mathfrak{M})$ consisting of a monoidal 2-category $\mathfrak{C}$ and a left $\mathfrak{C}$-module 2-category $\mathfrak{M}$ is equivalent to a pair in which both $\mathfrak{C}$ and $\mathfrak{M}$ are strict cubical (see definition \ref{def:module2cat} below). Thus, there is no loss of generality in assuming that $\mathfrak{C}$ and $\mathfrak{M}$ are strict cubical. In fact, by remark 2.2.9 of \cite{D4}, this strictification procedure holds for any set of module 2-categories.

\begin{Definition} \label{def:module2cat}
Let $\mathfrak{M}$ be a strict 2-category. A strict cubical left $\mathfrak{C}$-module 2-category structure on $\mathfrak{M}$ is a strict cubical 2-functor $\Box:\mathfrak{C}\times \mathfrak{M}\rightarrow \mathfrak{M}$ such that:
\begin{enumerate}
    \item The induced 2-functor $I\Box (-):\mathfrak{M}\rightarrow \mathfrak{M}$ is exactly the identity 2-functor,
    \item The two 2-functors $$\big((-)\Box (-)\big)\Box (-):\mathfrak{C}\times \mathfrak{C}\times\mathfrak{M}\rightarrow\mathfrak{M},\ \mathrm{and}\ (-)\Box \big((-)\Box (-)\big):\mathfrak{C}\times \mathfrak{C}\times\mathfrak{M}\rightarrow\mathfrak{M}$$ are equal on the nose.
\end{enumerate}
\end{Definition}

\begin{Notation}
It is straightforward to extend the graphical conventions introduced in \ref{sub:calculus} for strict cubical monoidal 2-categories to strict cubical left $\mathfrak{C}$-module 2-categories. Throughout this section, we use this extended graphical language.
\end{Notation}

\begin{Remark}
If $\mathds{k}$ is a field, and $\mathfrak{C}$ is a monoidal $\mathds{k}$-linear 2-category, then, by definition, $\Box:\mathfrak{C}\Box\mathfrak{C}\rightarrow \mathfrak{C}$ is a bilinear 2-functor. Likewise, if $\mathfrak{M}$ is a $\mathds{k}$-linear 2-category left $\mathfrak{C}$-module 2-category, we require that $\Box:\mathfrak{C}\times\mathfrak{M}\rightarrow \mathfrak{M}$ is a bilinear 2-functor.
\end{Remark}

\begin{Definition}\label{def:module2fun}
Let $\mathfrak{M}$ and $\mathfrak{N}$ be two strict cubical left $\mathfrak{C}$-module 2-categories. A left $\mathfrak{C}$-module 2-functor is a (not necessarily strict) 2-functor $F:\mathfrak{M}\rightarrow \mathfrak{N}$ together with:
\begin{enumerate}
    \item An adjoint 2-natural equivalence $k^F$ given on $A$ in $\mathfrak{C}$, and $M$ in $\mathfrak{M}$ by $$k^F_{A,M}:A\Box F(M)\rightarrow F(A\Box M);$$
    \item Two invertible modifications $\omega^F$, and $\gamma^F$ given on $A,B$ in $\mathfrak{C}$ and $M$ in $\mathfrak{M}$ by 
\end{enumerate}

\begin{center}
 $$\begin{tikzcd}
                                                                                          & A\Box F(B\Box M) \arrow[rd, "{k^F_{A,B\Box M}}"] \arrow[d, Rightarrow, "\omega^F_{A,B,M}", shorten > = 1ex, shorten < = 1ex] &                    \\
A\Box B\Box F(M) \arrow[rr, "{k^F_{A\Box B,M}}"'] \arrow[ru, "{Id_A\Box k^F_{B,M}}"] & {}                                                            & F(A\Box B\Box M),
\end{tikzcd}$$

$$\gamma^F_M: k^F_{I,M}\Rightarrow Id_{F(M)};$$
\end{center}
Subject to the following relations:
\begin{enumerate}
    \item[a.] For every $A,B,C$ in $\mathfrak{C}$, and $M$ in $\mathfrak{M}$, the equality
\end{enumerate}

\settoheight{\prelim}{\includegraphics[width=45mm]{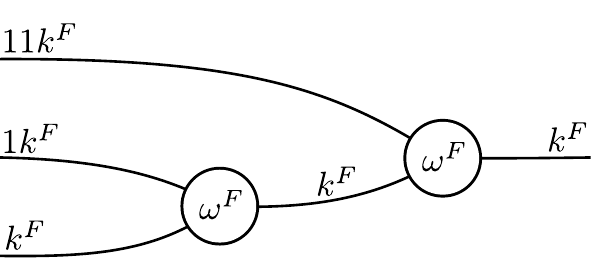}}

\begin{center}
\begin{tabular}{@{}ccc@{}}

\includegraphics[width=45mm]{Pictures/mod2fun/module2funpent1.pdf} & \raisebox{0.45\prelim}{$=$} &

\includegraphics[width=45mm]{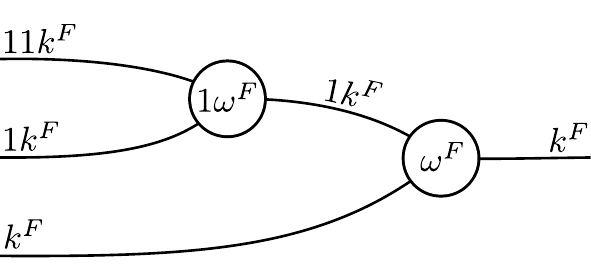}

\end{tabular}
\end{center}

\begin{enumerate}
    \item[] holds in $Hom_{\mathfrak{N}}(A\Box B\Box C\Box F(M), F(A\Box B\Box C\Box M))$,

    \item[b.] For every $A$ in $\mathfrak{C}$, and $M$ in $\mathfrak{M}$, the equality
    
\settoheight{\prelim}{\includegraphics[width=22.5mm]{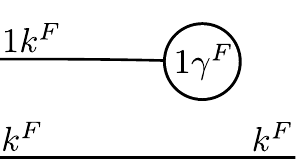}}

\begin{center}
\begin{tabular}{@{}ccc@{}}

\includegraphics[width=22.5mm]{Pictures/mod2fun/module2fununit1.pdf} & \raisebox{0.45\prelim}{$=$} &

\includegraphics[width=22.5mm]{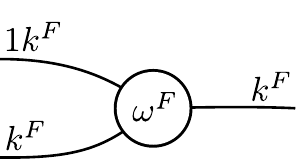}

\end{tabular}
\end{center}
    
    holds in $Hom_{\mathfrak{N}}(A\Box I\Box F(M), F(A\Box M))$;
    
    \item[c.] For every $B$ in $\mathfrak{C}$, and $M$ in $\mathfrak{M}$, the equality

\settoheight{\prelim}{\includegraphics[width=22.5mm]{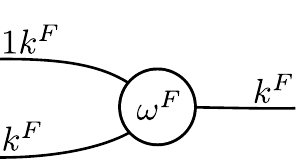}}

\begin{center}
\begin{tabular}{@{}ccc@{}}

\includegraphics[width=22.5mm]{Pictures/mod2fun/module2fununit3.pdf} & \raisebox{0.45\prelim}{$=$} &

\includegraphics[width=22.5mm]{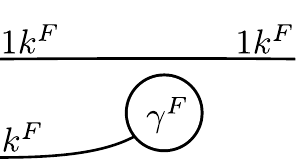}

\end{tabular}
\end{center}

holds in $Hom_{\mathfrak{N}}(I\Box B\Box F(M), F(B\Box M)).$

\end{enumerate}
\end{Definition}

\begin{Definition}\label{def:module2nat}
Let $F,G:\mathfrak{M}\rightarrow \mathfrak{N}$ be two left $\mathfrak{C}$-module 2-functors as in definition \ref{def:module2fun}. A left $\mathfrak{C}$-module 2-natural transformation is a 2-natural transformation $\theta:F\Rightarrow G$ equipped with an invertible modification $\Pi^{\theta}$ given on $A$ in $\mathfrak{C}$, and $M$ in $\mathfrak{M}$ by $$\begin{tikzcd}[sep=tiny]
A G(M) \arrow[ddd, "k^G"']\arrow[Rightarrow, rrrddd, "\Pi^{\theta}", shorten > = 2ex, shorten < = 2ex]  &                                        &    & A F(M) \arrow[ddd, "k^F"] \arrow[lll, "1 \theta"'] \\
 &  &    & \\
  &  &  &  \\
G(AM)                                            &                                        &    &  F(AM); \arrow[lll, "\theta"]
\end{tikzcd}$$

Subject to the following relations:

\begin{enumerate}
    \item[a.] For every $A,B$ in $\mathfrak{C}$, and $M$ in $\mathfrak{M}$, the equality
\end{enumerate}

\settoheight{\prelim}{\includegraphics[width=52.5mm]{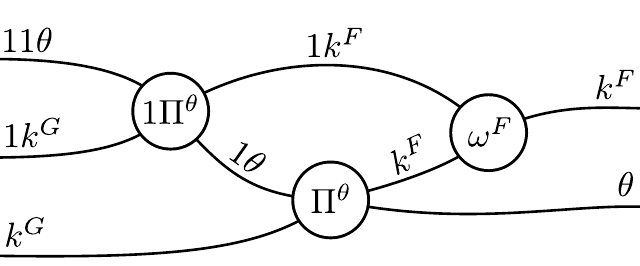}}

\begin{center}
\begin{tabular}{@{}ccc@{}}

\includegraphics[width=52.5mm]{Pictures/mod2fun/natmod/module2nat1.pdf} & \raisebox{0.45\prelim}{$=$} &

\includegraphics[width=45mm]{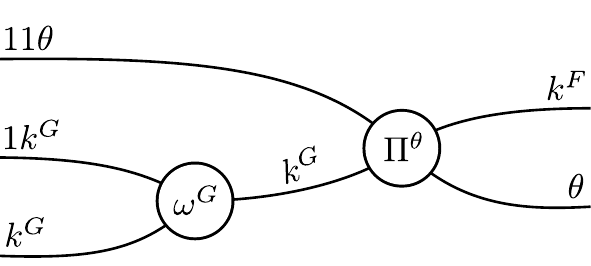}

\end{tabular}
\end{center}

\begin{enumerate}
    
    \item[] holds in $Hom_{\mathfrak{N}}(A\Box B\Box F(M), G(A\Box B\Box M))$;
    
    \item[b.] For every $M$ in $\mathfrak{M}$, the equality

\settoheight{\prelim}{\includegraphics[width=37.5mm]{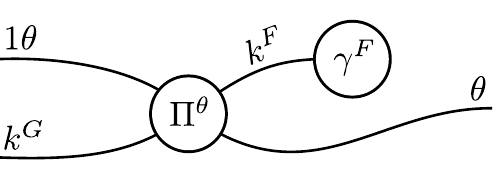}}

\begin{center}
\begin{tabular}{@{}ccc@{}}

\includegraphics[width=37.5mm]{Pictures/mod2fun/natmod/module2nat3.pdf} & \raisebox{0.45\prelim}{$=$} &

\includegraphics[width=30mm]{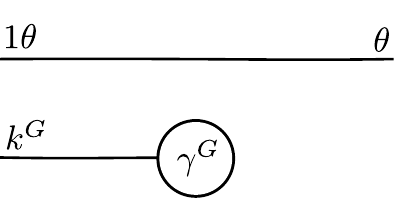}

\end{tabular}
\end{center}

holds in $Hom_{\mathfrak{N}}(I\Box F(M), G(M)).$

\end{enumerate}
\end{Definition}

\begin{Definition}\label{def:module2modif}
Let $\theta,\tau:F\Rightarrow G$ be two left $\mathfrak{C}$-module 2-natural transformations. A left $\mathfrak{C}$-module modification is a modification $\Xi:\theta\Rrightarrow \tau$ such that for every $A$ in $\mathfrak{C}$, and $M$ in $\mathfrak{M}$ the equality

\settoheight{\prelim}{\includegraphics[width=30mm]{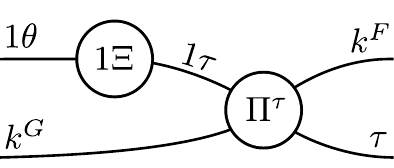}}

\begin{center}
\begin{tabular}{@{}ccc@{}}

\includegraphics[width=30mm]{Pictures/mod2fun/natmod/modulemodif1.pdf} & \raisebox{0.4\prelim}{$=$} &

\includegraphics[width=30mm]{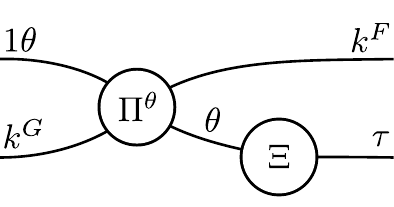}

\end{tabular}
\end{center}

holds in $Hom_{\mathfrak{N}}(A\Box F(M), G(A\Box M)).$
\end{Definition}

Fixing two strict cubical left $\mathfrak{C}$-module 2-categories $\mathfrak{M}$ and $\mathfrak{N}$, it was shown in proposition 2.2.1 of \cite{D4} that left $\mathfrak{C}$-module 2-functors, left $\mathfrak{C}$-module 2-natural transformation, and left $\mathfrak{C}$-module modifications form a 2-category, which we denote by $\mathbf{Fun}_{\mathfrak{C}}(\mathfrak{M},\mathfrak{N})$. In particular, given $\theta,\overline{\theta},\overline{\overline{\theta}}:F\Rightarrow G$ left $\mathfrak{C}$-module 2-natural transformations, and two left $\mathfrak{C}$-module modifications $\Xi:\theta\Rrightarrow\overline{\theta}$, $Z:\overline{\theta}\Rrightarrow\overline{\overline{\theta}}$, the vertical composite $Z\bullet \Xi$ is a left $\mathfrak{C}$-module modification. Further, given two left $\mathfrak{C}$-module 2-natural transformations $\theta:F\Rightarrow G$ and $\tau:G\Rightarrow H$, their composite is the 2-natural transformation $\tau \cdot \theta$ equipped with the invertible modification $$\Pi^{\tau \cdot \theta}:= (\Pi^{\tau}\cdot \theta)\bullet (\tau \cdot \Pi^{\theta}).$$ Thanks to our strictness hypotheses, the above composition of left $\mathfrak{C}$-module 2-natural transformations is in fact strictly associative and unital. Thence, $\mathbf{Fun}_{\mathfrak{C}}(\mathfrak{M},\mathfrak{N})$ is in fact a strict 2-category. For later use, we now assemble all of these 2-categories together.

\begin{Theorem}\label{thm:leftCmodule3category}
Let $\mathfrak{C}$ be a monoidal 2-category. Left $\mathfrak{C}$-module 2-categories, left $\mathfrak{C}$-module 2-functors, left $\mathfrak{C}$-module 2-natural transformations, and left $\mathfrak{C}$-module modifications form a 3-category, which we denote by $\mathbf{LMod}(\mathfrak{C})$.
\end{Theorem}
\begin{proof}
In section 5.1 of \cite{Gur}, the author constructs a 3-category of 2-categories. Our proof is precisely a left $\mathfrak{C}$-module version of this argument. In order to do this, it is enough to upgrade the structures defined in section 5.1 of \cite{Gur} with suitable left $\mathfrak{C}$-module actions. Furthermore, thanks to proposition 2.2.8 and remark 2.2.9 of \cite{D4}, we may assume without loss of generality that $\mathfrak{C}$ and every left $\mathfrak{C}$-module 2-category is strict cubical.

Let $\mathfrak{M}$, $\mathfrak{N}$, and $\mathfrak{P}$ be strict cubical left $\mathfrak{C}$-module 2-categories. We begin by constructing the 2-functor $$\circ:\mathbf{Fun}_{\mathfrak{C}}(\mathfrak{N},\mathfrak{P})\times \mathbf{Fun}_{\mathfrak{C}}(\mathfrak{M},\mathfrak{N})\rightarrow \mathbf{Fun}_{\mathfrak{C}}(\mathfrak{M},\mathfrak{P})$$ providing us with the composition of left $\mathfrak{C}$-module 2-functors. Given two left $\mathfrak{C}$-module 2-functors $F:\mathfrak{M}\rightarrow \mathfrak{N}$ and $G:\mathfrak{N}\rightarrow \mathfrak{P}$, we endow their composite $G\circ F$ with a left $\mathfrak{C}$-module structure using the following assignments. We define the adjoint 2-natural equivalence $k^{G\circ F}$ by $$k^{G\circ F}_{A,M} := G(k^F_{A,M})\circ k^G_{A,F(M)},$$ for every $A$ in $\mathfrak{C}$ and $M$ in $\mathfrak{M}$, and the two invertible modifications $\omega^{G\circ F}$, and $\gamma^{G\circ F}$ by 

\begin{center}
\settoheight{\prelim}{\includegraphics[width=45mm]{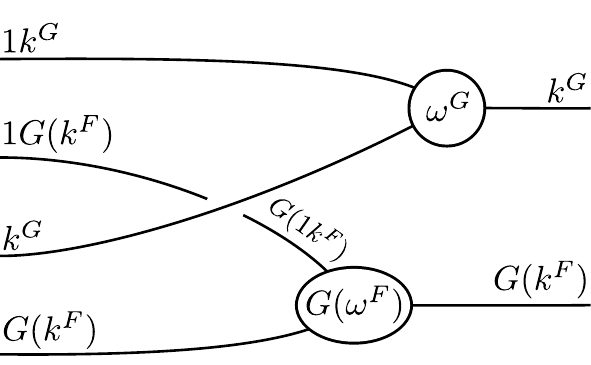}}

\begin{tabular}{@{}cc@{}}

 \raisebox{0.45\prelim}{$\omega^{G\circ F}:=$} &

\includegraphics[width=45mm]{Pictures/mod2fun/natmod/omegacomposite.pdf},

\end{tabular}

\settoheight{\prelim}{\includegraphics[width=30mm]{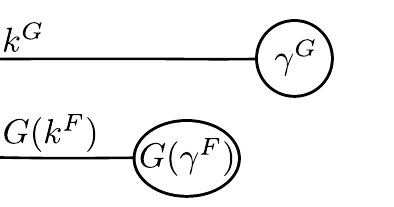}}

\begin{tabular}{@{}cc@{}}

 \raisebox{0.45\prelim}{$\gamma^{G\circ F}:=$} &

\includegraphics[width=30mm]{Pictures/mod2fun/natmod/gammacomposite.pdf}.

\end{tabular}
\end{center}

\noindent It is not hard to show that the above data satisfies the axioms of definition \ref{def:module2fun}.

Then, given a left $\mathfrak{C}$-module 2-natural transformation $\theta:F_1\Rightarrow F_2$ between two left $\mathfrak{C}$-module 2-functors $F_1,F_2:\mathfrak{M}\rightarrow \mathfrak{N}$, we endow the 2-natural transformation $G\circ \theta$ with a left $\mathfrak{C}$-module structure by setting

\begin{center}

\settoheight{\prelim}{\includegraphics[width=45mm]{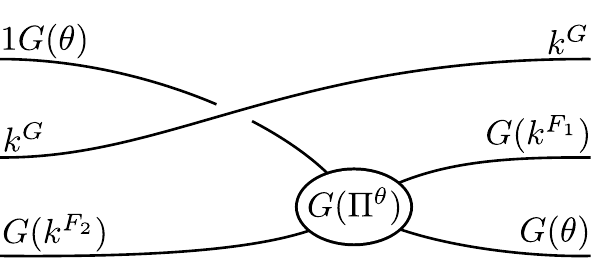}}

\begin{tabular}{@{}cc@{}}

\raisebox{0.45\prelim}{$\Pi^{G\circ\theta}:=$} &

\includegraphics[width=45mm]{Pictures/mod2fun/natmod/thetacomposite.pdf}.

\end{tabular}
\end{center}

Likewise, given a left $\mathfrak{C}$-module 2-natural transformation $\tau:G_1\Rightarrow G_2$ between two left $\mathfrak{C}$-module 2-functors $G_1,G_2:\mathfrak{N}\rightarrow \mathfrak{P}$, we may similarly define a left $\mathfrak{C}$-module structure on the 2-natural transformation $\tau\circ F$ by 

\begin{center}

\settoheight{\prelim}{\includegraphics[width=45mm]{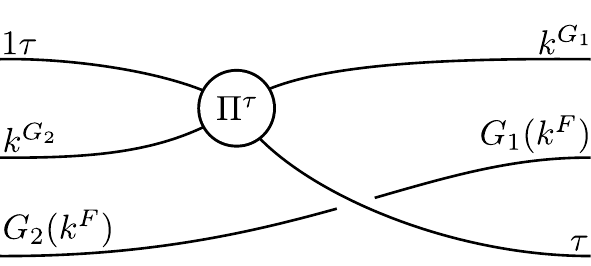}}

\begin{tabular}{@{}cc@{}}

\raisebox{0.45\prelim}{$\Pi^{\tau\circ F}:=$} &

\includegraphics[width=45mm]{Pictures/mod2fun/natmod/taucomposite.pdf}.

\end{tabular}
\end{center}

\noindent Now, recall from the proof of proposition 5.1 of \cite{Gur} that $\tau\circ\theta = (G_2\circ\theta)\cdot (\tau\circ F_1)$, so that the 2-natural transformation $\tau\circ\theta$ inherits a $\mathfrak{C}$-module structure. These assignments can be straightforwardly extended to left $\mathfrak{C}$-module modifications, so that we obtain a functor $$Nat_{\mathfrak{C}}(G_1,G_2)\times Nat_{\mathfrak{C}}(F_1,F_2)\rightarrow Nat_{\mathfrak{C}}(G_1\circ F_1,G_2\circ F_2)$$ between the 1-categories of left $\mathfrak{C}$-module 2-natural transformations. The additional structure constraints needed to define a 2-functor are the invertible modifications given in proposition 5.1 of \cite{Gur}, and one checks easily that they respect the relevant left $\mathfrak{C}$-module structures defined above. Thus, we obtain the desired 2-functor $\circ:\mathbf{Fun}_{\mathfrak{C}}(\mathfrak{N},\mathfrak{P})\times \mathbf{Fun}_{\mathfrak{C}}(\mathfrak{M},\mathfrak{N})\rightarrow \mathbf{Fun}_{\mathfrak{C}}(\mathfrak{M},\mathfrak{P})$. Moreover, the unit on $\mathfrak{M}$ for the composition of left $\mathfrak{C}$-module 2-functors is given by the identity 2-functor $Id:\mathfrak{M}\rightarrow \mathfrak{M}$ with its canonical left $\mathfrak{C}$-module structure.

Proposition 5.3 of \cite{Gur} defines an adjoint 2-natural equivalence $\alpha$ witnessing the associativity of the composition of (plain) 2-functors. Now, let $\mathfrak{M}$, $\mathfrak{N}$, $\mathfrak{P}$, and $\mathfrak{Q}$ be strict cubical left $\mathfrak{C}$-module 2-categories, and let $F:\mathfrak{M}\rightarrow\mathfrak{N}$, $G:\mathfrak{N}\rightarrow\mathfrak{P}$, and $H:\mathfrak{P}\rightarrow\mathfrak{Q}$ be left $\mathfrak{C}$-module 2-functors. It follows from proposition 5.3 of \cite{Gur} that $\alpha_{H,G,F}:(H\circ G)\circ F\simeq H\circ (G\circ F)$ is the identity 2-natural transformation. Further, the left $\mathfrak{C}$-module structures of $(H\circ G)\circ F$ and $H\circ (G\circ F)$ are equal, so that we can upgrade $\alpha_{H,G,F}$ to a left $\mathfrak{C}$-module adjoint 2-natural equivalence using the identity modification. The collection of these assignments promote $\alpha$ to an adjoint 2-natural equivalence witnessing the associativity of the composition of left $\mathfrak{C}$-module 2-functors.

Analogously, the adjoint 2-natural equivalences $l$, and $r$ constructed in proposition 5.5 of \cite{Gur}, witnessing that composition of 2-functors is unital, can be promoted to left $\mathfrak{C}$-module adjoint 2-natural equivalences. Namely, as we are working with strict 2-categories, these adjoint 2-natural equivalences are in fact both given by the identity 2-natural adjoint equivalence. Thus, given a $\mathfrak{C}$-module 2-functor $F:\mathfrak{M}\rightarrow \mathfrak{N}$ between strict cubical left $\mathfrak{C}$-module 2-categories, the 2-natural transformations $l_F$ and $r_F$ can canonically be upgraded to left $\mathfrak{C}$-module adjoint 2-natural equivalences. With these additional pieces of data, $l$ and $r$ define adjoint 2-natural equivalence witnessing the unitality of the composition of left $\mathfrak{C}$-module 2-functors. The proof is then completed by checking that the invertible modification $\pi$, $\mu$, $\lambda$, and $\rho$ given in proposition 5.6 of \cite{Gur} are compatible with the left $\mathfrak{C}$-module structures we have defined. This is immediate as it follows from our strictness assumptions that $\pi$, $\mu$, $\lambda$, and $\rho$ are all identity modifications.
\end{proof}

\begin{Remark}
It follows immediately from our proof of theorem \ref{thm:leftCmodule3category} that there is a forgetful 3-functor $\mathbf{LMod}(\mathfrak{C})\rightarrow \mathbf{2Cat}$ to the 3-category of 2-categories.
\end{Remark}

\begin{Corollary}\label{cor:bimodule2functorcategory}
Let $\mathfrak{M}$ be a left $\mathfrak{C}$-module 2-category. Then, $\mathbf{End}_{\mathfrak{C}}(\mathfrak{M})$ is a monoidal 2-category. Further, given $\mathfrak{M}$ be any left $\mathfrak{C}$-module 2-category, the 2-category $\mathbf{Fun}_{\mathfrak{C}}(\mathfrak{M},\mathfrak{N})$ is an $\mathbf{End}_{\mathfrak{C}}(\mathfrak{N})$-$\mathbf{End}_{\mathfrak{C}}(\mathfrak{M})$-bimodule 2-category.
\end{Corollary}

We end this first section on module 2-categories with the following proposition, which will constitute a key ingredient in our study of the Morita theory of fusion 2-categories.

\begin{Proposition}\label{prop:dualmoduleaction}
Let $\mathfrak{M}$ be a left $\mathfrak{C}$-module 2-category. The action of $\mathbf{End}_{\mathfrak{C}}(\mathfrak{M})$ on $\mathfrak{M}$ given by evaluation defines a left $\mathbf{End}_{\mathfrak{C}}(\mathfrak{M})$-module structure on $\mathfrak{M}$. Further, this structure is compatible with the left $\mathfrak{C}$-module one, so that $\mathfrak{M}$ is a left $\mathbf{End}_{\mathfrak{C}}(\mathfrak{M})\times\mathfrak{C}$-module 2-category.
\end{Proposition}
\begin{proof}
One may directly check that evaluation of 2-functors $\mathbf{End}_{\mathfrak{C}}(\mathfrak{M})\times\mathfrak{M}\rightarrow \mathfrak{M}$ provides $\mathfrak{M}$ with a left $\mathbf{End}_{\mathfrak{C}}(\mathfrak{M})$-module structure. By definition, this left $\mathbf{End}_{\mathfrak{C}}(\mathfrak{M})$-module structure on $\mathfrak{M}$ is compatible with the left $\mathfrak{C}$-module structure, so that $\mathfrak{M}$ has a left $\mathbf{End}_{\mathfrak{C}}(\mathfrak{M})\times\mathfrak{C}$-module structure.
\end{proof}

\subsection{Module 2-Functors and 2-Adjunctions}\label{sub:2adjunction}

The goal of this section is study the interaction between the notion of a module 2-functor recalled above, and that of a 2-adjunction between 2-functors, which we now recall.

\begin{Definition}\label{def:2adjunction}
Let $\mathfrak{M}$ and $\mathfrak{N}$ be two 2-categories, and $F:\mathfrak{M}\rightarrow \mathfrak{N}$ and $G:\mathfrak{N}\rightarrow \mathfrak{M}$ be two 2-functors. A 2-adjunction between $F$ and $G$ consists of two 2-natural transformations $u^F$, called the unit, and $c^F$, called the counit, given on $M$ in $\mathfrak{M}$ and $N$ in $\mathfrak{N}$ by $$u^F_M:M\rightarrow G(F(M)),\textrm{  and  }c^F_N:F(G(N))\rightarrow N,$$ together with two invertible modifications $\Phi^F$ and $\Psi^F$, called triangulators, given on $M$ in $\mathfrak{M}$ and $N$ in $\mathfrak{N}$ by $$\Phi^F_{M}:c^F_{F(M)}\circ F(u^F_M)\cong Id_{F(M)},$$ $$\Psi^F_N:G(c^F_N)\circ u^F_{G(N)}\cong Id_{G(N)}.$$ We also say that $F$ is a left 2-adjoint to $G$, or that $G$ is a right 2-adjoint to $F$.
\end{Definition}

\begin{Definition}
Let $\mathfrak{C}$ be a monoidal 2-category, $\mathfrak{M}$ and $\mathfrak{N}$ be two left $\mathfrak{C}$-module 2-categories, and $F:\mathfrak{M}\rightarrow \mathfrak{N}$ and $G:\mathfrak{N}\rightarrow \mathfrak{M}$ be two left $\mathfrak{C}$-module 2-functors. A left $\mathfrak{C}$-module 2-adjunction between $F$ and $G$ is a 2-adjunction between $F$ and $G$ as in definition \ref{def:2adjunction} such that $u^F$ and $c^F$ are left $\mathfrak{C}$-module 2-natural transformations, and $\Phi^F$ and $\Psi^F$ are left $\mathfrak{C}$-module modifications.
\end{Definition}

Let $\mathfrak{C}$ be a rigid monoidal 2-category, and assume that both $\mathfrak{M}$ and $\mathfrak{N}$ are left $\mathfrak{C}$-module 2-categories. The categorified version of corollary 2.13 of \cite{DSPS14} holds, as we show in the next two propositions. In fact, our proof also establishes the categorifications of their lemmas 2.10 and 2.11.

\begin{Proposition}\label{prop:C2Funrigid}
Let $\mathfrak{C}$ be a rigid monoidal 2-category, and let $F:\mathfrak{M}\rightarrow\mathfrak{N}$ be a left $\mathfrak{C}$-module 2-functor between left $\mathfrak{C}$-module 2-categories. If $F$ has a right 2-adjoint $G$, then $G$ can canonically be upgraded to a left $\mathfrak{C}$-module right 2-adjoint to $F$.
\end{Proposition}
\begin{proof}
Thank to proposition 2.2.8 and remark 2.2.9 of \cite{D4}, we may assume that $\mathfrak{C}$ is strict cubical, and that both $\mathfrak{M}$ and $\mathfrak{N}$ are strict cubical left $\mathfrak{C}$-module 2-categories. Further, let us denote by $(\omega^{F^{-1}})^{\bullet}$ and $(\gamma^{F^{-1}})^{\bullet}$ the 2-isomorphisms given by

\settoheight{\prelim}{\includegraphics[width=30mm]{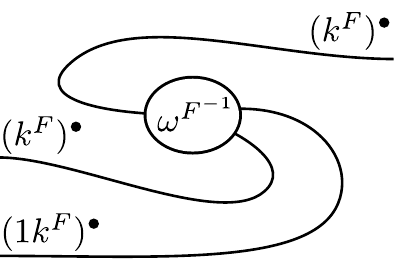}}

\begin{center}
\raisebox{0.45\prelim}{$(\omega^{F^{-1}})^{\bullet}:=\ $}
\includegraphics[width=30mm]{Pictures/mod2fun/2adjunction/omegaFbullet.pdf}, \raisebox{0.45\prelim}{$\ (\gamma^{F^{-1}})^{\bullet}:=\ $}\includegraphics[width=22.5mm]{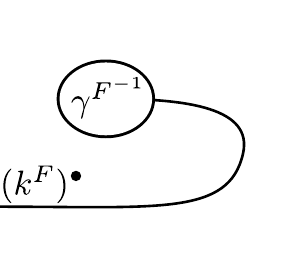},
\end{center}

\noindent where the cups and the caps denote the unit and counit 2-isomorphisms witnessing that $k^F$ and $(k^F)^{\bullet}$ form an adjoint 2-natural equivalence.

We begin by proving that $G$ can be endowed with a lax left $\mathfrak{C}$-module structure. Given $A$ in $\mathfrak{C}$, and $M$ in $\mathfrak{M}$, we let the 2-natural transformation $k^G$ be given by $$k^G_{A,M}:A\Box G(M)\xrightarrow{u^F} G(F(A\Box G(M)))\xrightarrow{G((k^F)^{\bullet})} G(A\Box F(G(M)))\xrightarrow{G(1c^F)} G(A\Box M),$$ where $(k^F)^{\bullet}$ is the pseudo-inverse of $k^F$ provided in the data of a left $\mathfrak{C}$-module 2-functor. The invertible modifications $\omega^G$ and $\gamma^G$ are given by

\settoheight{\prelim}{\includegraphics[width=90mm]{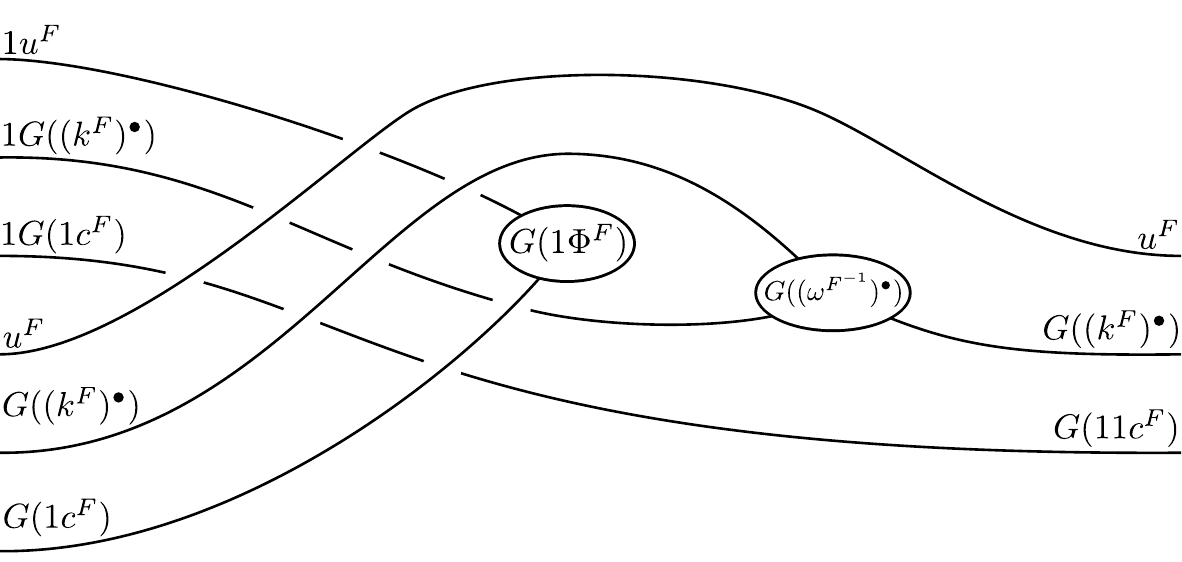}}

\begin{center}
\begin{tabular}{@{}cc@{}}
\raisebox{0.45\prelim}{$\omega^G:=$} &
\includegraphics[width=90mm]{Pictures/mod2fun/2adjunction/omegaG.pdf},
\end{tabular}
\end{center}

\settoheight{\prelim}{\includegraphics[width=37.5mm]{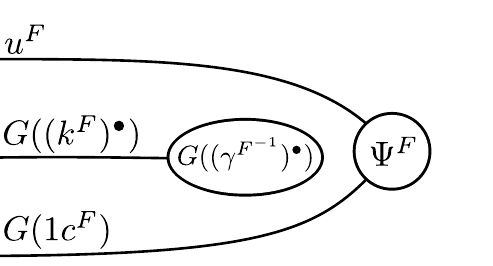}}

\begin{center}
\begin{tabular}{@{}cc@{}}
\raisebox{0.45\prelim}{$\gamma^G:=$} &
\includegraphics[width=37.5mm]{Pictures/mod2fun/2adjunction/gammaG.pdf}.
\end{tabular}
\end{center}

\noindent Using the axioms of definition \ref{def:module2fun} for $F$, it is easy to check that $\omega^G$ and $\gamma^G$ satisfy the axioms of \ref{def:module2fun}.

We now show that $k^G$ can be upgraded to an adjoint 2-natural equivalence. As every 2-natural equivalence can be upgraded to an adjoint 2-natural equivalence (see section 1 of \cite{Gur2}), it is enough to exhibit for every $A$ in $\mathfrak{C}$ and $M$ in $\mathfrak{M}$, a pseudo-inverse $(k^G)^{\bullet}_{A,M}$ for the 1-morphism $k^G_{A,M}$. Let $^{\sharp}A$ be a left dual for $A$ in $\mathfrak{C}$ with unit 1-morphism $i_A:I\rightarrow A\Box {^{\sharp}A}$, counit 1-morphism $e_A:{^{\sharp}A}\Box A\rightarrow I$ and 2-isomorphisms $C_A:(e_A\Box {^{\sharp}A})\circ ((^{\sharp}A)\Box i_A)\Rightarrow Id_{^{\sharp}A}$, and $D_A:Id_A\Rightarrow (A\Box e_A)\circ (i_A\Box A)$. We define \begin{align*}(k^G)^{\bullet}_{A,M}:G(A\Box M)&\xrightarrow{i_A1} A\Box (^{\sharp} A)\Box G(A\Box M)\xrightarrow{1u^F} A\Box GF({^{\sharp} A}\Box G(A\Box M))\\ & \xrightarrow{1 G((k^F)^{\bullet})}A\Box G({^{\sharp} A}\Box FG(A\Box M))\\ &\xrightarrow{1G(1c^F)} A\Box G({^{\sharp} A}\Box A\Box M) \xrightarrow{1G(e_A1)} A\Box G(M),\end{align*} where $(k^F)^{\bullet}$ denotes the canonical pseudo-inverse of $k^F$ supplied by the definition of a module 2-functor. The two 2-isomorphisms

\begin{center}
\includegraphics[width=112.5mm]{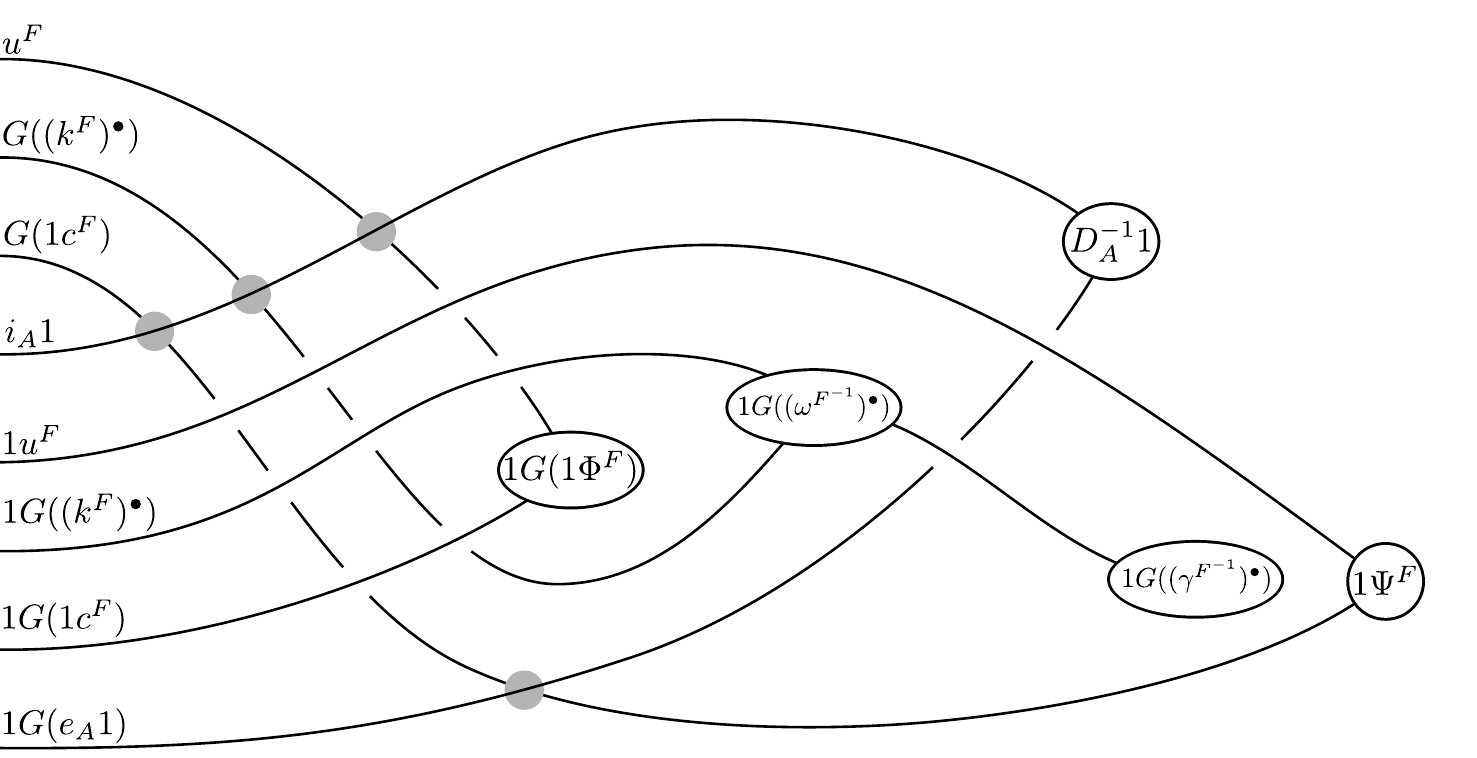},
\includegraphics[width=112.5mm]{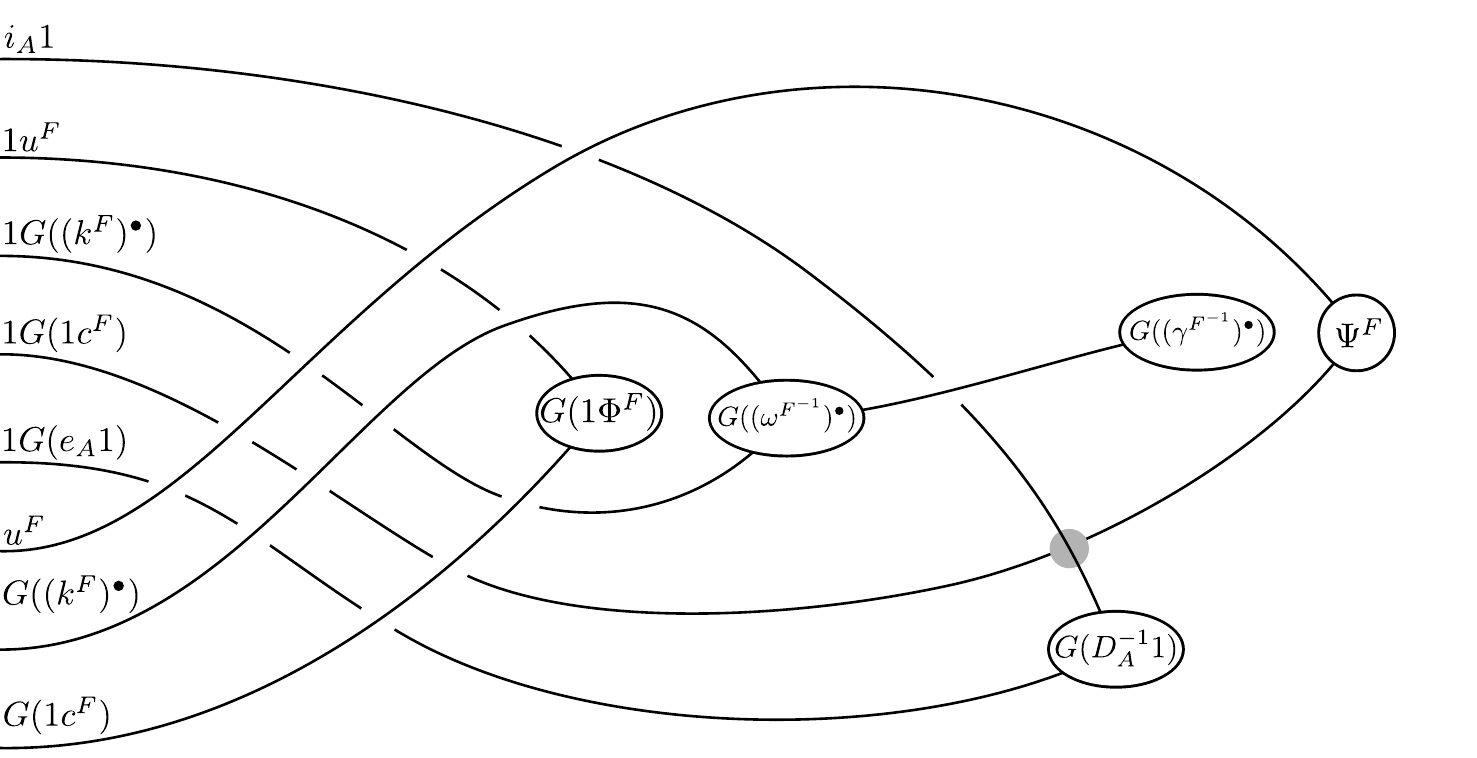},
\end{center}

\noindent witness that $(k^G)^{\bullet}_{A,M}$ is a pseudo-inverse for $k^G_{A,M}$ as desired.

It remains to upgrade $u^F$ and $c^F$ to left $\mathfrak{C}$-module 2-natural transformations, and show that $\Phi^F$ and $\Psi^F$ define invertible left $\mathfrak{C}$-module modifications. For the first part, we endow $u^F$ and $c^F$ with left $\mathfrak{C}$-module structures using the modifications $\Pi^{u^F}$ and $\Pi^{c^F}$ specified by

\settoheight{\prelim}{\includegraphics[width=75mm]{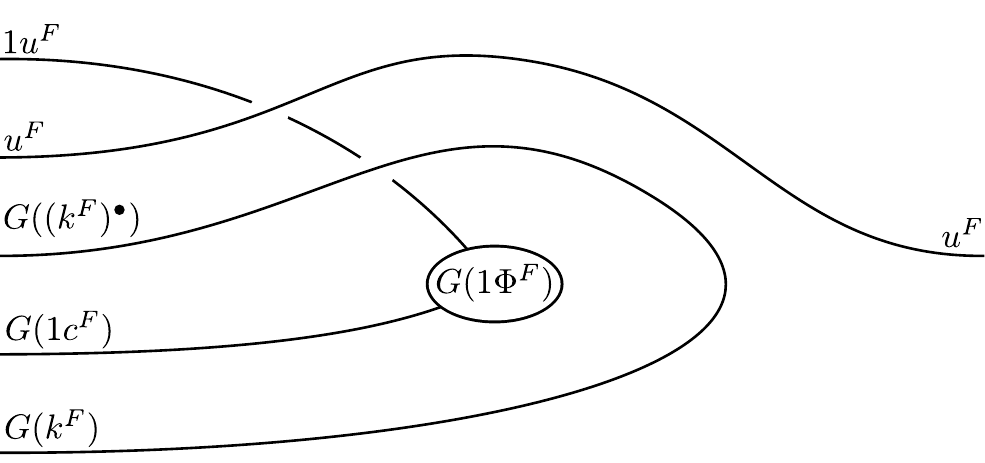}}

\begin{center}
\begin{tabular}{@{}cc@{}}
\raisebox{0.45\prelim}{$\Pi^{u^F}_A:=$} &
\includegraphics[width=75mm]{Pictures/mod2fun/2adjunction/uFmodule.pdf},\\
\raisebox{0.45\prelim}{$\Pi^{c^F}_A:=$} &
\includegraphics[width=75mm]{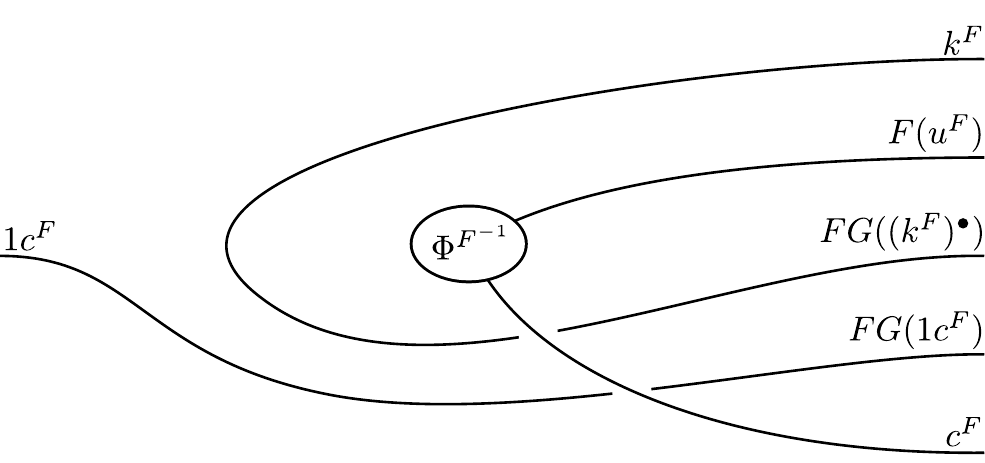}.
\end{tabular}
\end{center}

\noindent It is easy to check that $\Pi^{u^F}$ and $\Pi^{c^F}$ satisfy the axioms of definition \ref{def:module2nat}. Finally, using naturality together with axioms a and b of definition \ref{def:module2fun} for $G$, one can readily check that $\Phi^F$ and $\Psi^F$ are compatible with the left $\mathfrak{C}$-module structure on $u^F$ and $c^F$ defined above, which finishes the proof of the result.
\end{proof}

The analogue of proposition \ref{prop:C2Funrigid} for left 2-adjoints also holds.

\begin{Proposition}
Let $G:\mathfrak{N}\rightarrow\mathfrak{M}$ be a left $\mathfrak{C}$-module 2-functor. If $G$ has a left 2-adjoint $F$, then $F$ can canonically be upgraded to left $\mathfrak{C}$-module left 2-adjoint to $G$.
\end{Proposition}
\begin{proof}
This follows by applying proposition \ref{prop:C2Funrigid} to $\mathfrak{C}^{1op}$, the monoidal 2-category obtained from $\mathfrak{C}$ by reversing the direction of the 1-morphisms.
\end{proof}

%% file: SeparableM2C.tex
\section{Dual Tensor 2-Categories and Morita Equivalence}

In general, it is difficult to work with arbitrary compact semisimple module 2-categories over a fixed compact semisimple tensor 2-category $\mathfrak{C}$. Motivated by the decategorified situation studied in \cite{DSPS13}, we therefore restrict our attention to a particularly nice class of compact semisimple module 2-categories called separable module 2-categories. We prove that the 3-category of separable algebras in $\mathfrak{C}$ is equivalent to the 3-category of separable left $\mathfrak{C}$-module 2-categories. Under the assumption that $\mathfrak{C}$ be locally separable, we then show that the monoidal 2-category of bimodules over a separable algebra in $\mathfrak{C}$ is a compact semisimple tensor 2-category. This allows us to define the dual tensor 2-category of $\mathfrak{C}$ with respect to a separable module 2-category. We end by giving three equivalent characterizations of Morita equivalence between locally separable compact semisimple tensor 2-categories. Throughout, we work over a fixed field $\mathds{k}$, meaning that all (monoidal) categories and functors under consideration are $\mathds{k}$-linear.

\subsection{Separable Module 2-Categories}

Let us fix $\mathfrak{C}$ a compact semisimple tensor 2-category.

\begin{Definition}
A compact semisimple left $\mathfrak{C}$-module 2-category $\mathfrak{M}$ is called separable if there exists a separable algebra $A$ in $\mathfrak{C}$ such that $$\mathfrak{M}\simeq\mathbf{Mod}_{\mathfrak{C}}(A)$$ as left $\mathfrak{C}$-module 2-categories.
\end{Definition}

In theorem \ref{thm:leftCmodule3category}, we have proven that left $\mathfrak{C}$-module 2-categories form a 3-category, which we denote by $\mathbf{LMod}(\mathfrak{C})$. We will write $\mathbf{LMod}^{sep}(\mathfrak{C})$ for the full sub-3-category whose objects are the separable module 2-categories. We are now ready to state our next theorem, of which a closely related variant was conjectured in remark 5.3.9 of \cite{D4}. Let us mention that if $\mathds{k}$ is algebraically closed of characteristic zero and $\mathfrak{C}=\mathbf{2Vect}$, we recover the main result of \cite{D1} as every multifusion 1-category is separable. More generally, if $\mathds{k}$ is perfect and $\mathfrak{C}=\mathbf{2Vect}$, the theorem below also recovers corollary 3.1.5 of \cite{D5} thanks to proposition 2.5.10 of \cite{DSPS13}.

\begin{Theorem}\label{thm:equivalencealgebrasmodules2categories}
Let $\mathfrak{C}$ be a compact semisimple tensor 2-category. There is a linear 3-functor, contravariant on 1-morphisms, $$\mathbf{Mod}_{\mathfrak{C}}:\mathbf{Mor}^{sep}(\mathfrak{C})\rightarrow\mathbf{LMod}^{sep}(\mathfrak{C})$$ that sends a separable algebra in $\mathfrak{C}$ to the associated separable left $\mathfrak{C}$-module 2-category of right modules. Moreover, this 3-functor is an equivalence.
\end{Theorem}
\begin{proof}
Without loss of generality, we may assume that $\mathfrak{C}$ is strict cubical. The monoidal unit $I$ of $\mathfrak{C}$ is canonically a separable algebra in $\mathfrak{C}$. Thanks to theorem \ref{thm:separablealgebra3category}, this yields a contravariant linear 3-functor $$\mathbf{Hom}_{\mathbf{Mor}^{sep}(\mathfrak{C})}(-,I):\mathbf{Mor}^{sep}(\mathfrak{C})\rightarrow\mathbf{2Cat}_{\mathds{k}}$$ to the 3-category of $\mathds{k}$-linear 2-categories. But, we have $\mathbf{End}_{\mathbf{Mor}^{sep}(\mathfrak{C})}(I,I)=\mathbf{Bimod}_{\mathfrak{C}}(I)=\mathfrak{C}$ as monoidal 2-categories thanks to our strictness hypothesis. In particular, for every separable algebra $A$ in $\mathfrak{C}$, the 2-category $\mathbf{Hom}_{\mathbf{Mor}^{sep}(\mathfrak{C})}(A,I)$ has a canonical left $\mathfrak{C}$-module structure, which is compatible with bimodule morphisms in the variable $A$. Further,
$\mathbf{Hom}_{\mathbf{Mor}^{sep}(\mathfrak{C})}(A,I)= \mathbf{Bimod}_{\mathfrak{C}}(I,A)=\mathbf{Mod}_{\mathfrak{C}}(A)$ is a separable left $\mathfrak{C}$-module 2-category. Consequently, the 3-functor $\mathbf{Hom}_{\mathbf{Mor}^{sep}(\mathfrak{C})}(-,I)$ can canonically be lifted to a 3-functor $$\mathbf{Mod}_{\mathfrak{C}}:\mathbf{Mor}^{sep}(\mathfrak{C})\rightarrow\mathbf{LMod}^{sep}(\mathfrak{C}).$$

It remains to prove that $\mathbf{Mod}_{\mathfrak{C}}$ is an equivalence of 3-categories, i.e.\ that it is essentially surjective and induces equivalences on $Hom$-2-categories. Essential surjectivty follows immediately from the definition of a separable left $\mathfrak{C}$-module 2-category. Therefore, it is only left to prove that for every separable algebras $A$, $B$ in $\mathfrak{C}$, the 2-functor $$\begin{tabular}{r c c c}
$\mathbf{F}:$&$\mathbf{Bimod}_{\mathfrak{C}}(A,B)$&$\rightarrow$ &$\mathbf{Fun}_{\mathfrak{C}}(\mathbf{Mod}_{\mathfrak{C}}(A),\mathbf{Mod}_{\mathfrak{C}}(B))$\\& $P$&$\mapsto$& $(-)\Box_A P$\end{tabular}$$ induced by $\mathbf{Mod}_{\mathfrak{C}}$ is an equivalence of 2-categories. In order to exhibit a pseudo-inverse to $\mathbf{F}$, consider the following 2-functor $$\begin{tabular}{r c c c}
   $\mathbf{B}:$ & $\mathbf{Fun}_{\mathfrak{C}}(\mathbf{Mod}_{\mathfrak{C}}(A),\mathbf{Mod}_{\mathfrak{C}}(B))$ & $\rightarrow$ & $\mathbf{Bimod}_{\mathfrak{C}}(A,B),$ \\
   & $F$ & $\mapsto$ & $F(A)$
\end{tabular}$$
where the left $A$-module structure on $F(A)$ arises from the canonical $A$-$A$-bimodule structure on $A$. This assignment can straightforwardly be extended to left $\mathfrak{C}$-module 2-natural transformations and left $\mathfrak{C}$-module modifications. Further, for any $A$-$B$-bimodule $P$ in $\mathfrak{C}$, we have that $\mathbf{B}\circ\mathbf{F}(P) = A\Box_A P$ as an $A$-$B$-bimodule in $\mathfrak{C}$. Thence, by lemma \ref{lem:relativetensorA}, we find that $\mathbf{B}\circ\mathbf{F}\simeq Id$ as desired.

Now, let $F:\mathbf{Mod}_{\mathfrak{C}}(A)\rightarrow\mathbf{Mod}_{\mathfrak{C}}(B)$ be a left $\mathfrak{C}$-module 2-functor. By definition, for any right $A$-module $M$, we have that $\big(\mathbf{F}\circ\mathbf{B}(F)\big)(M) = M\Box_A F(A)$. We claim that $M\Box_A F(A)\simeq F(M)$ as $A$-$B$-bimodules. Namely, as splittings of 2-condensation monads are preserved by all 2-functors, it follows from the last part of theorem \ref{thm:tensormodules} that $F(A\Box M)\rightarrow F(A\Box_A M)$ is 2-universal with respect to $A$-balanced $A$-$B$-bimodule morphisms. Then, by comparing the 2-universal properties, we find that the 1-morphism $k^F_{M,A}:M\Box F(A)\simeq F(M\Box A)$ induces an equivalence $M\Box_A F(A)\simeq F(M\Box_A A)$ in $\mathbf{Bimod}_{\mathfrak{C}}(A,B)$. Thanks to lemma \ref{lem:relativetensorA}, we have $F(M\Box_A A)\simeq F(M)$ as $A$-$B$-bimodules, which established the claim. Finally, it follows from its construction that the equivalence $M\Box_A F(A)\simeq F(M)$ is 2-natural both in $M$ and in $F$, so that we get $\mathbf{F}\circ\mathbf{B}\simeq Id$. This finishes proving that $\mathbf{F}$ and $\mathbf{B}$ are pseudo-inverses.
\end{proof}

The above theorem yields two equivalent characterizations of Morita equivalence for separable algebras in $\mathfrak{C}$. In specific examples, the second one can be unfolded further as we explain below.

\begin{Corollary}\label{cor:internalMoritaequivalence}
Let $A$ and $B$ be two separable algebras in $\mathfrak{C}$, then the following are equivalent:
\begin{enumerate}
    \item The left $\mathfrak{C}$-module 2-categories $\mathbf{Mod}_{\mathfrak{C}}(A)$ and $\mathbf{Mod}_{\mathfrak{C}}(B)$ are equivalent.
    \item The separable algebras $A$ and $B$ are equivalent as objects of $\mathbf{Mor}^{sep}(\mathfrak{C})$.
\end{enumerate}
If either of these conditions is satisfied, we say that $A$ and $B$ are Morita equivalent.
\end{Corollary}

\begin{Example}
Let $\mathfrak{C}=\mathbf{2Vect}_G$ for some finite group $G$, and, for simplicity, let us also assume that $\mathds{k}$ is algebraically closed of characteristic zero. Then, two $G$-graded multifusion 1-categories $\mathcal{C}$ and $\mathcal{D}$ are Morita equivalent in the sense of 2 of corollary \ref{cor:internalMoritaequivalence} if and only if there exists an invertible $G$-graded finite semisimple $\mathcal{C}$-$\mathcal{D}$-bimodule 1-category $\mathcal{M}$. It follows from proposition 4.2 of \cite{ENO2} that $\mathcal{M}$ is invertible if and only if $\mathcal{D}^{mop}\simeq End_{\mathcal{C}}(\mathcal{M})$ as $G$-graded multifusion 1-categories.

In particular, if $\mathcal{C}$ and $\mathcal{D}$ are faithfully $G$-graded, then they are Morita equivalent in the sense of 2 of corollary \ref{cor:internalMoritaequivalence} if and only if they are graded Morita equivalent in the sense of definition 4.10 of \cite{GJS}. Then, the finite semisimple case of theorem 4.16 of \cite{GJS} asserts that $\mathcal{C}$ and $\mathcal{D}$ are graded Morita equivalent if and only if $\mathbf{Mod}_{\mathbf{2Vect}_G}(\mathcal{C})$ and $\mathbf{Mod}_{\mathbf{2Vect}_G}(\mathcal{D})$ are equivalent as 2-categories with a $G$-action. This is exactly the content of corollary \ref{cor:internalMoritaequivalence} with $\mathfrak{C}=\mathbf{2Vect}_G$.
\end{Example}

\begin{Example}
If $\mathfrak{C}=\mathbf{2Rep}(G)$ for some finite group $G$, and  $\mathds{k}$ is algebraically closed of characteristic zero, then point 2 of corollary \ref{cor:internalMoritaequivalence} can be unpacked further. Namely, let $\mathcal{C}$ be a multifusion 1-category $\mathcal{C}$ with $G$-action, and let $\mathcal{M}$ be finite semisimple 1-category equipped with a $G$-action and a compatible left $\mathcal{C}$-module structure. Then, the multifusion 1-category $End_{\mathcal{C}}(\mathcal{M})$ of all left $\mathcal{C}$-module endofunctors on $\mathcal{M}$ carries a canonical $G$-action. It follows from proposition 4.2 of \cite{ENO2} that two multifusion 1-categories $\mathcal{C}$ and $\mathcal{D}$ with $G$-actions are Morita equivalent if and only if there exists a left $\mathcal{C}$-module 1-category $\mathcal{M}$ as above such that $\mathcal{D}^{mop}\simeq End_{\mathcal{C}}(\mathcal{M})$ as multifusion 1-categories with a $G$-action.
\end{Example}

We now prove an alternative characterization of separability for a compact semisimple left $\mathfrak{C}$-module 2-category $\mathfrak{M}$ under mild assumptions on the underlying 2-categories of $\mathfrak{C}$ and $\mathfrak{M}$. More precisely, following \cite{D5}, if $\mathds{k}$ is perfect, we say that a compact semisimple 2-category $\mathfrak{A}$ is locally separable if for every simple object $A$ of $\mathfrak{A}$, the finite semisimple tensor 1-category $End_{\mathfrak{A}}(A)$ is separable in the sense of \cite{DSPS13}.

\begin{Proposition}\label{prop:characterizationseparablemodule}
Let $\mathds{k}$ be a perfect field, and $\mathfrak{C}$ a locally separable compact semisimple tensor 2-category. The locally separable compact semisimple left $\mathfrak{C}$-module 2-category $\mathfrak{M}$ is separable if and only if $\mathbf{End}_{\mathfrak{C}}(\mathfrak{M})$ is a compact semisimple 2-category.
\end{Proposition}
\begin{proof}
The forward direction follows by combining theorem \ref{thm:equivalencealgebrasmodules2categories} above with proposition 3.1.3 of \cite{D7}. Conversely, let us assume that $\mathbf{End}_{\mathfrak{C}}(\mathfrak{M})$ is compact semisimple. Thanks to theorem 5.3.4 and remark 5.3.10 of \cite{D4}, there exists an algebra $A$ in $\mathfrak{C}$ such that $\mathbf{Mod}_{\mathfrak{C}}(A)\simeq\mathfrak{M}$ as left $\mathfrak{C}$-module 2-categories. We will use the 2-functor $$\mathbf{B}:\mathbf{End}_{\mathfrak{C}}(\mathbf{Mod}_{\mathfrak{C}}(A))\rightarrow \mathbf{Bimod}_{\mathfrak{C}}(A),$$ sending a left $\mathfrak{C}$-module 2-functor to its value on the canonical $A$-$A$-bimodule $A$. Firstly, observe that the image under $\mathbf{B}$ of the identity 2-functor $Id$ on $\mathbf{Mod}_{\mathfrak{C}}(A)$ is given by $A$. Further, if we write $F:\mathbf{Mod}_{\mathfrak{C}}(A)\rightarrow \mathbf{Mod}_{\mathfrak{C}}(A)$ for the canonical left $\mathfrak{C}$-module 2-functor given by $M\mapsto M\Box A$, then we have $\mathbf{B}(F)=F(A)=A\Box A$. Secondly, observe that for any right $A$-module $M$, $n^M:M\Box A\rightarrow M$ defines a left $\mathfrak{C}$-module 2-natural transformation $n:F\Rightarrow Id$ such that $\mathbf{B}(n) = m:A\Box A\rightarrow A$ with its canonical $A$-$A$-bimodule structure. But $\mathbf{End}_{\mathfrak{C}}(\mathfrak{M})$ has right adjoints for 1-morphisms by hypothesis, so that $n$ has a right adjoint $n^{*}$ with counit $\epsilon^n$. As right adjoint are preserved by 2-functors, $\mathbf{B}(n^*)$ is a right adjoint for $m$ as an $A$-$A$-bimodule 1-morphism with counit $\mathbf{B}(\epsilon^n)$. This implies that $A$ is rigid. Thirdly, note that the 2-morphism $\epsilon^n:n\cdot n^*\Rightarrow Id$ is surjective. Namely, for every simple object $M$ of $\mathfrak{M}\simeq \mathbf{Mod}_{\mathfrak{C}}(A)$, the 2-morphism $\epsilon^n_M$ is surjective as $n^M:M\Box A\rightarrow M$ is a non-zero 1-morphism. But $\mathbf{End}_{\mathfrak{C}}(\mathfrak{M})$ is compact semisimple, so that $\epsilon^n$ has a section $\gamma^n$ as a left $\mathfrak{C}$-module modification. Then, $\mathbf{B}(\gamma^n)$ is a section of $\mathbf{B}(\epsilon^n)$ as an $A$-$A$-bimodule 2-morphism  so that $A$ is separable, and the proof is complete.
\end{proof}

The proof of the above proposition also establishes the following result (over any field $\mathds{k}$ and compact semisimple tensor 2-category $\mathfrak{C}$).

\begin{Corollary}\label{cor:Moritainvarianceseparability}
Assume that $\mathfrak{M}$ is a separable left $\mathfrak{C}$-module 2-category, and let $B$ be any algebra such that $\mathfrak{M}\simeq \mathbf{Mod}_{\mathfrak{C}}(B)$ as left $\mathfrak{C}$-module 2-categories, then $B$ is separable.
\end{Corollary}

\begin{Remark}
Let us call two arbitrary algebras $A$, and $B$ in $\mathfrak{C}$ Morita equivalent if $\mathbf{Mod}_{\mathfrak{C}}(A)\simeq\mathbf{Mod}_{\mathfrak{C}}(B)$ as left $\mathfrak{C}$-module 2-categories. Corollary \ref{cor:Moritainvarianceseparability} above may then be succinctly reformulated as the statement that separability is a Morita invariant property. Further, rigidity is also a Morita invariant property. On one hand, if $A$ is a rigid algebra, then $\mathbf{Mod}_{\mathfrak{C}}(A)$ has right adjoints by theorem 2.2.8 of \cite{D7}, so that $\mathbf{End}_{\mathfrak{C}}(\mathbf{Mod}_{\mathfrak{C}}(A))$ has right adjoints. On the other hand, the proof of proposition \ref{prop:characterizationseparablemodule} above, proves that if $\mathbf{End}_{\mathfrak{C}}(\mathbf{Mod}_{\mathfrak{C}}(A))$ has right adjoints, then $A$ is rigid.
\end{Remark}

We end this section by examining an example in detail.

\begin{Example}\label{ex:finitegroupgradeddual}
Let $G$ be a finite group whose order is coprime to $char(\mathds{k})$. Then, the monoidal forgetful 2-functor $\mathbf{2Vect}_G\rightarrow\mathbf{2Vect}$ provides $\mathbf{2Vect}$ with a canonical left $\mathbf{2Vect}_G$-module structure. We claim that $\mathbf{End}_{\mathbf{2Vect}_G}(\mathbf{2Vect})\simeq \mathbf{2Rep}(G)$ as 2-categories. Namely, let $F:\mathbf{2Vect}\rightarrow \mathbf{2Vect}$ be a left $\mathbf{2Vect}_G$-module 2-functor. As $\mathbf{2Vect}$ is generated by $\mathbf{Vect}$ under direct sums and splittings of 2-condensation monads, the underlying linear 2-functor $F$ is determined by $V:=F(\mathbf{Vect})$, a perfect 1-category. Further, unfolding the definition, we find that the left $\mathbf{2Vect}_G$-module structure on $F$ yields a $G$-action on $V$. But $\mathbf{2Vect}_G$ is the Cauchy completion of the monoidal 2-category $G\times \mathbf{2Vect}$, so that this $G$-action on $V$ characterizes $F$ completely up to equivalence. A similar argument deals with $\mathbf{2Vect}_G$-module 2-natural transformations and $\mathbf{2Vect}_G$-module modifications, establishing the desired equivalence $\mathbf{End}_{\mathbf{2Vect}_G}(\mathbf{2Vect})\simeq \mathbf{2Rep}(G)$ of 2-categories. An immediate consequence of the above equivalence is that $\mathbf{2Vect}$ is a separable $\mathbf{2Vect}_G$-module 2-category. Over an algebraically closed field of characteristic zero, this equivalence was first observed in section 3.2 of \cite{Delc}.

For later use, we now wish to upgrade this to an equivalence of monoidal 2-categories. Observe that the identity $\mathbf{2Vect}_G$-module 2-endofunctor of $\mathbf{2Vect}$ corresponds to the object $I=\mathbf{Vect}$ of $\mathbf{2Rep}(G)$ under the above equivalence. It follows from the proof of lemma \ref{lem:2repG} that $\mathbf{2Rep}(G)$ is a connected compact semisimple 2-category. By proposition 3.3.4 of \cite{D5}, the monoidal structure on $\mathbf{2Rep}(G)$ is therefore completely determined by the braiding $\beta$ on the finite semisimple tensor 1-category $End_{\mathbf{2Rep}(G)}(I)\simeq \mathbf{Rep}(G)$. But, by proposition \ref{prop:dualmoduleaction}, $\mathbf{2Vect}$ is a left $\mathbf{Mod}(\mathbf{Rep}^{\beta}(G))$-module 2-category. Thus, by definition, there exists a braided monoidal functor $\mathbf{Rep}^{\beta}(G)\rightarrow \mathcal{Z}(\mathbf{Vect})=\mathbf{Vect}$. As this functor is necessarily faithful, this forces the braiding $\beta$ to be the trivial one, so that $\mathbf{End}_{\mathbf{2Vect}_G}(\mathbf{2Vect})\simeq \mathbf{2Rep}(G)$ as monoidal 2-categories.
\end{Example}

\subsection{Indecomposable Module 2-Categories}
Let $\mathfrak{C}$ be a compact semisimple tensor 2-category, and $\mathfrak{M}$ a compact semisimple left $\mathfrak{C}$-module 2-category. It is useful to know when the compact semisimple monoidal 2-category $\mathbf{End}_{\mathfrak{C}}(\mathfrak{M})$ has simple monoidal unit. We now explain when this is the case.

\begin{Definition}\label{def:indecomposablemodule}
A compact semisimple left $\mathfrak{C}$-module 2-category $\mathfrak{M}$ is indecomposable if there exists a simple object $M$ of $\mathfrak{M}$ such that for any simple object $N$ of $\mathfrak{M}$, there exists an object $C$ in $\mathfrak{C}$ and a non-zero 1-morphisms $C\Box M\rightarrow N$.
\end{Definition}

\begin{Example}
A left $\mathbf{2Vect}$-module 2-category is indecomposable if and only if the underlying compact semisimple 2-category is connected. More generally, if $\mathfrak{C}$ is a connected compact semisimple tensor 2-category, then a left $\mathfrak{C}$-module 2-category is indecomposable if and only if the underlying compact semisimple 2-category is connected.
\end{Example}

\begin{Lemma}\label{lem:decompositionindecomposable}
Let $\mathfrak{M}$ be a compact semisimple left $\mathfrak{C}$-module 2-category. There exists a decomposition $$\mathfrak{M}\simeq \boxplus_{i=1}^n\mathfrak{M}_i$$ into a direct sum of indecomposable compact semisimple left $\mathfrak{C}$-module 2-categories.
\end{Lemma}
\begin{proof}
Let $\mathscr{O}(\mathfrak{M})$ denote the finite set of equivalence classes of simple objects of $\mathfrak{M}$. Let $M$, $N$ be two (equivalence classes of) simple objects in $\mathfrak{M}$, we write $M\sim N$ if there exists an object $C$ of $\mathfrak{C}$ and a non-zero 1-morphism $f:C\Box M\rightarrow N$. This relation is evidently reflexive, symmetry follows from lemma 2.2.10 of \cite{D4}, and transitivity from lemma 2.2.11 of \cite{D4}. Let us write $\mathscr{O}(\mathfrak{M})/\sim\ =\{X_1,...,X_n\}$, and let $\mathfrak{M}_i$ be the full compact semisimple sub-2-category of $\mathfrak{M}$ generated under direct sums and splittings of 2-condensation monads by the simple objects in $X_i$. As the relation $\sim$ is coarser than that of being connected, the sub-2-categories $\mathfrak{M}_i$ and $\mathfrak{M}_j$ do not contain any common simple object. Furthermore, it is immediate from the definition of $\sim$ that $\mathfrak{M}_i$ inherits a left $\mathfrak{C}$-module structure, under which it is indecomposable. Thence, we find $\mathfrak{M}\simeq \boxplus_{i=1}^n\mathfrak{M}_i$ as left $\mathfrak{C}$-module 2-categories.
\end{proof}

\begin{Lemma}
Let $\mathfrak{M}$ be a compact semisimple left $\mathfrak{C}$-module 2-category. Then, the identity left $\mathfrak{C}$-module 2-functor on $\mathfrak{M}$ splits as the direct sum of the projectors onto the $\mathfrak{M}_i$. Further, if $\mathfrak{M}$ is separable, every such projector is a simple object of $\mathbf{End}_{\mathfrak{C}}(\mathfrak{M})$.
\end{Lemma}
\begin{proof}
The first assertion is immediate. Let us assume that $\mathfrak{M}$ is separable, and $\mathfrak{M}\simeq \boxplus_{i=1}^n\mathfrak{M}_i$ be a decomposition of $\mathfrak{M}$ as a direct sum of indecomposable compact semisimple left $\mathfrak{C}$-module 2-categories. We wish to prove that the projection $P_i:\mathfrak{M}\twoheadrightarrow \mathfrak{M}_i\hookrightarrow \mathfrak{M}$ is a simple object of $\mathbf{End}_{\mathfrak{C}}(\mathfrak{M})$. To this end, let $Q,R:\mathfrak{M}\rightarrow \mathfrak{M}$ be two $\mathfrak{C}$-module 2-functor such that $P_i=Q\boxplus R$. Let us additionally assume that $Q(M)$ is non-zero for some $M$ (necessarily in $\mathfrak{M}_i$). Then, it follows from the proof of lemma \ref{lem:decompositionindecomposable} that given any simple object $N$ of $\mathfrak{M}_i$, there exists an object $C$ of $\mathfrak{C}$ and a non-zero 1-morphism $C\Box N\rightarrow M$. In particular, $M$ is the splitting of a 2-condensation monad supported on $C\Box N$. But splittings of 2-condensation monads are preserved by all 2-functors, so that $M$ is the splitting of a 2-condensation monad on $Q(C\Box N)$. As $M$ is non-zero, so must be $Q(C\Box N)$. Now, $Q$ is a left $\mathfrak{C}$-module 2-functor, so that $Q(C\Box N)\simeq C\Box Q(N)$, which implies that $Q(N)$ is non-zero. Finally, we have $N=P_i(N)=Q(N)\boxplus R(N)$, and $N$ is simple, so that $R(N)=0$ by proposition 1.1.7 of \cite{D5}. As $N$ was arbitrary, we find that $R=0$, which finishes the proof of the lemma.
\end{proof}

\begin{Corollary}
Let $\mathfrak{M}$ be a separable left $\mathfrak{C}$-module 2-category. Then $\mathfrak{M}$ is indecomposable if and only if $\mathbf{End}_{\mathfrak{C}}(\mathfrak{M})$ has simple monoidal unit.
\end{Corollary}

Given the equivalence of 3-categories established in theorem \ref{thm:equivalencealgebrasmodules2categories}, it is only natural to examine what property of a separable algebra corresponds to indecomposability of the associated module 2-category.

\begin{Definition}
Let $A$ be a separable algebra. We say that $A$ is indecomposable if $A$ is simple as an $A$-$A$-bimodule.
\end{Definition}

\begin{Corollary}
A separable algebra $A$ is indecomposable if and only if $\mathbf{Mod}_{\mathfrak{C}}(A)$ is indecomposable.
\end{Corollary}

\begin{Remark}
In particular, this shows that being indecomposable is a Morita invariant property of separable algebras in $\mathfrak{C}$. Further, it follows from lemma 5.2.6 that any separable algebra $A$ may be split into a direct sum of indecomposable separable algebras. A direct proof of this fact is given in the proof of theorem 3.1.6 of \cite{D7}.
\end{Remark}

Finally, definition \ref{def:indecomposablemodule} admits an obvious analogue for bimodule 2-categories. This yields a notion of indecomposability for compact semisimple tensor 2-categories.

\begin{Definition}
We say that the compact semisimple tensor 2-category $\mathfrak{C}$ is indecomposable if it is indecomposable as a $\mathfrak{C}$-$\mathfrak{C}$-bimodule 2-category.
\end{Definition}

The proof of lemma \ref{lem:decompositionindecomposable} can be adapted in the obvious so as to give the following result.

\begin{Lemma}\label{lem:indecomposabletensor2categories}
Every compact semisimple tensor 2-category $\mathfrak{C}$ splits as a direct sum of finitely many indecomposable compact semisimple tensor 2-categories.
\end{Lemma}

\subsection{Dual Tensor 2-Categories}

In this section, we assume throughout that $\mathds{k}$ is perfect. The following technical result is needed to prove our main theorem over fields of positive characteristic. Before stating it, we need to recall some terminology from \cite{D5}. We say that a compact semisimple 2-category $\mathfrak{C}$ is locally separable if for every simple object $C$ of $\mathfrak{C}$, the finite semisimple tensor 1-category $End_{\mathfrak{C}}(C)$ is separable in the sense of \cite{DSPS13}. We remark that, over fields of characteristic zero, they show that every finite semisimple tensor 1-category is separable, so that every compact semisimple 2-category is locally separable over such fields.

\begin{Proposition}\label{prop:localseparabilitymodule2categories}
Let $\mathfrak{C}$ be a locally separable compact semisimple monoidal 2-category, and $A$ a separable algebra in $\mathfrak{C}$. Then, $\mathbf{Mod}_{\mathfrak{C}}(A)$ is locally separable. 
\end{Proposition}
\begin{proof}
By theorem 1.4.7 of \cite{D5}, there exists a separable finite semisimple tensor 1-category $\mathcal{C}$ such that $\mathbf{Mod}(\mathcal{C})\simeq \mathfrak{C}$ as 2-categories. It follows from theorem \ref{thm:equivalencealgebrasmodules2categories} that $$\mathbf{End}(\mathfrak{C})\simeq \mathbf{Bimod}(\mathcal{C})^{mop}.$$ Further, as $\mathcal{C}$ is locally separable, theorem 3.1.6 and corollary 3.1.7 of \cite{D7} imply that $\mathbf{Bimod}(\mathcal{C})$ is compact semisimple. Then, thanks to corollary \ref{cor:algebrabimodule}, we find that separable algebras in $\mathbf{Bimod}(\mathcal{C})$ are precisely given by separable finite semisimple tensor 1-categories $\mathcal{D}$ equipped with a monoidal functor $\mathcal{C}\rightarrow \mathcal{D}$. Now, observe that the separable algebra $A$ in $\mathfrak{C}$ yields a separable algebra $\mathcal{A}$ in $\mathbf{End}(\mathfrak{C})\simeq \mathbf{Bimod}(\mathcal{C})^{mop}$ via $C\mapsto C\Box A$. Further, it follows from the definition that $$\mathbf{Mod}_{\mathfrak{C}}(A)\simeq \mathbf{Mod}_{\mathfrak{C}}(\mathcal{A}),$$ where, on the right hand-side, we use the canonical right $\mathbf{End}(\mathfrak{C})^{mop}$-module structure on $\mathfrak{C}$.

Now, note that $End(\mathcal{C})$, the finite semisimple tensor 1-category of linear endofunctors of $\mathcal{C}$, is a separable algebra in $\mathbf{Bimod}(\mathcal{C})$ via the left action of $\mathcal{C}$ on itself. Further, there are equivalences of right $\mathbf{Bimod}(\mathcal{C})$-module 2-categories $$\mathbf{LMod}_{\mathbf{Bimod}(\mathcal{C})}(End(\mathcal{C}))\simeq \mathbf{Bimod}(End(\mathcal{C}), \mathcal{C})\simeq \mathbf{Mod}(\mathcal{C})\simeq \mathfrak{C},$$ as $End(\mathcal{C})$ and $\mathbf{Vect}$ are Morita equivalent finite semisimple tensor 1-categories. Putting everything together, we find that there are equivalences of 2-categories \begin{align*}\mathbf{Mod}_{\mathfrak{C}}(A)&\simeq \mathbf{Bimod}_{\mathbf{Bimod}(\mathcal{C})}(End(\mathcal{C}),\mathcal{A})\\ &\simeq \mathbf{Bimod}(End(\mathcal{C}),\mathcal{A})\simeq \mathbf{Mod}(End(\mathcal{C})^{mop}\boxtimes\mathcal{A}).\end{align*} But it follows from corollary 2.5.11 of \cite{DSPS13} that $End(\mathcal{C})^{mop}\boxtimes\mathcal{A}$ is a separable finite semisimple tensor 1-category, so that $\mathbf{Mod}(End(\mathcal{C})^{mop}\boxtimes\mathcal{A})$ is locally separable by theorem 1.4.6 of \cite{D5}.
\end{proof}

We are now in the position to prove the following theorem.

\begin{Theorem}\label{thm:bimodrigidmonoidal}
Let $\mathds{k}$ be a perfect field, and $A$ a separable algebra in a locally separable compact semisimple tensor 2-category $\mathfrak{C}$. Then, $$\mathbf{End}_{\mathfrak{C}}(\mathbf{Mod}_{\mathfrak{C}}(A))\simeq \mathbf{Bimod}_{\mathfrak{C}}(A)^{mop}$$ is a compact semisimple tensor 2-category.
\end{Theorem}
\begin{proof}
The equivalence of monoidal 2-categories is an immediate consequence of theorem \ref{thm:equivalencealgebrasmodules2categories}. Furthermore, it follows from theorem 3.1.6 of \cite{D7} that the underlying 2-category of $\mathbf{Bimod}_{\mathfrak{C}}(A)$ is compact semisimple. Thus, it only remains to establish the existence of duals. We will show that $\mathbf{End}_{\mathfrak{C}}(\mathbf{Mod}_{\mathfrak{C}}(A))$ satisfies this property. Namely, it follows from proposition \ref{prop:localseparabilitymodule2categories} that $\mathbf{Mod}_{\mathfrak{C}}(A)$ is locally separable. Then, corollary 3.2.3 of \cite{D5} shows that every linear 2-functor $\mathbf{Mod}_{\mathfrak{C}}(A)\rightarrow \mathbf{Mod}_{\mathfrak{C}}(A)$ has a right 2-adjoint 2-functor. In particular, proposition \ref{prop:C2Funrigid} applies to every object of $\mathbf{End}_{\mathfrak{C}}(\mathbf{Mod}_{\mathfrak{C}}(A))$, which proves that $\mathbf{End}_{\mathfrak{C}}(\mathbf{Mod}_{\mathfrak{C}}(A))$ has right duals. By corollary 1.3.4 of \cite{D2}, $\mathbf{End}_{\mathfrak{C}}(\mathbf{Mod}_{\mathfrak{C}}(A))$ also has left duals, which concludes the proof of the result.
\end{proof}

\begin{Remark}
The assumption that $\mathfrak{C}$ be locally separable in theorem \ref{thm:bimodrigidmonoidal} might not be necessary. Namely, we believe that it is possible to show directly that for any separable algebra $A$ in a compact semisimple tensor 2-category, the monoidal 2-category $\mathbf{Bimod}_{\mathfrak{C}}(A)$ has duals. We leave it to the interested reader to pursue this line of investigation.
\end{Remark}

Thanks to the above theorem, the following definition is sensible.

\begin{Definition}
Let $\mathfrak{C}$ be a locally separable compact semisimple tensor 2-category, and let $\mathfrak{M}$ be a separable left $\mathfrak{C}$-module 2-category. We write $\mathfrak{C}_{\mathfrak{M}}^*$ for the compact semisimple tensor 2-category $\mathbf{End}_{\mathfrak{C}}(\mathfrak{M})$, and call it the dual tensor 2-category to $\mathfrak{C}$ with respect to $\mathfrak{M}$.
\end{Definition}

Combining the theorem \ref{thm:bimodrigidmonoidal} above with corollary 2.2.4 of \cite{D5}, we obtain the following result.

\begin{Corollary}\label{cor:dualfusion2category}
Let $\mathds{k}$ be an algebraically closed field of characteristic zero, $\mathfrak{C}$ a multifusion 2-category, and $\mathfrak{M}$ a separable left $\mathfrak{C}$-module 2-category. Then, $\mathfrak{C}_{\mathfrak{M}}^*$, the dual tensor 2-category to $\mathfrak{C}$ with respect to $\mathfrak{M}$, is a multifusion 2-category.
\end{Corollary}

The following corollary follows from the discussion given in example \ref{ex:algebrastwistedgroupgraded}.

\begin{Corollary}
Let $G$ be a finite group, and $\pi$ a 4-cocycle for $G$ with coefficients in $\mathds{k}^{\times}$. Further, let $H\subseteq G$ be a subgroup of order coprime to $char(\mathds{k})$ and $\gamma$ a 3-cochain for $H$ with coefficients in $\mathds{k}^{\times}$. Then, $\mathbf{Bimod}_{\mathbf{2Vect}^{\pi}_G}(\mathbf{Vect}^{\gamma}_H)$ is a compact semisimple tensor 2-category.
\end{Corollary}

We end this section by examining some examples of dual tensor 2-categories.

\begin{Example}\label{ex:dual2Vect}
Let $\mathds{k}$ be an algebraically closed field, and let $\mathcal{C}$ be a separable fusion 1-category. The locally separable compact semisimple 2-category $\mathbf{Mod}(\mathcal{C})$ admits a canonical left $\mathbf{2Vect}$-module structure. We claim that $\mathbf{2Vect}_{\mathbf{Mod}(\mathcal{C})}^*\simeq \mathbf{Mod}(\mathcal{Z}(\mathcal{C}))$ as monoidal 2-categories. Note that $\mathcal{Z}(\mathcal{C})$ is a finite semisimple tensor 1-category thanks to corollary 2.5.9 of \cite{DSPS13}. Theorem \ref{thm:equivalencealgebrasmodules2categories} provides us with an equivalence of monoidal 2-categories $\mathbf{2Vect}_{\mathbf{Mod}(\mathcal{C})}^*\simeq \mathbf{Bimod}(\mathcal{C})^{mop}$. Provided $\mathds{k}$ has characteristic zero, theorem 1.3 of \cite{Gre} gives an equivalence of monoidal 2-categories $\mathbf{Bimod}(\mathcal{C})\simeq \mathbf{LMod}(\mathcal{Z}(\mathcal{C}))$, to the monoidal 2-category of finite semisimple left $\mathcal{Z}(\mathcal{C})$-module 1-categories. Finally, we have $\mathbf{LMod}(\mathcal{Z}(\mathcal{C}))\simeq \mathbf{Mod}(\mathcal{Z}(\mathcal{C})^{\beta op})$ as monoidal 2-categories, which concludes the proof.

Alternatively, over an arbitrary algebraically closed field $\mathds{k}$, as $\mathcal{C}$ is a separable fusion 1-category, $\mathcal{C}^{mop}\boxtimes\mathcal{C}$ is also a separable fusion 1-category, so that the compact semisimple 2-category $\mathbf{Bimod}(\mathcal{C})\simeq \mathbf{Mod}(\mathcal{C}^{mop}\boxtimes\mathcal{C})$ is connected. Thus, by proposition 3.3.4 of \cite{D5}, in order to determine the monoidal structure on $\mathbf{Bimod}(\mathcal{C})$, it is enough to understand the braiding on the fusion 1-category $End_{\mathcal{C}\textrm{-}\mathcal{C}}(\mathcal{C})$ of endomorphisms the monoidal unit. Inspection shows that it is equivalent to $\mathcal{Z}(\mathcal{C})^{mop}$ as a braided monoidal 1-category. But, there is an equivalence of braided monoidal 1-categories $\mathcal{Z}(\mathcal{C})^{mop}\simeq \mathcal{Z}(\mathcal{C})^{\beta op}$. Thus, we have $End_{\mathcal{C}\textrm{-}\mathcal{C}}(\mathcal{C})\simeq \mathcal{Z}(\mathcal{C})^{\beta op}$ as braided fusion 1-categories.
\end{Example}

\begin{Example}\label{ex:dualnondegenerate}
We now discuss a generalization of example \ref{ex:dual2Vect}. For simplicity, we will assume that $\mathds{k}$ is an algebraically closed field of characteristic zero. Let $\mathcal{B}$ be a non-degenerate braided fusion 1-category, and $\mathcal{C}$ a $\mathcal{B}$-central fusion 1-category. We claim that there is an equivalence of monoidal 2-categories $\mathbf{Mod}(\mathcal{B})_{\mathbf{Mod}(\mathcal{C})}^*\simeq \mathbf{Mod}(\mathcal{A})$, where $\mathcal{A}$ is the centralizer of $\mathcal{B}$ in $\mathcal{Z}(\mathcal{C})$, which is non-degenerate by theorem 3.13 of \cite{DGNO}. By corollary 5.9 of \cite{DMNO}, this implies further that $\mathcal{B}$ and $\mathcal{A}$ are Witt equivalent non-degenerate braided fusion 1-categories.

In order to prove the claim, observe that, as $\mathcal{C}$ is a separable algebra in $\mathbf{Mod}(\mathcal{B})$, corollary \ref{cor:dualfusion2category} establishes that $\mathbf{Mod}(\mathcal{B})_{\mathbf{Mod}(\mathcal{C})}^*\simeq \mathbf{Bimod}_{\mathbf{Mod}(\mathcal{B})}(\mathcal{C})^{mop}$ is a fusion 2-category. Moreover, after having unfolded the definitions, we find that $\mathbf{Bimod}_{\mathbf{Mod}(\mathcal{B})}(\mathcal{C})\simeq \mathbf{Mod}(\mathcal{C}^{mop}\boxtimes_{\mathcal{B}}\mathcal{C})$ as finite semisimple 2-categories. As $\mathcal{B}$ is non-degenerate, it follows from theorems 2.26 and 3.20 of \cite{BJSS} that there is an equivalence $$\mathcal{C}^{mop}\boxtimes_{\mathcal{B}}\mathcal{C}\simeq End_{\mathcal{A}}(\mathcal{C})$$ of multifusion 1-categories. In particular, $\mathbf{Bimod}_{\mathbf{Mod}(\mathcal{B})}(\mathcal{C})\simeq \mathbf{Mod}(\mathcal{C}^{mop}\boxtimes_{\mathcal{B}}\mathcal{C})$ is a connected finite semisimple 2-category. The above equivalence of multifusion 1-categories implies that $\mathcal{A}\simeq End_{\mathcal{C}^{mop}\boxtimes_{\mathcal{B}}\mathcal{C}}(\mathcal{C})$ as fusion 1-categories, so that the endomorphism fusion 1-category of $\mathcal{C}$ in $\mathbf{Bimod}_{\mathbf{Mod}(\mathcal{B})}(\mathcal{C})$ is given by $\mathcal{A}$. But $\mathcal{C}$ is the monoidal unit of $\mathbf{Bimod}_{\mathbf{Mod}(\mathcal{B})}(\mathcal{C})$. Thence, appealing to proposition 2.4.7 of \cite{D2}, it is enough to understand the braiding on the fusion 1-category $\mathcal{A}$ of endomorphisms of $\mathcal{C}$ in $\mathbf{Bimod}_{\mathbf{Mod}(\mathcal{B})}(\mathcal{C})$.

To this end, note that the forgetful 2-functor $\mathbf{Bimod}_{\mathbf{Mod}(\mathcal{B})}(\mathcal{C})\rightarrow \mathbf{Bimod}(\mathcal{C})$ induces the canonical inclusion of fusion 1-categories $\mathcal{A}\hookrightarrow \mathcal{Z}(\mathcal{C})^{\beta op}$. But the forgetful monoidal 2-functor $$\mathbf{End}_{\mathbf{Mod}(\mathcal{B})}(\mathbf{Mod}(\mathcal{C}))^{mop}\rightarrow \mathbf{End}(\mathbf{Mod}(\mathcal{C}))^{mop}$$ is identified via theorem \ref{thm:equivalencealgebrasmodules2categories} to $\mathbf{Bimod}_{\mathbf{Mod}(\mathcal{B})}(\mathcal{C})\rightarrow \mathbf{Bimod}(\mathcal{C})$. This shows that the monoidal inclusion $\mathcal{A}\hookrightarrow \mathcal{Z}(\mathcal{C})^{\beta op}$ is in fact braided, thereby establishing the desired equivalence $\mathbf{Bimod}_{\mathbf{Mod}(\mathcal{B})}(\mathcal{C})\simeq \mathbf{Mod}(\mathcal{A}^{\beta op})$ of monoidal 2-categories.
\end{Example}

\begin{Remark}
As one can readily observe from example \ref{ex:finitegroupgradeddual}, the non-degeneracy condition in example \ref{ex:dualnondegenerate} above can not be omitted in general. We will return to this point in \cite{D9}.
\end{Remark}

\subsection{Morita Equivalence}

Let us fix $\mathds{k}$ a perfect field. We give three equivalent characterizations of Morita equivalence between locally separable compact semisimple tensor 2-categories. Before doing so, we need a definition.

Let $\mathfrak{C}$ be a compact semisimple tensor 2-category. It follows from proposition 1.1.7 of \cite{D5} that there is a splitting $I\simeq \boxplus I_i$ of the monoidal unit of $\mathfrak{C}$ into a finite direct sum of simple objects.

\begin{Definition}
Let $\mathfrak{M}$ be a left $\mathfrak{C}$-module 2-category. We say that $\mathfrak{M}$ is faithful if the action of $I_i$ on $\mathfrak{M}$ is non-zero for every $i$. Let $A$ be an algebra in $\mathfrak{C}$, we say that $A$ is faithful if $\mathbf{Mod}_{\mathfrak{C}}(A)$ is a faithful left $\mathfrak{C}$-module 2-category.
\end{Definition}

\begin{Remark}\label{rem:technicalfaithful}
If $\mathfrak{C}$ has simple monoidal unit, then every non-zero module 2-category is faithful. This holds more generally if $\mathfrak{C}$ is indecomposable. In fact, if $\mathfrak{C}\simeq \boxplus\mathfrak{C}_n$ is a splitting of $\mathfrak{C}$ into a finite direct sum of indecomposable compact semisimple tensor 2-categories, then an algebra $A$ in $\mathfrak{C}$ is faithful if and only if its underlying object has a summand in $\mathfrak{C}_n$ for every $n$.
\end{Remark}

\begin{Theorem}\label{thm:Moritaequivalence}
For any two locally separable compact semisimple tensor 2-categories $\mathfrak{C}$ and $\mathfrak{D}$ over a perfect field $\mathds{k}$, the following are equivalent:
\begin{enumerate}
    \item The 3-categories $\mathbf{LMod}^{sep}(\mathfrak{C})$ and $\mathbf{LMod}^{sep}(\mathfrak{D})$ are equivalent.
    \item There exists a faithful separable left $\mathfrak{C}$-module 2-category $\mathfrak{M}$, and an equivalence of monoidal 2-categories $\mathfrak{D}^{mop}\simeq \mathfrak{C}_{\mathfrak{M}}^*$.
    \item There exists a faithful separable algebra $A$ in $\mathfrak{C}$, and an equivalence of monoidal 2-categories $\mathfrak{D}\simeq \mathbf{Bimod}_{\mathfrak{C}}(A)$.
\end{enumerate}
If either of the above conditions is satisfied, we say that $\mathfrak{C}$ and $\mathfrak{D}$ are Morita equivalent.
\end{Theorem}
\begin{proof}
The equivalence between 2 and 3 follows from theorem \ref{thm:equivalencealgebrasmodules2categories}. Let us now assume that there is an equivalence of 3-categories $\mathbf{F}:\mathbf{LMod}^{sep}(\mathfrak{C})\simeq \mathbf{LMod}^{sep}(\mathfrak{D})$, with pseudo-inverse $\mathbf{G}$. As there is a canonical equivalence of monoidal 2-categories $\mathbf{End}_{\mathfrak{D}}(\mathfrak{D})\simeq \mathfrak{D}^{mop}$, and $\mathbf{G}$ induces an equivalence $$\mathbf{End}_{\mathfrak{D}}(\mathfrak{D})\simeq \mathbf{End}_{\mathfrak{C}}(\mathbf{G}(\mathfrak{D}))$$ of monoidal 2-categories, we find that $\mathfrak{D}^{mop}\simeq \mathfrak{C}_{\mathbf{G}(\mathfrak{D})}^*$. Now, suppose that there is a simple summand $I_i$ of the monoidal unit $I$ of $\mathfrak{C}$ that acts as zero on $\mathbf{G}(\mathfrak{D})$. Then, we have $\mathbf{Fun}_{\mathfrak{C}}(\mathbf{Mod}_{\mathfrak{C}}(I_i), \mathbf{G}(\mathfrak{D})) = 0$. On the other hand, as $\mathbf{F}$ is an equivalence of 3-categories, we have that $$\mathbf{Fun}_{\mathfrak{C}}(\mathbf{Mod}_{\mathfrak{C}}(I_i), \mathbf{G}(\mathfrak{D}))\simeq \mathbf{Fun}_{\mathfrak{D}}(\mathbf{F}(\mathbf{Mod}_{\mathfrak{C}}(I_i)),\mathfrak{D}),$$ and the right hand-side is non-zero by proposition \ref{prop:C2Funrigid}, so that $\mathbf{G}(\mathfrak{D})$ is a faithful left $\mathfrak{C}$-module 2-category.

Conversely, let $A$ be a faithful separable algebra $A$ in $\mathfrak{C}$ such that $\mathfrak{D}:=\mathbf{Bimod}_{\mathfrak{C}}(A)$ is a locally separable compact semisimple tensor 2-category. Without loss of generality, we may assume that $\mathfrak{C}$ is strict cubical. Firstly, we claim that there is a faithful separable algebra $B$ in $\mathfrak{D}$ such that $$\mathbf{Bimod}_{\mathfrak{D}}(B)\simeq \mathfrak{C}$$ as monoidal 2-categories.

To see this, recall from corollary \ref{cor:algebrabimodule} that separable algebras in $\mathfrak{D}$ are precisely separable algebras in $\mathfrak{C}$ equipped with an algebra 1-homomorphism from $A$. Now, we fix $C$ an object of $\mathfrak{C}$ that has a simple summand in every connected component of $\mathfrak{C}$, and write $(^{\sharp}C,C,i_C,e_C,E_C,F_C)$ for a coherent left dual for $C$ in the sense of \cite{Pstr} (see also \cite{D2}). We take $B$ to be the algebra in $\mathfrak{C}$ whose underlying object is $A\Box C\Box (^{\sharp}C)\Box A$, and with unit and mutliplication 1-morphisms given by $$i^{B}:I\xrightarrow{i} A\xrightarrow{m^*}A A\xrightarrow{1i_C1} A C (^{\sharp}C) A,$$ \begin{align*}m^{B}:A C (^{\sharp}C) A A C (^{\sharp}C) A&\xrightarrow{111m111}A C (^{\sharp}C) A C (^{\sharp}C) A\xrightarrow{111i^*111}A C (^{\sharp}C)C (^{\sharp}C) A\\ &\xrightarrow{11e_C11}A C (^{\sharp}C) A,\end{align*} where $i^*$ is a right adjoint in $\mathfrak{C}$ for the unit $i$ of $A$. The coherence 2-isomorphisms are defined in the obvious way. Further, it follows from the definition that $(1i_C1)\circ m^*:A\rightarrow B$ is an algebra 1-homomorphism in $\mathfrak{C}$. Moreover, we have that $\mathfrak{C}\simeq \mathbf{Mod}_{\mathfrak{C}}(B)$ as left $\mathfrak{C}$-module 2-categories via $D\mapsto D\Box (^{\sharp}C)\Box A$ for every $D$ in $\mathfrak{C}$. Namely, as $A$ is faithful, $A\Box C$ is a $\mathfrak{C}$-generator of $\mathfrak{C}$ in the sense of definition 5.3.1 of \cite{D4}. In particular, we can apply theorem 5.3.4 of \cite{D4} with $M:=A\Box C$, or more precisely the generalization given in remark 5.3.10 therein. The assertion then follows by combining example 4.1.3 of \cite{D4} with the fact that $A$ is a self-dual object of $\mathfrak{C}$ with coevaluation 1-morphism $m^*\circ i$ and evaluation 1-morphism $i^*\circ m$. For later use, we also record that the 2-functor $\mathfrak{C}\rightarrow\mathbf{LMod}_{\mathfrak{C}}(B)$ given by $D\mapsto A\Box C\Box D$ is an equivalence of right $\mathfrak{C}$-module 2-categories. This follows from the above argument applied to $\mathfrak{C}^{mop}$. Finally, theorem \ref{thm:equivalencealgebrasmodules2categories} implies that there is an equivalence of monoidal 2-categories $\mathfrak{C}\rightarrow\mathbf{Bimod}_{\mathfrak{C}}(B)$, which is given by $D\mapsto A\Box C\Box D \Box (^{\sharp}C)\Box A$. This concludes the proof of the claim.

Secondly, observe that there is a 3-functor $$\mathbf{Fun}_{\mathfrak{C}}(\mathbf{Mod}_{\mathfrak{C}}(A),-):\mathbf{LMod}^{sep}(\mathfrak{C})\rightarrow \mathbf{LMod}(\mathfrak{D}).$$ Up to the equivalence of 3-categories of theorem \ref{thm:equivalencealgebrasmodules2categories}, this 3-functor is equivalent to $\mathbf{Bimod}_{\mathfrak{C}}(A,-):\mathbf{Mor}^{sep}(\mathfrak{C})\rightarrow \mathbf{LMod}(\mathfrak{D}).$ Thanks to the claim above, there is also a 3-functor $$\mathbf{Fun}_{\mathfrak{D}}(\mathbf{Mod}_{\mathfrak{D}}(B),-):\mathbf{LMod}^{sep}(\mathfrak{D})\rightarrow \mathbf{LMod}(\mathfrak{C}).$$ This last 3-functor is equivalent to $\mathbf{Bimod}_{\mathfrak{D}}(B,-):\mathbf{Mor}^{sep}(\mathfrak{D})\rightarrow \mathbf{LMod}(\mathfrak{C})$ via the equivalence of 3-categories of theorem \ref{thm:equivalencealgebrasmodules2categories}. Further, up to the identification $\mathfrak{C}\simeq\mathbf{Bimod}_{\mathfrak{C}}(B)$, it follows from the claim above and its proof that $\mathbf{Bimod}_{\mathfrak{D}}(B,-)\simeq \mathbf{Mod}_{\mathfrak{C}}(-)$ as left $\mathfrak{C}$-module 2-categories. In particular, $\mathbf{Fun}_{\mathfrak{D}}(\mathbf{Mod}_{\mathfrak{D}}(B),-)$ factors through $\mathbf{LMod}^{sep}(\mathfrak{C})$. Putting the above discussion together, we find that there are equivalences of 3-functors \begin{align*}\mathbf{Fun}_{\mathfrak{C}}(\mathbf{Mod}_{\mathfrak{C}}(A),\mathbf{Fun}_{\mathfrak{D}}(&\mathbf{Mod}_{\mathfrak{D}}(B),\mathbf{Mod}_{\mathfrak{D}}(-)))\\ &\simeq  \mathbf{Fun}_{\mathfrak{C}}(\mathbf{Mod}_{\mathfrak{C}}(A),\mathbf{Bimod}_{\mathfrak{D}}(B,-))\\ &\simeq  \mathbf{Fun}_{\mathfrak{C}}(\mathbf{Mod}_{\mathfrak{C}}(A),\mathbf{Mod}_{\mathfrak{C}}(-))\\& \simeq \mathbf{Bimod}_{\mathfrak{C}}(A,-)\simeq \mathbf{Mod}_{\mathfrak{D}}(-).\end{align*} This shows that the composite of the two 3-functors $\mathbf{Fun}_{\mathfrak{D}}(\mathbf{Mod}_{\mathfrak{D}}(B),-)$ and $\mathbf{Fun}_{\mathfrak{C}}(\mathbf{Mod}_{\mathfrak{C}}(A),-)$ is equivalent to the identity on $\mathbf{LMod}^{sep}(\mathfrak{D})$.

Finally, one can run the above argument starting with $\mathfrak{D}$ and $B$. This shows that the 3-functor $\mathbf{Fun}_{\mathfrak{D}}(\mathbf{Mod}_{\mathfrak{D}}(B),-)$ has both a left and a right pseudo-inverse, so that it induces an equivalences of 3-categories $\mathbf{LMod}^{sep}(\mathfrak{D})\simeq \mathbf{LMod}^{sep}(\mathfrak{C})$ as desired.
\end{proof}

Let us record the following corollary of the proof of theorem \ref{thm:Moritaequivalence}.

\begin{Corollary}
Let $\mathfrak{C}$ be a locally separable compact semisimple tensor 2-category, and let $\mathfrak{M}$ be a separable left $\mathfrak{C}$-module 2-category. Then, $\mathfrak{M}$ is a separable left $\mathfrak{C}^*_{\mathfrak{M}}$-module 2-category. Furthermore, if $\mathfrak{M}$ is faithful, the canonical monoidal 2-functor $\mathfrak{C}\rightarrow (\mathfrak{C}^*_{\mathfrak{M}})^*_{\mathfrak{M}}$ is an equivalence.
\end{Corollary}
\begin{proof}
The first part is immediate. For the second part, we use the notations of the proof of theorem \ref{thm:Moritaequivalence}. In particular, there is an equivalence of monoidal 2-categories $$(\mathfrak{C}^*_{\mathfrak{M}})^*_{\mathfrak{M}}\simeq \mathbf{Bimod}_{\mathfrak{C}}(B).$$ Under this equivalence, the canonical monoidal 2-functor $\mathfrak{C}\rightarrow (\mathfrak{C}^*_{\mathfrak{M}})^*_{\mathfrak{M}}$ is identified with the 2-functor $\mathfrak{C}\rightarrow \mathbf{Bimod}_{\mathfrak{C}}(B)$ given by $D\mapsto A\Box C\Box D\Box (^{\sharp}C)\Box A$. But, this 2-functor was shown to be an equivalence in the course of the proof of theorem \ref{thm:Moritaequivalence}.
\end{proof}

We end by examining two examples.

\begin{Example}
Let $G$ be a finite group whose order is coprime to $char(\mathds{k})$. It follows from example \ref{ex:finitegroupgradeddual} above that the locally separable compact semisimple tensor 2-categories $\mathbf{2Vect}_G$ and $\mathbf{2Rep}(G)$ are Morita equivalent. In particular, the 3-categories of separable left $\mathbf{2Vect}_G$-module 2-categories and of separable left $\mathbf{2Rep}(G)$-module 2-categories are equivalent.
\end{Example}

\begin{Example}
Let $\mathds{k}$ be algebraically closed of characteristic zero, and  $\mathcal{B}_1$, $\mathcal{B}_2$ be two non-degenerate braided fusion 1-categories. It follows from example \ref{ex:dualnondegenerate} that the fusion 2-categories $\mathbf{Mod}(\mathcal{B}_1)$ and $\mathbf{Mod}(\mathcal{B}_2)$ are Morita equivalent if and only if the non-degenerate braided fusion 1-categories $\mathcal{B}_1$ and $\mathcal{B}_2$ are Witt equivalent in the sense of \cite{DMNO}.
\end{Example}

%% file: Appendix.tex
\begin{landscape}
    \section{Appendix}

    \subsection{Diagrams for the proof of theorem \ref{thm:tensormodules}}\label{sub:condensation}

    \vspace*{\fill}
    \begin{figure}[!hbt]
    \centering
    \includegraphics[width=180mm]{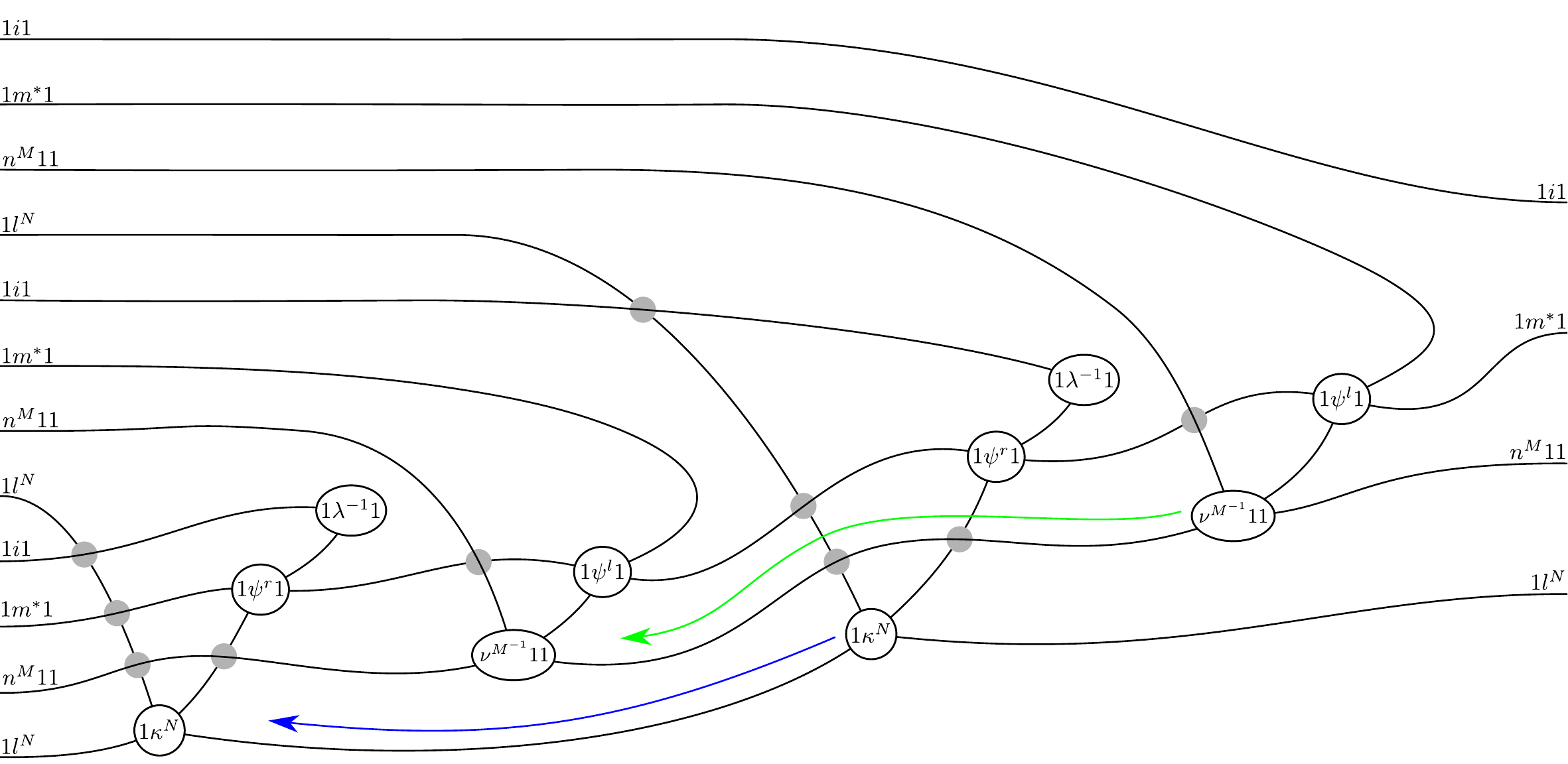}
    \caption{Associativity (Part 1)}
    \label{fig:muassociativity1}
    \end{figure}
    \vfill
\end{landscape}

\begin{landscape}
    \vspace*{\fill}
    \begin{figure}[!hbt]
    \centering
    \includegraphics[width=180mm]{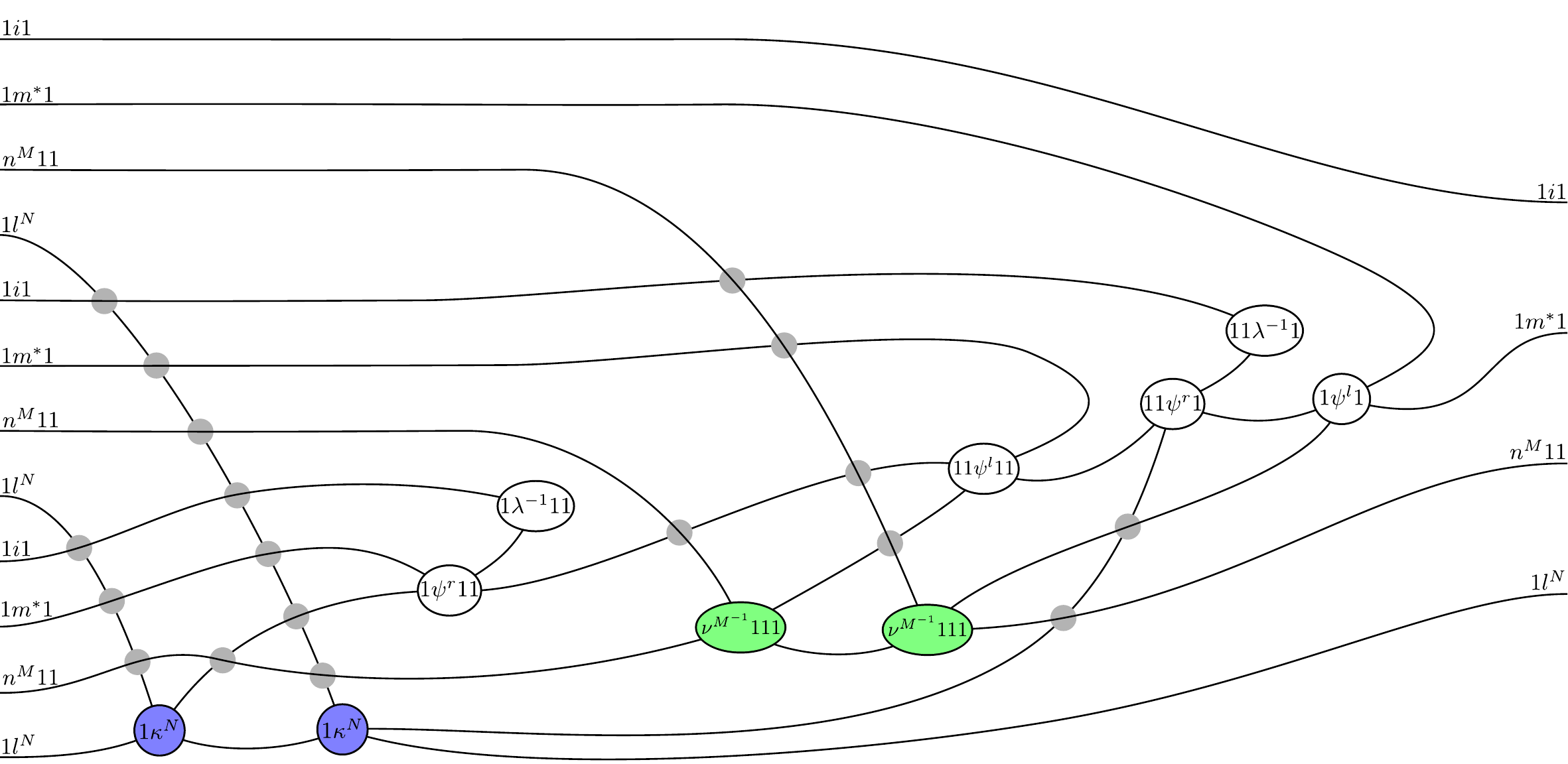}
    \caption{Associativity (Part 2)}
    \label{fig:muassociativity2}
    \end{figure}
    \vfill
\end{landscape}

\begin{landscape}
    \vspace*{\fill}
    \begin{figure}[!hbt]
    \centering
    \includegraphics[width=180mm]{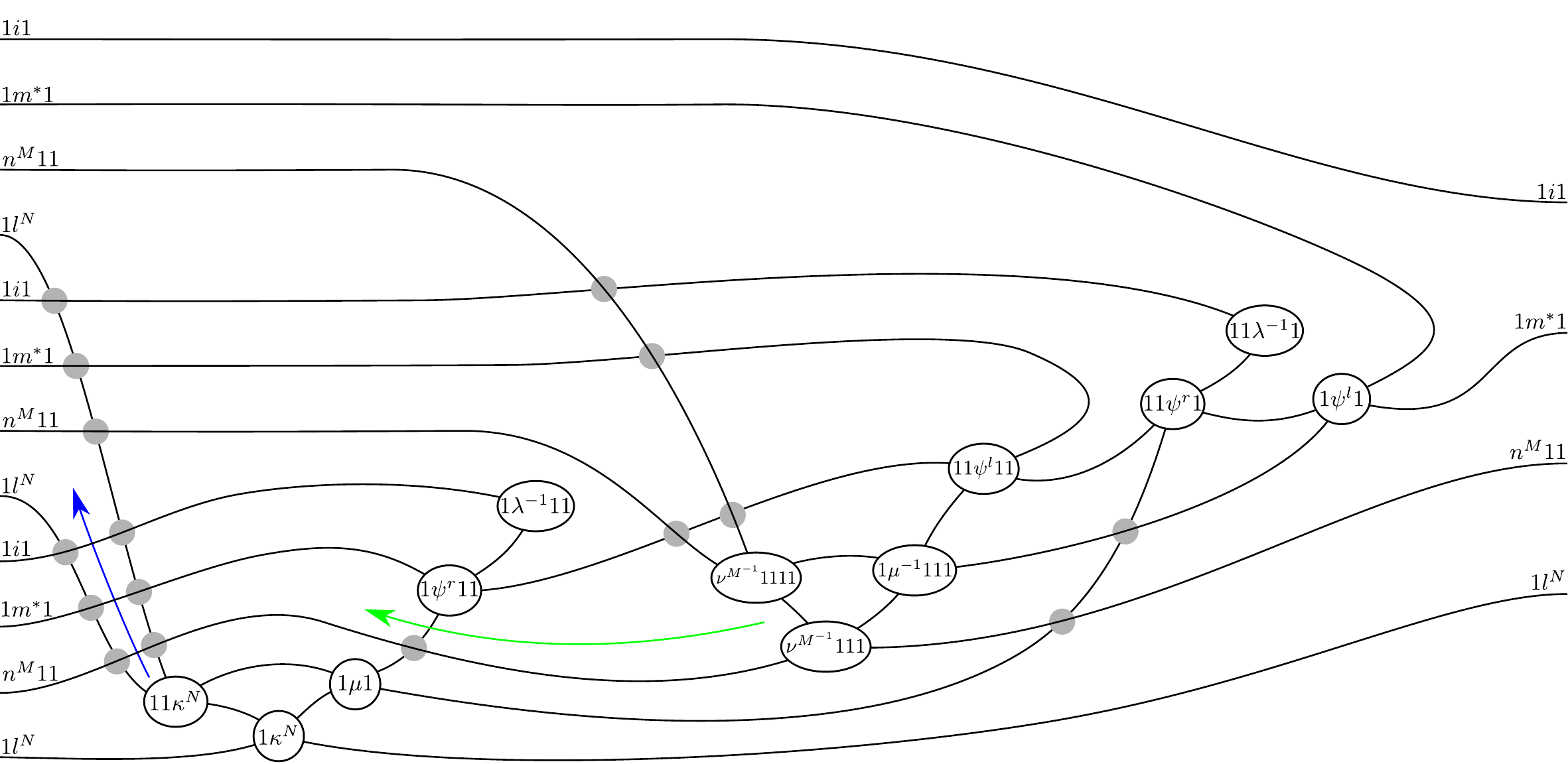}
    \caption{Associativity (Part 3)}
    \label{fig:muassociativity3}
    \end{figure}
    \vfill
\end{landscape}

\begin{landscape}
    \vspace*{\fill}
    \begin{figure}[!hbt]
    \centering
    \includegraphics[width=180mm]{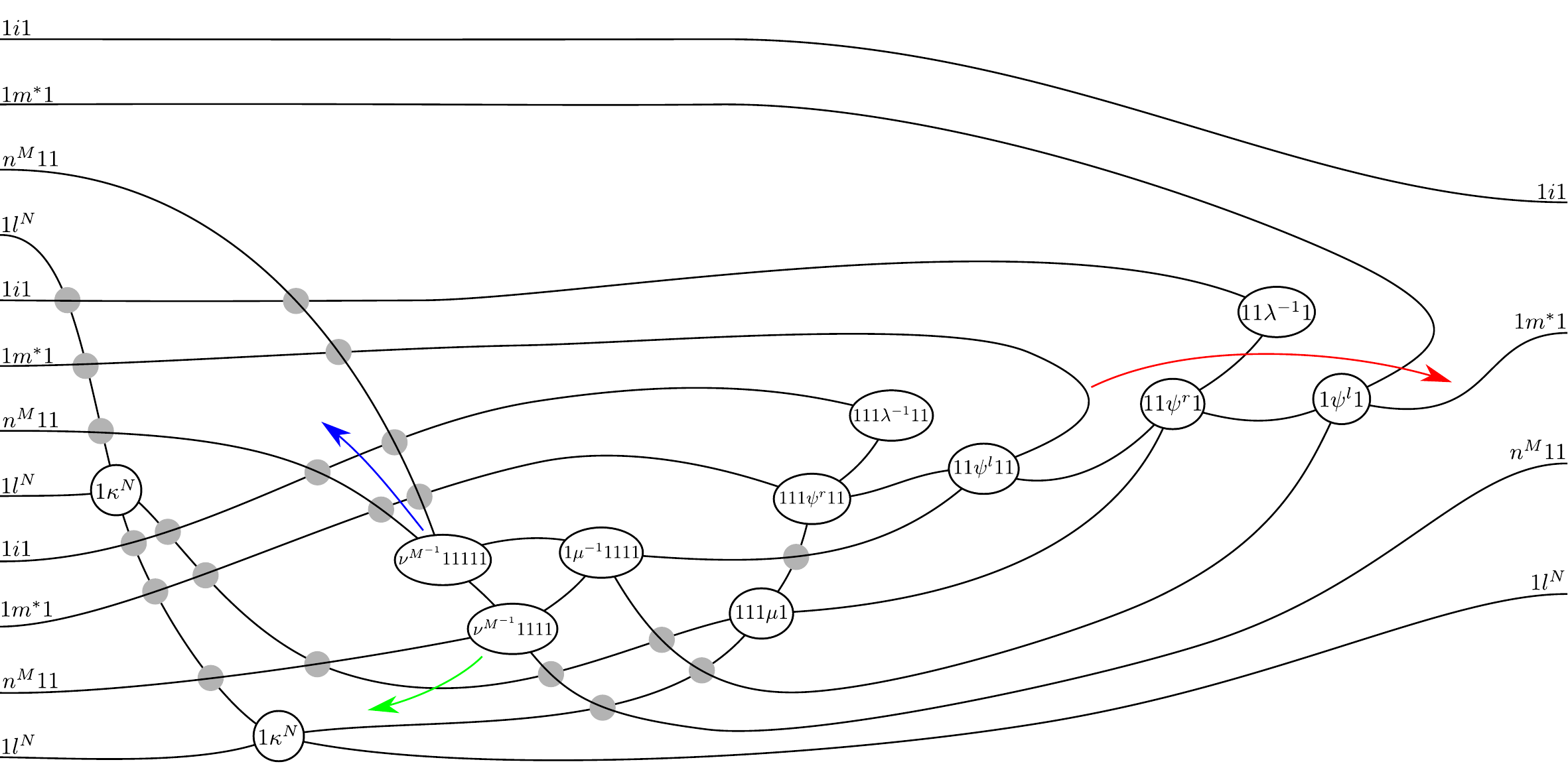}
    \caption{Associativity (Part 4)}
    \label{fig:muassociativity4}
    \end{figure}
    \vfill
\end{landscape}

\begin{landscape}
    \vspace*{\fill}
    \begin{figure}[!hbt]
    \centering
    \includegraphics[width=180mm]{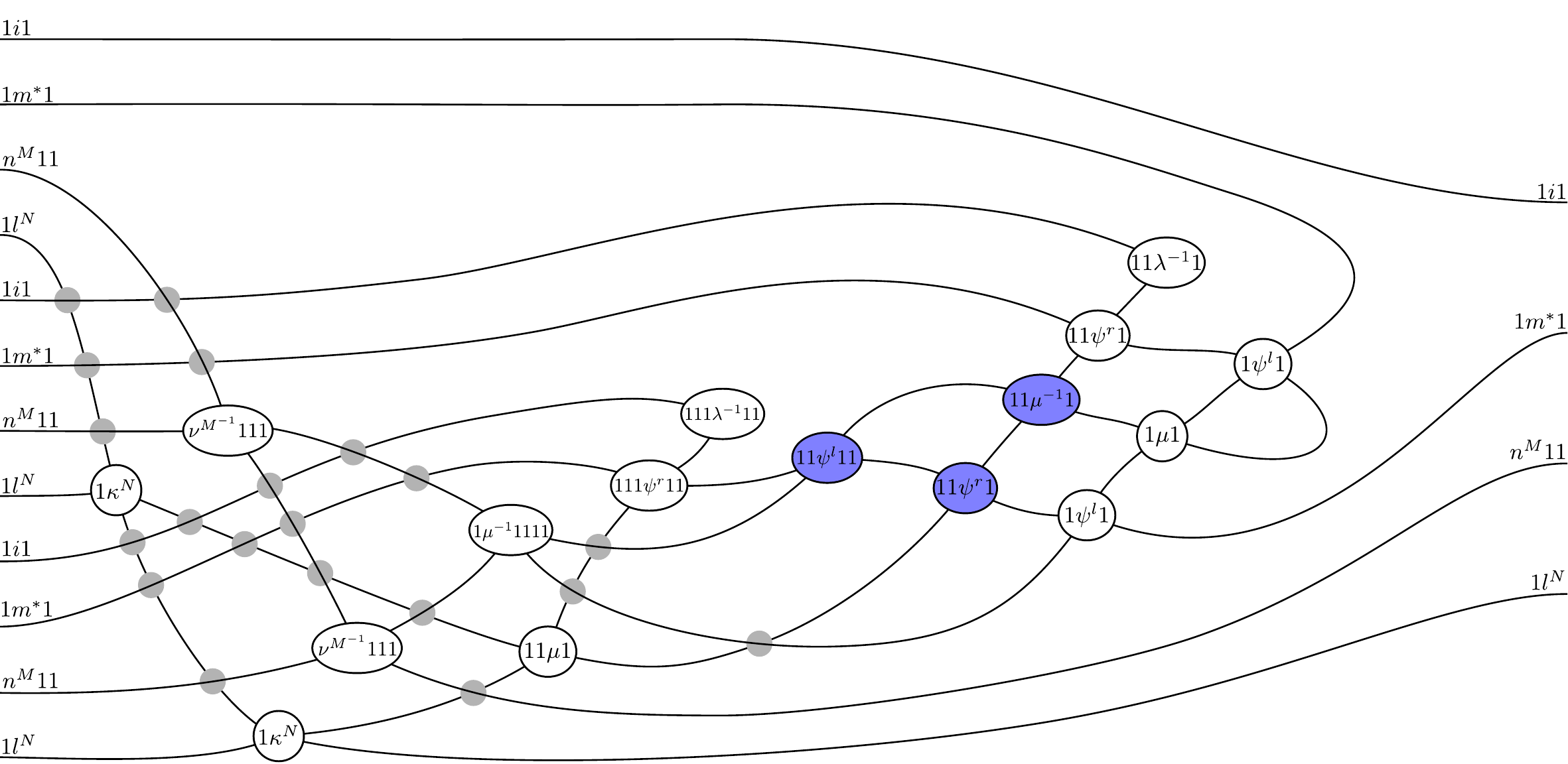}
    \caption{Associativity (Part 5)}
    \label{fig:muassociativity5}
    \end{figure}
    \vfill
\end{landscape}

\begin{landscape}
    \vspace*{\fill}
    \begin{figure}[!hbt]
    \centering
    \includegraphics[width=180mm]{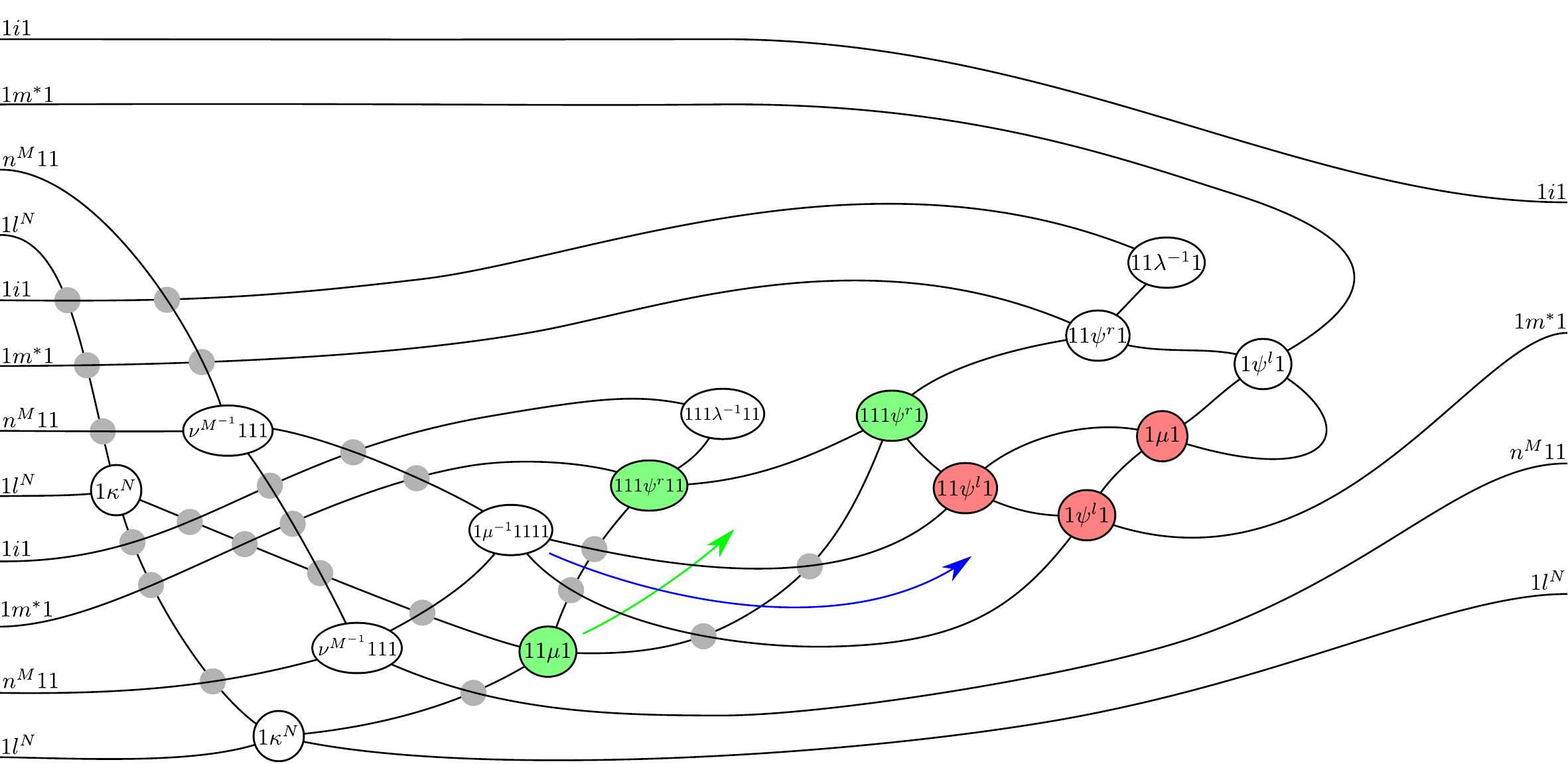}
    \caption{Associativity (Part 6)}
    \label{fig:muassociativity6}
    \end{figure}
    \vfill
\end{landscape}

\begin{landscape}
    \vspace*{\fill}
    \begin{figure}[!hbt]
    \centering
    \includegraphics[width=180mm]{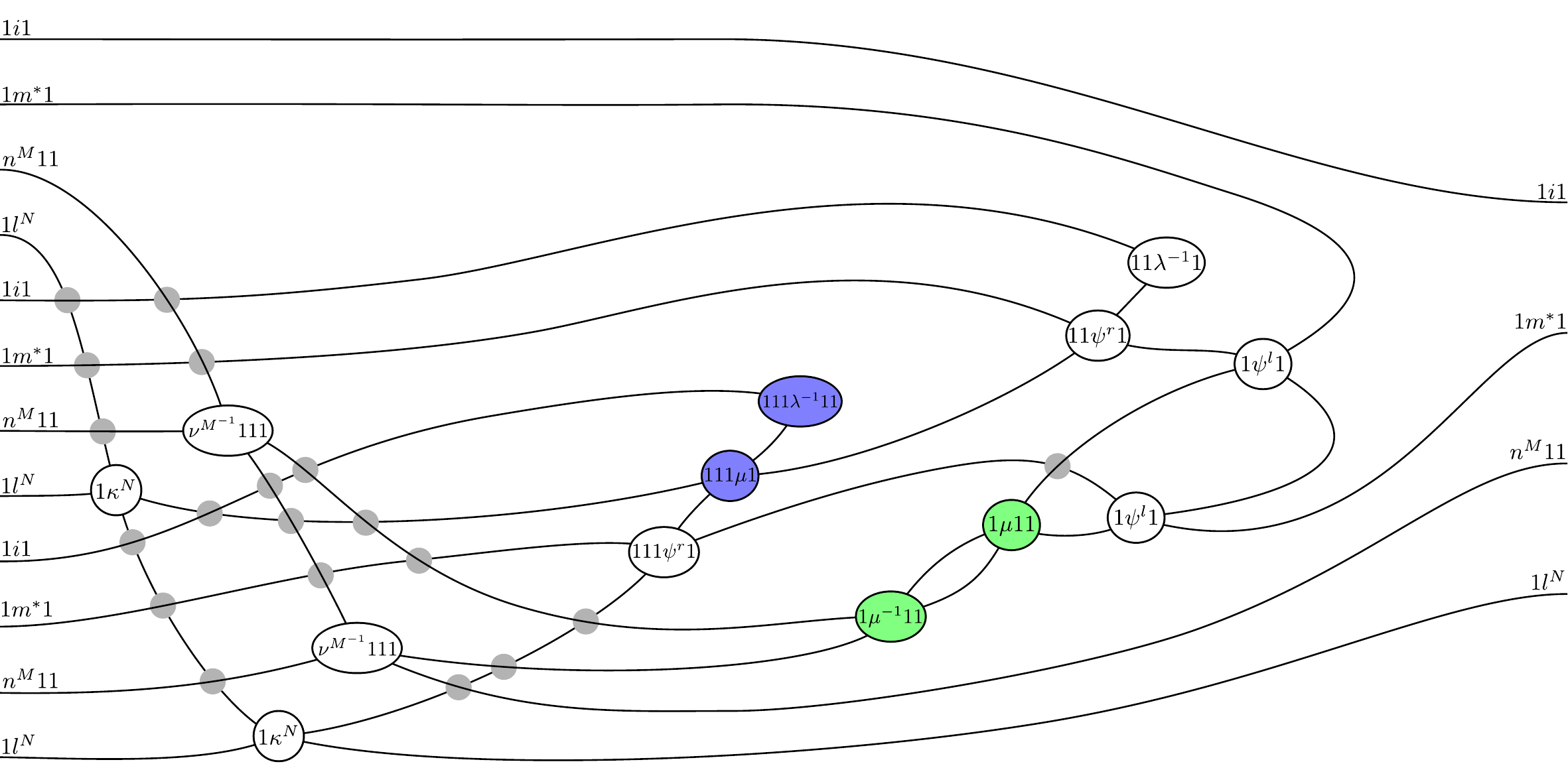}
    \caption{Associativity (Part 7)}
    \label{fig:muassociativity7}
    \end{figure}
    \vfill
\end{landscape}

\begin{landscape}
    \vspace*{\fill}
    \begin{figure}[!hbt]
    \centering
    \includegraphics[width=180mm]{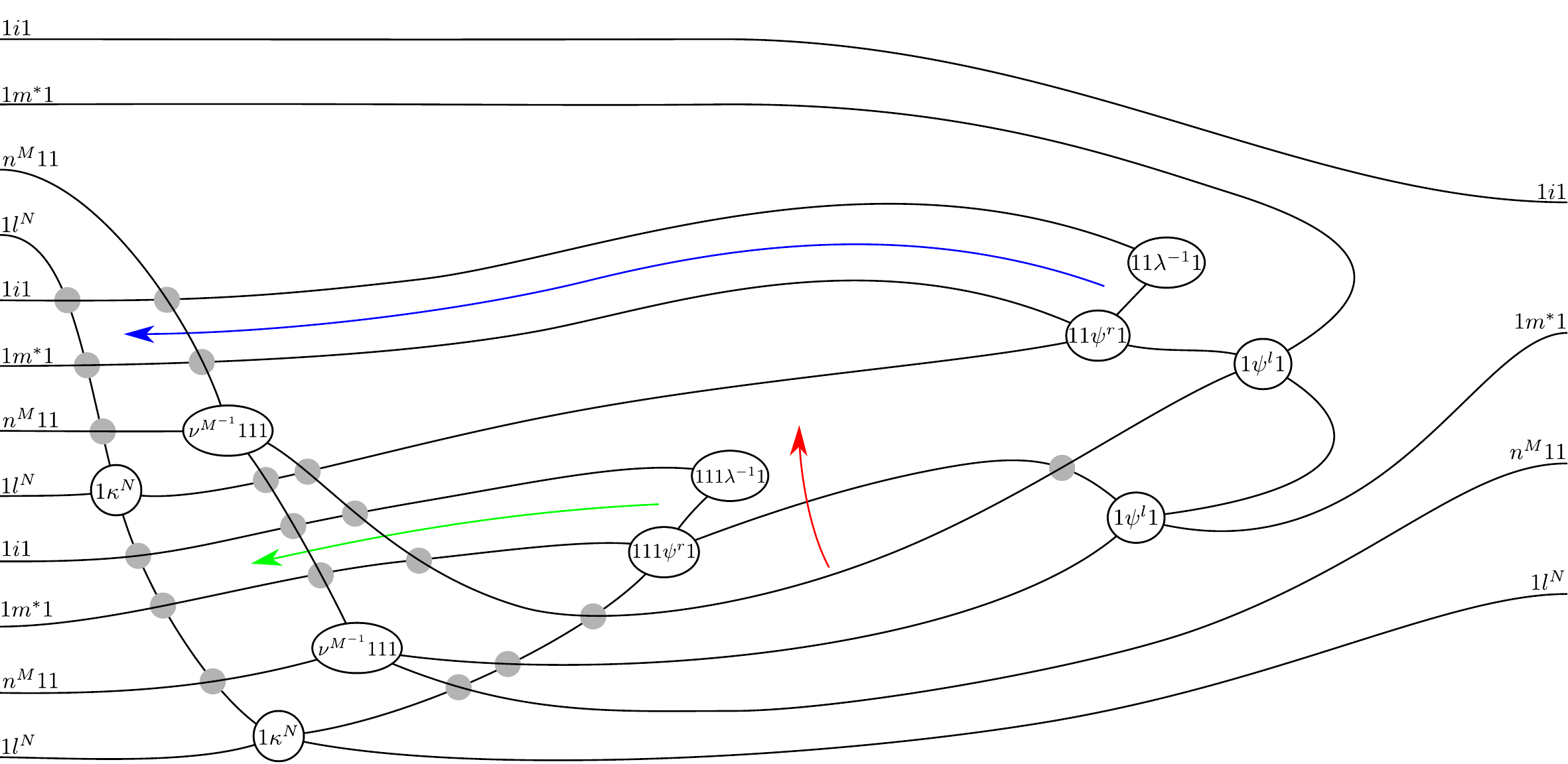}
    \caption{Associativity (Part 8)}
    \label{fig:muassociativity8}
    \end{figure}
    \vfill
\end{landscape}

\begin{landscape}
    \vspace*{\fill}
    \begin{figure}[!hbt]
    \centering
    \includegraphics[width=180mm]{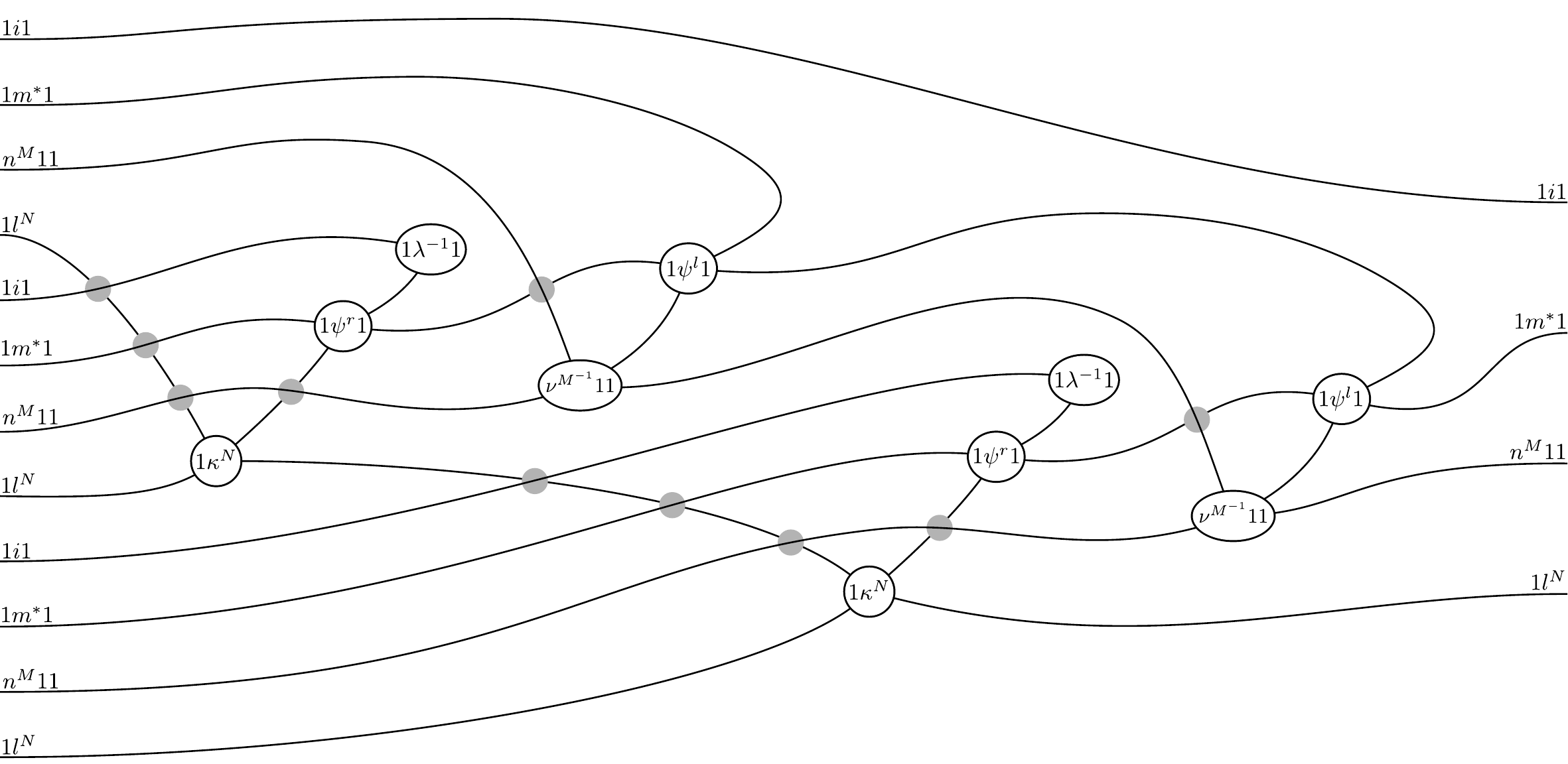}
    \caption{Associativity (Part 9)}
    \label{fig:muassociativity9}
    \end{figure}
    \vfill
\end{landscape}

\FloatBarrier

\begin{landscape}
    \vspace*{\fill}
    \begin{figure}[!hbt]
    \centering
    \includegraphics[width=180mm]{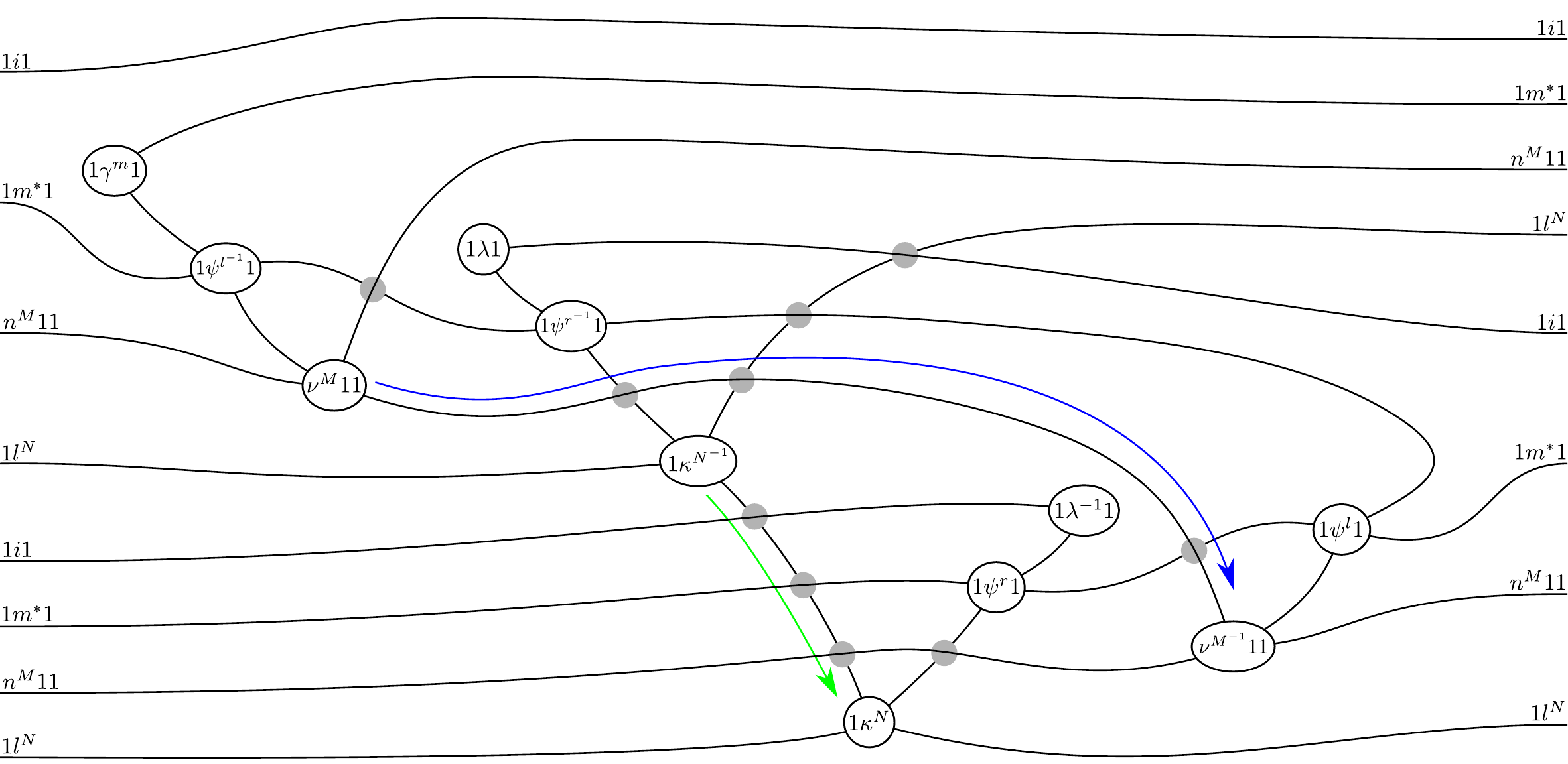}
    \caption{Left Frobenius (Part 1)}
    \label{fig:frobenius1}
    \end{figure}
    \vfill
\end{landscape}

\begin{landscape}
    \vspace*{\fill}
    \begin{figure}[!hbt]
    \centering
    \includegraphics[width=180mm]{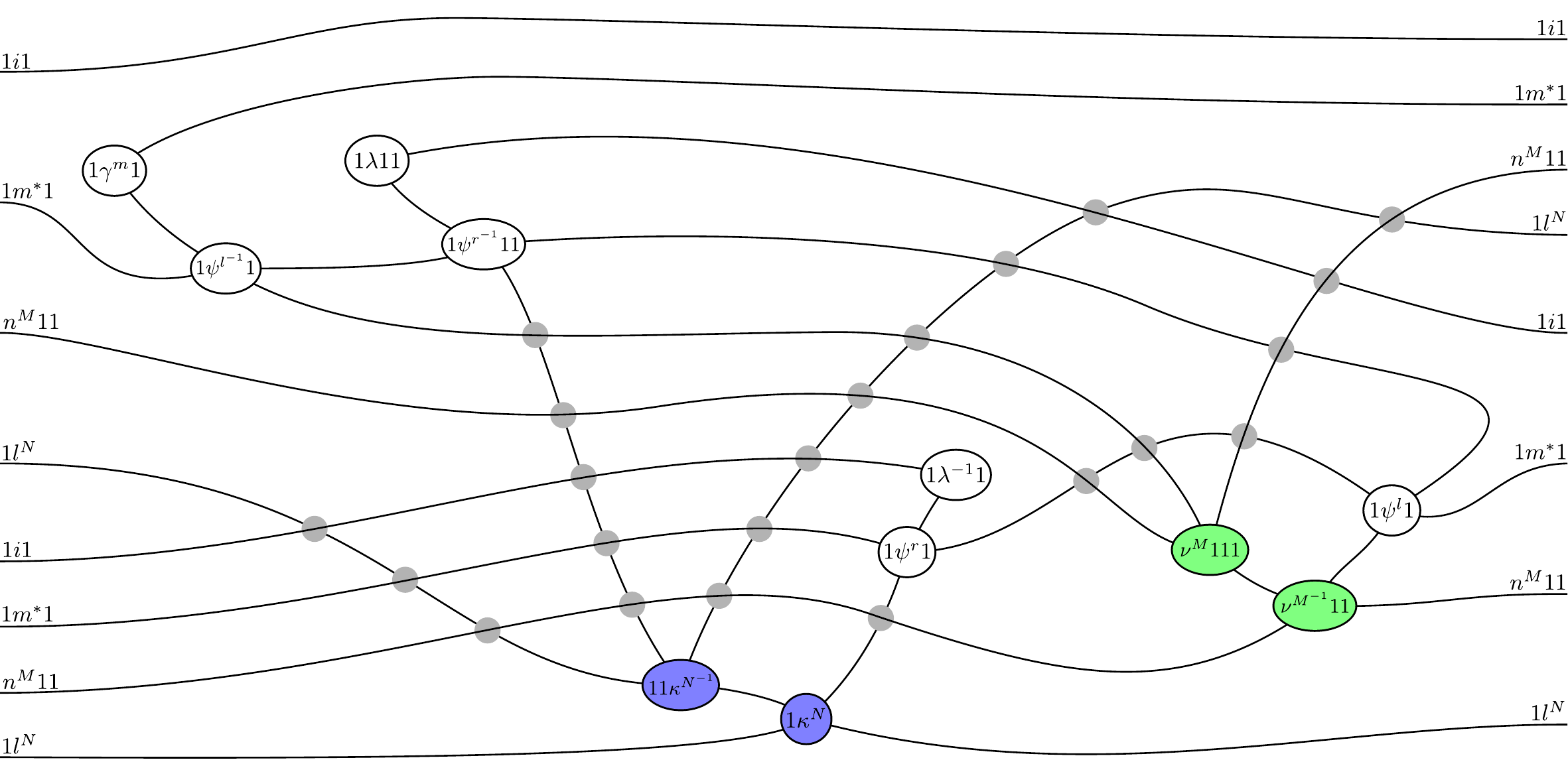}
    \caption{Left Frobenius (Part 2)}
    \label{fig:frobenius2}
    \end{figure}
    \vfill
\end{landscape}

\begin{landscape}
    \vspace*{\fill}
    \begin{figure}[!hbt]
    \centering
    \includegraphics[width=180mm]{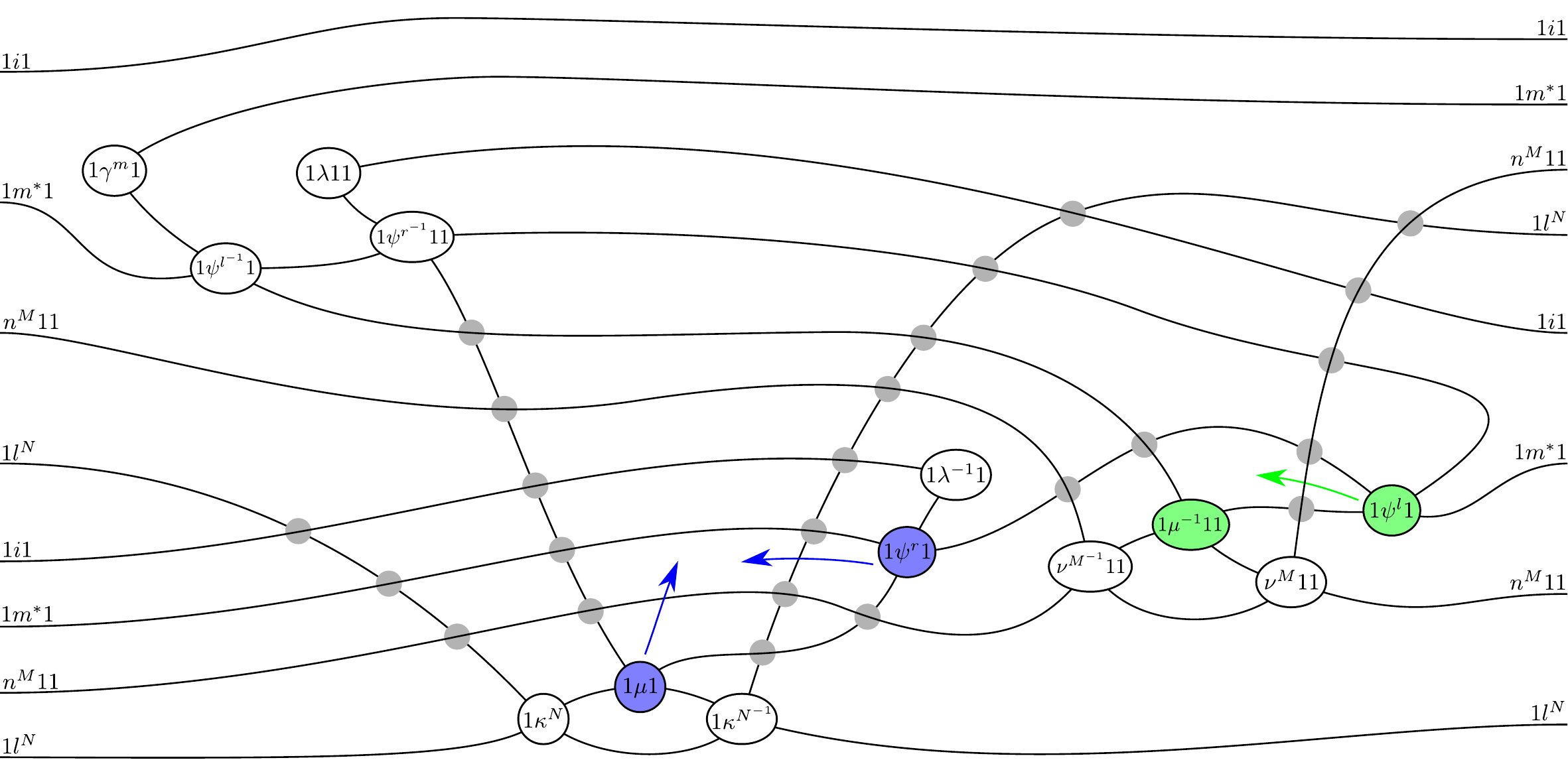}
    \caption{Left Frobenius (Part 3)}
    \label{fig:frobenius3}
    \end{figure}
    \vfill
\end{landscape}

\begin{landscape}
    \vspace*{\fill}
    \begin{figure}[!hbt]
    \centering
    \includegraphics[width=180mm]{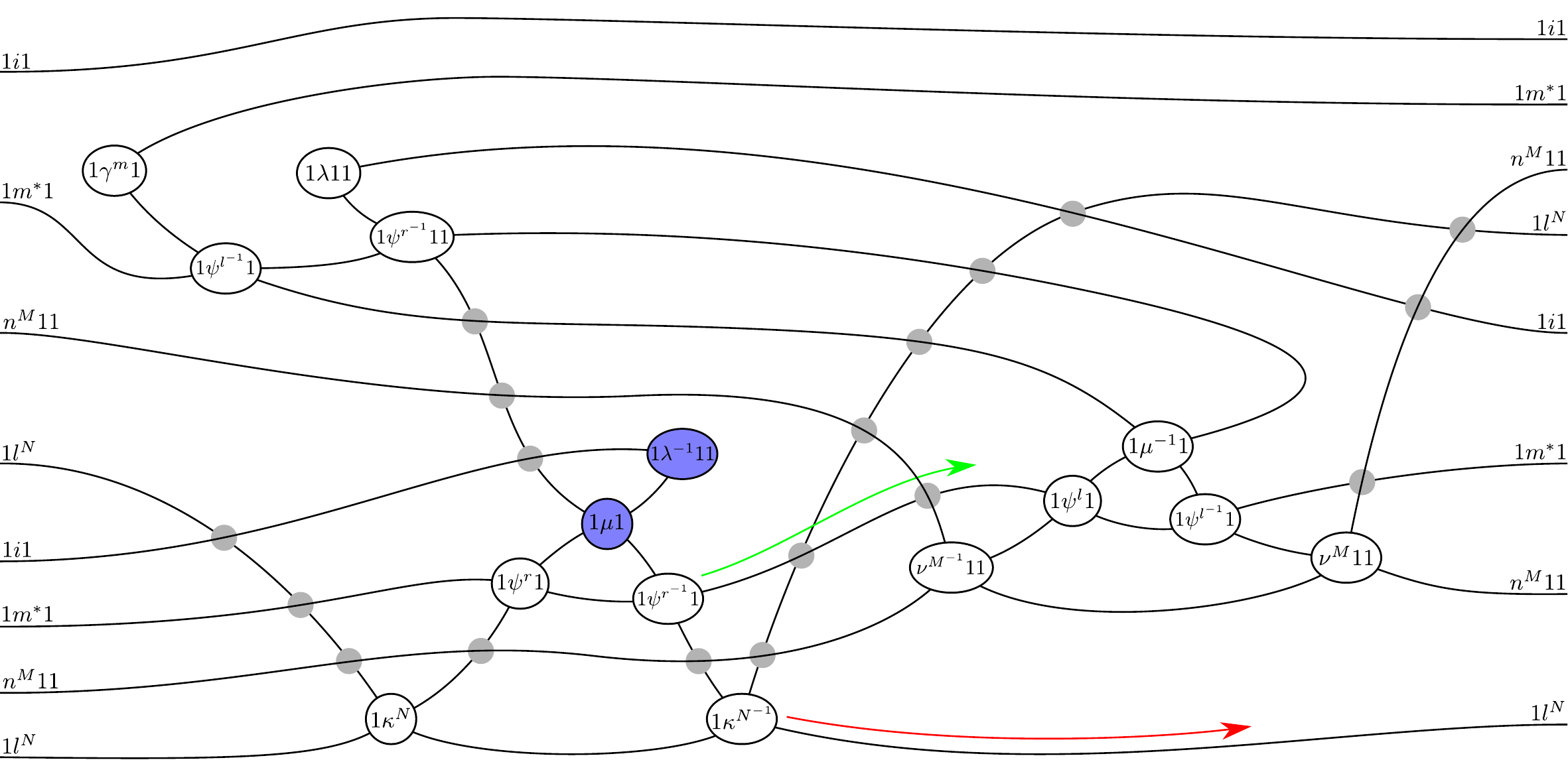}
    \caption{Left Frobenius (Part 4)}
    \label{fig:frobenius4}
    \end{figure}
    \vfill
\end{landscape}

\begin{landscape}
    \vspace*{\fill}
    \begin{figure}[!hbt]
    \centering
    \includegraphics[width=180mm]{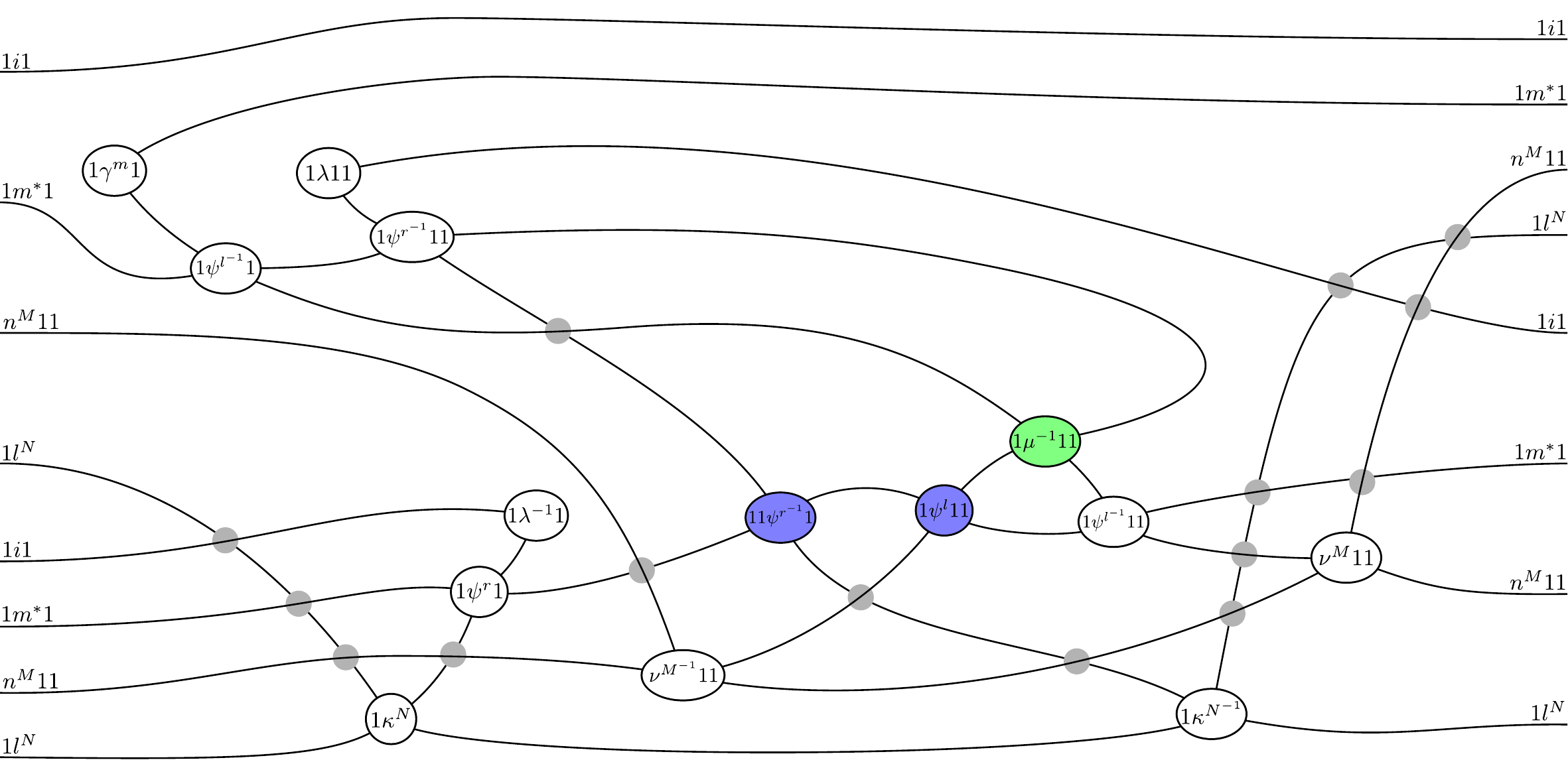}
    \caption{Left Frobenius (Part 5)}
    \label{fig:frobenius5}
    \end{figure}
    \vfill
\end{landscape}

\begin{landscape}
    \vspace*{\fill}
    \begin{figure}[!hbt]
    \centering
    \includegraphics[width=180mm]{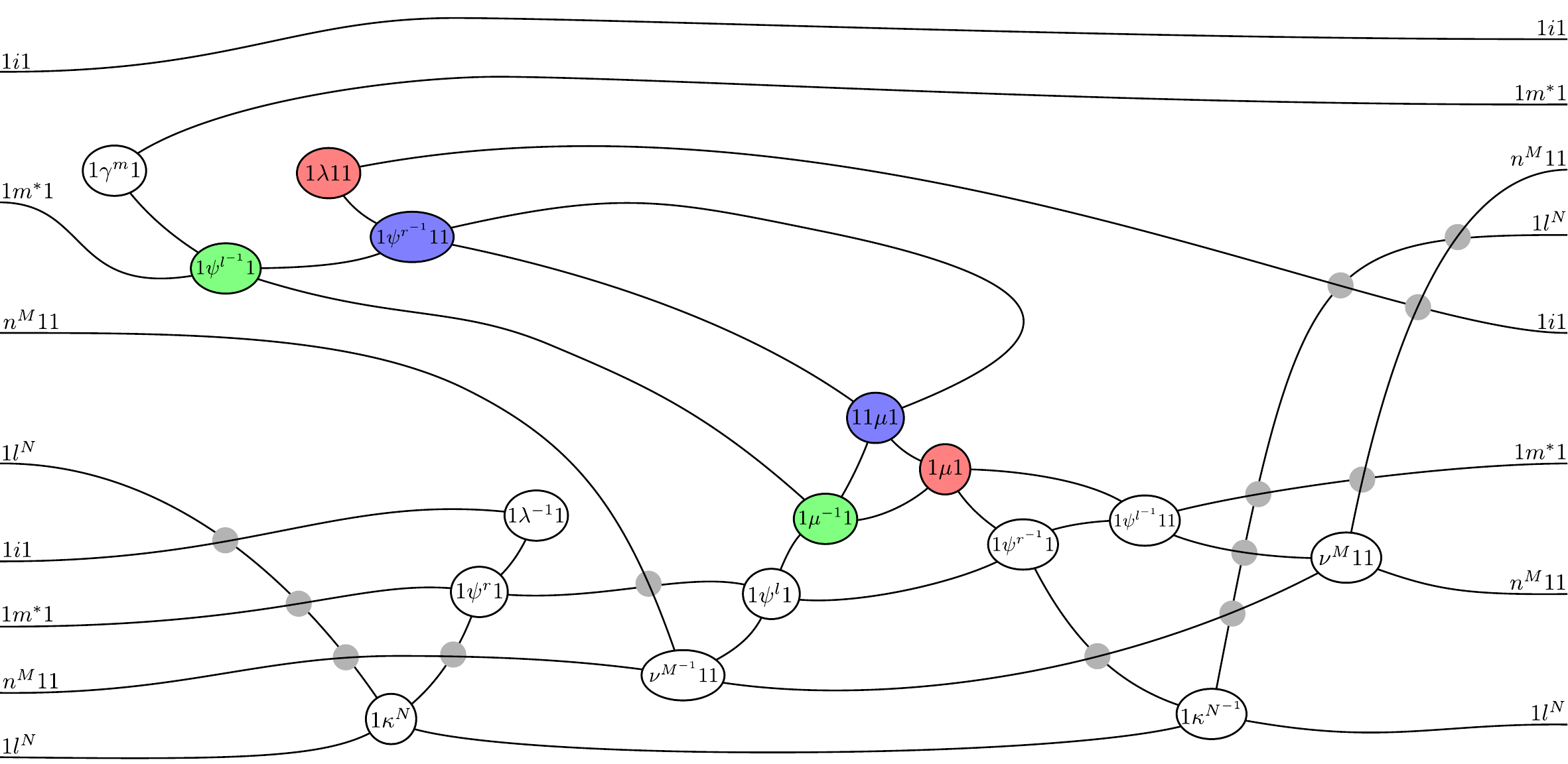}
    \caption{Left Frobenius (Part 6)}
    \label{fig:frobenius6}
    \end{figure}
    \vfill
\end{landscape}

\begin{landscape}
    \vspace*{\fill}
    \begin{figure}[!hbt]
    \centering
    \includegraphics[width=180mm]{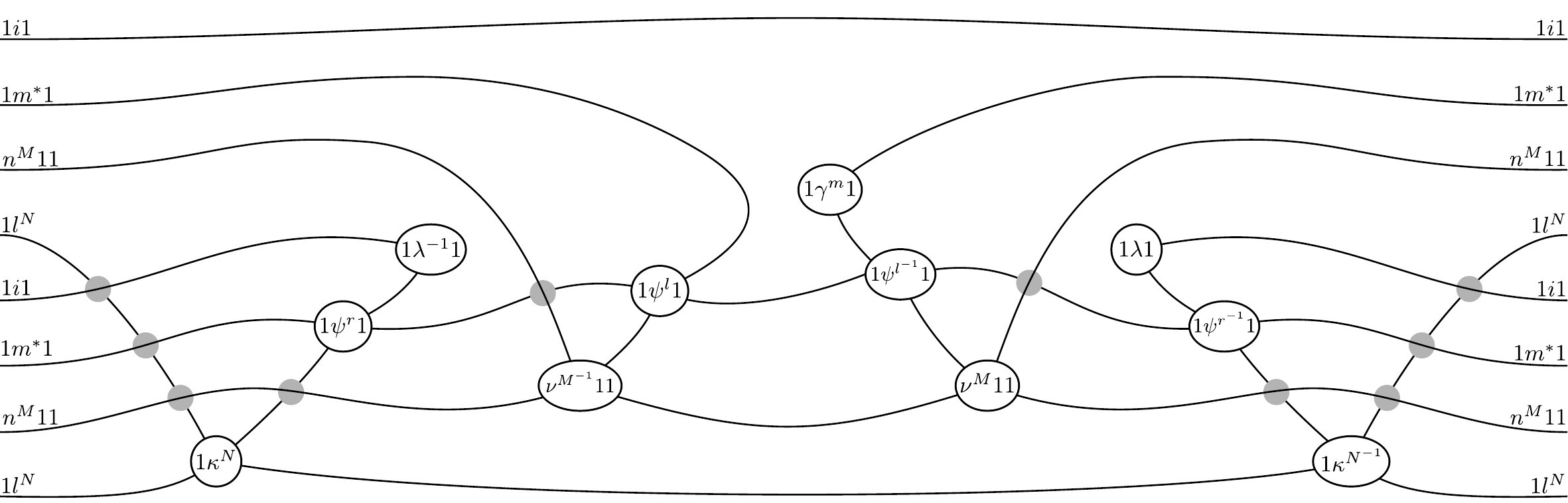}
    \caption{Left Frobenius (Part 7)}
    \label{fig:frobenius7}
    \end{figure}
    \vfill
\end{landscape}

\FloatBarrier